\documentclass{gtart}

\def\ifplaintex{\expandafter\ifx\csname documentclass\endcsname\relax}

\ifplaintex 
\else
\fi

\expandafter\ifx\csname beginpicture\endcsname\relax
\expandafter\ifx\csname documentclass\endcsname\relax
\input pictex \else
\input prepictex \input pictex \input postpictex \fi\fi

\def\gt{{\mathsurround=0pt\it $\cal G\mskip-2mu$eometry \&\ 
$\cal T\!\!$opology}}        

\def\gtp{{\mathsurround=0pt\it $\cal G\mskip-2mu$eometry \&\ 
$\cal T\!\!$opology $\cal P\!$ublications}}  

\def\lognumber#1{\def\thelognumber{#1}}
\def\volumenumber#1{\def\thevolumenumber{#1}}
\def\papernumber#1{\def\thepapernumber{#1}}
\def\volumeyear#1{\def\thevolumeyear{#1}}

\def\pagenumbers#1#2{\def\startpage{#1}\def\finishpage{#2}}
\def\published#1{\def\publishdate{#1}}
\def\proposed#1{\def\theproposer{#1}}
\def\seconded#1{\def\theseconders{#1}}
\def\received#1{\def\receiveddate{#1}}

\def\accepted#1{\def\accepteddate{#1}}

\def\asciiaddress#1{\def\theasciiaddress{#1}}
\def\asciiemail#1{\def\theasciiemail{#1}}
\long\def\asciiabstract#1{\long\def\theasciiabstract{#1}}

\def\shortauthors#1{\def\theshortauthors{#1}}

\let\\\par\let\thelognumber\relax
\let\thevolumenumber\relax\let\thepapernumber\relax
\let\thevolumeyear\relax\let\thesamplenumber\relax\let\startpage\relax
\let\finishpage\relax\let\publishdate\relax\let\receiveddate\relax
\let\reviseddate\relax\let\accepteddate\relax\let\theasciititle\relax
\let\theasciiauthors\relax\let\theasciiaddress\relax
\let\theasciiabstract\relax
\let\theasciiemail\relax\let\theshortauthors\relax\let\theshorttitle\relax

\long\def\maketitlep{   

\count0=\startpage

\gt\hfill      
\beginpicture
\setcoordinatesystem units <0.33truein, 0.33truein> point at 2.2 0.9
\setplotsymbol ({$\cal G$})
\plotsymbolspacing=9truept
\circulararc 315 degrees from 0 1 center at 0 0
\setplotsymbol ({$\cal T$})
\circulararc 315 degrees from 1 -1 center at 1 0
\endpicture
\break
{\small\ifx\thesamplenumber\relax 
Volume \else Sample
\fi\thevolumenumber\ (\thevolumeyear)
\startpage--\finishpage\nl
Published: \publishdate}
\vglue 0.5truein plus 0.4fil minus 0.1truein

{\parskip=0pt\leftskip 0pt plus 1fil\def\\{\par\smallskip}{\ifplaintex\large
\else\Large\fi\bf\thetitle}\par\medskip}   

\vglue 0pt plus 0.1fil 

{\parskip=0pt\leftskip 0pt plus 1fil\def\\{\par}{\sc\theauthors}
\par\medskip}

\vglue 0pt plus 0.1fil 

{\small\parskip=0pt\let\newline\\
{\leftskip 0pt plus 1fil\def\\{\par}{\sl\theaddress}\par}
\expandafter\ifx\theemail\relax    
\relax\else\vglue 5pt plus 0.02fil minus 2pt\def\\{\stdspace{\rm 
and}\stdspace} 
\cl{Email:\stdspace\tt\theemail}\fi
\ifx\theurl\relax                  
\relax\else\vglue 5pt plus 0.02fil minus 2pt\def\\{\stdspace{\rm 
and}\stdspace}
\cl{URL:\stdspace\tt\theurl}\fi\par}

\vglue 7pt plus 0.3fil minus 3pt

{\bf Abstract}
\vglue 5pt plus 0.1fil minus 2pt

\theabstract

\vglue 7pt plus 0.3fil minus 3pt

{\bf AMS Classification numbers}\quad Primary:\quad \theprimaryclass

Secondary:\quad \thesecondaryclass

\vglue 5pt plus 0.3fil minus 2pt

{\bf Keywords}\quad \thekeywords

\vglue 10pt plus 0.5fil minus 5pt

{\small  Proposed: \theproposer\hfill Received: \receiveddate\nl
Seconded: \theseconders\hfill 
\ifx\reviseddate\relax                         
Accepted: \accepteddate                        
\else
Revised: \reviseddate                          
\fi}
\eject
}       

\let\maketitlepage\maketitlep
\let\maketitle\maketitlepage

\font\phead=cmsl9 scaled 950
\font=cmsl9 scaled 1050
\font\pnum=cmbx10 scaled 913
\font\lnum=cmbx10 
\font\pfoot=cmsl9 scaled 950
\font=cmsl9 scaled 1050
\ifplaintex
\headline{\vbox to 0pt{\vskip -4.5mm\line{\small\phead\ifnum
\count0=\startpage ISSN 1364-0380 (on line)
1465-3060 (printed) \hfill {\pnum\folio}\else\ifodd\count0\def\\{ }%
\ifx\theshorttitle\relax\thetitle\else\theshorttitle\fi\hfill{\pnum\folio}
\else\def\\{ and }{\pnum\folio}\hfill\ifx\theshortauthors\relax\theauthors
\else\theshortauthors\fi\fi\fi}\vss}}
\footline{\vbox to 0pt{\vglue 0mm\line{\small\pfoot\ifnum\count0=\startpage
\copyright\ \gtp\hfill\else
\gt, Volume \thevolumenumber\ (\thevolumeyear)\hfill\fi}\vss
}}
\else
\makeatletter
\def\@oddhead{{\small\lhead\ifnum\count0=\startpage ISSN 1364-0380 (on line)
1465-3060 (printed) \hfill {\lnum\number\count0}\else\ifodd\count0
\def\\{ }\ifx\theshorttitle\relax \thetitle \else\theshorttitle\fi\hfill
{\lnum\number\count0}\else\def\\{ and }{\lnum\number\count0}
\hfill\ifx\theshortauthors\relax 
\theauthors\else\theshortauthors\fi\fi\fi}}\def\@evenhead{\@oddhead}
\def\@oddfoot{\small\lfoot\ifnum\count0=\startpage\copyright\ \gtp\hfill\else
\gt, Volume \thevolumenumber\ (\thevolumeyear)\hfill\fi}
\def\@evenfoot{\@oddfoot}
\makeatother
\fi

\newwrite\gtoutfile
\long\gdef\makeheadfile{  
{\def\\{, }\def\s{ }
\immediate\openout\gtoutfile head.xxx
\immediate\write\gtoutfile{Proxy-for: \ifx\theasciiauthors\relax
\theauthors\else\theasciiauthors\fi\s<\ifx\theasciiemail\relax\theemail\else\theasciiemail\fi>}
\immediate\write\gtoutfile{\noexpand\\}
\immediate\write\gtoutfile{Authors: \ifx\theasciiauthors\relax
\theauthors\else\theasciiauthors\fi}
{\def\\{ }\immediate\write\gtoutfile{Title: \ifx\theasciititle\relax
\thetitle\else\theasciititle\fi}}
\immediate\write\gtoutfile{Subj-class: GT or SG or MG etc}
\immediate\write\gtoutfile{MSC-class: \theprimaryclass\ifx\thesecondaryclass\relax\else, \thesecondaryclass\fi}
\immediate\write\gtoutfile{Journal-ref: Geom. Topol. \thevolumenumber
(\thevolumeyear) \startpage-\finishpage}
\immediate\write\gtoutfile{Comments: Published by Geometry and Topology at}
\immediate\write\gtoutfile{\s\s http://www.maths.warwick.ac.uk/gt/GTVol\thevolumenumber/paper\thepapernumber.abs.html}
\immediate\write\gtoutfile{\noexpand\\}
\immediate\write\gtoutfile{}
\ifx\theasciiabstract\relax
\immediate\write\gtoutfile{\theabstract}\else
\immediate\write\gtoutfile{\theasciiabstract}\fi
\immediate\write\gtoutfile{}
\immediate\write\gtoutfile{\noexpand\\}
\immediate\write\gtoutfile{}
\immediate\closeout\gtoutfile}}  

\def\maketitlepage{\maketitlep\makeheadfile}
\let\maketitle\maketitlepage

\lognumber{353}
\volumenumber{7}\papernumber{25}\volumeyear{2003}
\pagenumbers{799}{888}
\received{19 August 2003}
\published{4 December 2003}
\accepted{13 November 2003}
\proposed{Leonid Polterovich}
\seconded{Robion Kirby, Simon Donaldson}

\usepackage{amsmath, amssymb, psfrag, graphicx}
\def\psfraga <#1,#2> #3#4 {%
\psfrag{#3}{\smash{\rlap{\kern #1 \raise #2\hbox{#4}}}}}

\def\abs#1{\mathopen|#1\mathclose|}
\def\Abs#1{\left|#1\right|}

\def\Pc{{\mathcal P}}
\def\Cc{\mathcal{C}}

\def\l{\lambda}
\def\e{\varepsilon}

\def\g{\gamma}
\def\d{\delta}

\newcommand{\diam}{\mathrm{diam}}
\newcommand{\mvee}{\mathop{\vee}}
\newcommand{\og}{\overline{\gamma}}
\newcommand{\ug}{\underline{\gamma}}

\newcommand{\uN}{\underline{N}}
\newcommand{\oN}{\overline{N}}

\newcommand{\circW}{\mathop{W}\limits^{\circ}}

\newcommand{\loc}{\mathrm {loc}}
\newcommand\hook{\mathbin{\hbox{\vrule height .5pt width 3.5pt depth 0pt
\vrule height 6pt width .5pt depth 0pt}}}

\newcommand{\bC}{{\mathbf C}}

\newcommand{\bF}{{\mathbf F}}

\newcommand{\bS}{{\mathbf S}}

\newcommand{\wt}{\widetilde}
\newcommand{\wh}{\widehat}

\newcommand{\Mc}{\mathcal{M}}

\newcommand{\cD}{\mathcal {D}}

\newcommand{\Hc}{\mathcal{H}}

  \newcommand{\dist}{{\mathrm {dist}}}

\newcommand{\Int}{{\mathrm {Int\,}}}

\newcommand{\Z}{{\mathbb Z}}
\newcommand{\ZZ}{{\mathbb Z}}
\newcommand{\RR}{\mathbb {R}}
\newcommand{\R}{\mathbb {R}}
\newcommand{\Q}{\mathbb {Q}\,}
\newcommand{\C}{\mathbb {C}\,}

\newcommand{\const}{{\mathrm const}}

\newcommand{\Tk}{\mathrm{Thick}}
\newcommand{\Tn}{\mathrm{Thin}}
\newcommand{\inj}{\mathrm{injrad}}
\newcommand{\st}{\mathrm{stable}}
\newcommand{\Reeb}{\mathbf{R}}

\newtheorem{theorem}{Theorem}[section]
\newtheorem{corollary}[theorem]{Corollary}
\newtheorem{lemma}[theorem]{Lemma}
\newtheorem{proposition}[theorem]{Proposition}

\newtheorem{example}[theorem]{Example}
\newtheorem{examples}[theorem]{Examples}

\newtheorem{remark}[theorem]{Remark}

\def\R{\mathbb{R}}
\def\Z{\mathbb{Z}}

\def\bR{\mathbb{R}}
\def\bC{\mathbb{C}}

\def\cD{{\cal D}}

\def\od{\overline{d}}
\def\ud{\underline{d}}
\def\of{\overline{f}}

\def\uz{\underline{z}}
\def\oz{\overline{z}}

\def\uZ{\underline{Z}}
\def\oZ{\overline{Z}}
\def\uG{\underline{\Gamma}}
\def\oG{\overline{\Gamma}}

\begin{document}
\title{Compactness results in  Symplectic Field Theory}
\author{F Bourgeois, Y Eliashberg, H Hofer\\K Wysocki, E Zehnder}
\shortauthors{Bourgeois, Eliashberg, Hofer, Wysocki and Zehnder}

\address{Universite Libre de Bruxelles,
B-1050 Bruxelles,
Belgium\\Stanford University, Stanford, CA 94305-2125 
USA\\Courant Institute, New York, NY 10012,
USA\\The University of Melbourne,
Parkville, VIC 3010, Australia\\ETH Zentrum,
CH-8092 Zurich,
Switzerland\\\medskip
\tt{fbourgeo@ulb.ac.be, eliash@math.stanford.edu, hofer@cims.nyu.edu\\K.Wysocki@ms.unimelb.edu.au, eduard.zehnder@math.ethz.ch}}
\asciiaddress{Universite Libre de Bruxelles,
B-1050 Bruxelles,
Belgium\\Stanford University, Stanford, CA 94305-2125 
USA\\Courant Institute, New York, NY 10012,
USA\\The University of Melbourne,
Parkville, VIC 3010, Australia\\ETH Zentrum,
CH-8092 Zurich,
Switzerland}
\asciiemail{fbourgeo@ulb.ac.be, eliash@math.stanford.edu, hofer@cims.nyu.edu,
K.Wysocki@ms.unimelb.edu.au, eduard.zehnder@math.ethz.ch}

\begin{abstract}
This is one in a series of papers devoted to the foundations of
Symplectic Field Theory sketched in \cite{EGH}.  We prove compactness
results for moduli spaces of holomorphic curves arising in Symplectic
Field Theory. The theorems generalize Gromov's compactness theorem in
\cite{Gr} as well as compactness theorems in Floer homology
theory, \cite{F1,F2}, and in contact geometry, \cite{H,HWZ8}.
\end{abstract}

\asciiabstract{
This is one in a series of papers devoted to the foundations of
Symplectic Field Theory sketched in [Y Eliashberg, A Givental and
H Hofer, Introduction to Symplectic Field Theory, Geom. Funct.
Anal. Special Volume, Part II (2000) 560--673]. We prove compactness
results for moduli spaces of holomorphic curves arising in Symplectic
Field Theory. The theorems generalize Gromov's compactness theorem in
[M Gromov, Pseudo-holomorphic curves in symplectic manifolds,
Invent. Math.  82 (1985) 307--347] as well as compactness theorems in
Floer homology theory, [A Floer, The unregularized gradient flow of
the symplectic action, Comm.  Pure Appl. Math.  41 (1988) 775--813 and
Morse theory for Lagrangian intersections, J. Diff.  Geom.  28
(1988) 513--547], and in contact geometry, [H Hofer,
Pseudo-holomorphic curves and Weinstein conjecture in dimension three,
Invent. Math.  114 (1993) 307--347 and H Hofer,
K Wysocki and E Zehnder, Foliations of the Tight Three Sphere,
Annals of Mathematics,  157 (2003) 125--255].}


\maketitlepage

\section{Introduction}
Starting with Gromov's work \cite{Gr}, pseudo-holomorphic curves
(holomorphic curves in short)
became a major tool in Hamiltonian systems, symplectic  and contact geometries.
Holomorphic curves  are smooth maps
$$F\co (S, j)\to (W, J)$$
from Riemann surfaces $(S, j)$ with a complex structure $j$ into almost
  complex manifolds $(W, J)$ with almost complex structure $J$ (ie,
$J^2=-\text{id}$) having the property that their linearized maps $dF$
are complex linear at every point. Hence $F$ satisfies the elliptic
system of first order partial differential equations
$$dF\circ j=J(F)\circ dF.$$ Our aim is to investigate the compactness
properties of the set of holomorphic curves varying not only the maps
but also their domains consisting of punctured and nodal Riemann
surfaces, {as well as (in certain cases) their targets.}  Another
major difference from Gromov's set-up is that the target manifolds of
the holomorphic curves under consideration are not necessarily
compact, and not necessarily having finite geometry at infinity. The
target manifolds are almost complex manifolds with {\it cylindrical
ends} and, in particular, {\it cylindrical almost complex
manifolds}. We also consider the effect of a deformation of the almost
complex structure on the target manifold leading to its {\it
splitting} into manifolds with cylindrical ends. This deformation is
an analogue of the ``stretching of the neck" operation, popular in
gauge theory.  Our analysis introduces the new concepts of stable
holomorphic buildings and makes use of the well known Deligne--Mumford
compactification of Riemann surfaces. The main results are formulated
and proven in Section \ref{sec:compact}. We will see that the moduli
spaces of stable buildings, whose holomorphic maps have uniformly
bounded energies and whose domains have a fixed arithmetic genus and a
fixed number of marked points, are compact metric spaces.  The
compactness results for holomorphic curves proved in this paper cover
a variety of applications, from the original Gromov compactness
theorem for holomorphic curves
\cite{Gr}, to Floer homology theory \cite{F1,F2}, and
to Symplectic Field Theory \cite{EGH}. In fact,
  all  compactness results for holomorphic curves without boundary known to us,
  including the compactness theorems  in  \cite{Ionel-Parker,
Ionel-Parker2}, \cite{Ruan-Li} and  \cite{Li},
   follow from the theorems  we shall prove here.
   Gromov's compactness theorem for closed holomorphic curves asserts
compactness under the condition
   of the  boundedness of the area.  The holomorphic curves
  we consider
  are proper holomorphic maps of Riemann surfaces with punctures into
non-compact manifolds,
   thus  they  usually have infinite area.
  The bound on the area as the condition for compactness is replaced
here by the bound on another quantity, called
  energy. In  applications of Gromov's theorem it  is  crucial to get
an a priori bound on the area. The only
  known case in which  the a priori estimates can be obtained is the
case in which
   the almost complex structure is calibrated or tamed by a symplectic form.
For this reason we decided  to include the appropriate taming
conditions in the  statement
of the theorems. This brings certain simplifications in the proofs
and  still allows to
cover all currently known applications.

\rk{Outline of the paper} We begin in Sections \ref{sec:cylindrical} 
and \ref{sec:ends} with basics about cylindrical almost complex
manifolds and almost complex manifolds with cylindrical ends. We also
describe the process of splitting of a complex manifolds along a real
hypersurface and appropriate symplectic taming conditions which will
enable us to prove the compactness results.  In Section \ref{sec:DM}
we recall standard facts about hyperbolic geometry of Riemann surfaces
and define the Deligne--Mumford compactification of the space of
punctured Riemann surfaces, and its slightly bigger {\it decorated}
version. In Sections
\ref{sec:holom-cylin} and \ref{sec:curves-ends} we define important
notions of {\it contact and symplectic energy} and discuss the
asymptotic behavior and other important analytic facts about
holomorphic curves satisfying appropriate energy bounds. Sections
\ref{buildings-cyl}--\ref{sec:buildings-split} are devoted to the
description of the compactification of the moduli spaces of
holomorphic curves in cylindrical manifolds, manifolds with
cylindrical ends and in the families of almost complex manifolds which
appear in the process of splitting.  Section \ref{sec:compact} is the
central section of the paper, where we prove our main compactness
results.  In Section \ref{sec:other} we formulate some other related
compactness theorems which can be proven using analytic techniques
developed in the paper. Finally, in Appendix A we prove the necessary
asymptotic convergence estimates, and in Appendix B describe metric
structures on the compactified moduli spaces of holomorphic curves.

\rk{Acknowledgements} Parts of this paper were written while some 
of the authors visited the Institute for Advanced Study at Princeton,
Courant Institute for Mathematical Sciences at NYU, New York and
Forschungsinstitut f\"ur Mathematik at the ETH, Zurich. The authors
thank all these institutions for the hospitality. The authors are also
grateful to the American Institute of Mathematics at Palo Alto for
organizing very stimulating workshops and a special program on contact
geometry. Finally, the authors want to thank the referee for his
critical remarks.

F Bourgeois is partially supported by FNRS, the CMI Liftoff program
and EDGE; Y Eliashberg is partially supported by an NSF grant and the
Oswald Veblen Foundation; H Hofer is partially supported by an NSF
grant, a Clay scholarship and the Wolfensohn Foundation; K Wysocki is
partially supported by an Australian Research Council grant; E Zehnder
is partially supported by the TH-project.

\section{Cylindrical almost complex manifolds}
  \label{sec:cylindrical}
  \subsection{Cylindrical almost complex structures}\label{sec:cyl-structures}
An almost complex structure $ J$ on $\R\times V$ is called
{\it cylindrical} if it is invariant under translations
$$(t,x)\mapsto (t+c,x)\;,
t,c\in\R,\;x\in V,$$
and the vector field $\Reeb=J\frac{\partial}{\partial t}$ is  {\it horizontal},
  ie, tangent to levels $t\times V$, $t\in \R$.
  Clearly, any cylindrical structure on $\R\times V$ is determined
  by the CR--structure $(\xi=JTV\cap TV,J_\xi=J|_{\xi})$ and
  the  restriction of the vector field $\Reeb$ to $V=\{0\}\times V$, which
      also will be denoted by $\Reeb$. Furthermore, the distribution
$\xi$ and the vector field $\Reeb$ uniquely determine
       a $1$--form $\l=\l_J$ on $V$ satisfying  $\l(\Reeb)\equiv 1$ and
$\l|_{\xi}=0$.
       The cylindrical structure $J$ is called {\it symmetric}
       if $\l$ is preserved by the flow of $\Reeb$, ie, if the Lie
derivative $L_{\Reeb}\l$ vanishes.
       This is equivalent to the condition $\Reeb\hook d\l=0$ in view
of the Cartan formula
       $L_{\Reeb}\l=\Reeb\hook d\l+d(\l(\Reeb))$.
       It is important to point out that the flow of $\Reeb$ {\it is
not required }
       to preserve $J$ itself.
Here are three important examples of symmetric cylindrical almost
complex structures.
\begin{example}\label{ex:contact}
  {\rm Suppose that $\l$ is a contact form and $\Reeb$ its Reeb vector field.
  Let us recall that this means that $\lambda\wedge (d\l)^{n-1}$,
where $\dim V=2n-1$, is a volume
  form on $V$,  $\Reeb$ generates the kernel of $d\l$ and is
normalized by the condition
  $\l(\Reeb)=1$. Then $L_{\Reeb}\lambda=0$, and hence,
  any $J$ with $\l=\l_J$ is automatically symmetric.
   We refer to this example as to
the {\it contact case}.
It is important to observe that in the contact case the levels
$c\times V,\,c\in\R,$
are } strictly pseudo-convex {\rm being co-oriented as boundaries of
the domains
$(-\infty,c]\times V$. }
   \end{example}
  \begin{example}\label{ex:circle-bundle}
  {\rm
   Let   $\pi\co V\to M$ be a  principal $S^1$--bundle over a closed
manifold $M$ and
       $\Reeb$ the vector field which is the infinitesimal
       generator of the $S^1$--action. Let $\lambda$ be an
        $S^1$--connection form on $V$ such that $\lambda(\Reeb)\equiv 1$.
       Then we have $L_{\Reeb}\lambda=0$  and hence any cylindrical
almost complex structure
       $J$ with $\l_J=\l$  is
       symmetric.        }
       \end{example}
\begin{example}
\label{ex:mapping-torus}
{\rm Suppose that the $1$--form $\l$   is closed. The equality $d\l=0$
trivially implies $\Reeb\hook d\lambda=0$, and hence the
corresponding  $J$ is symmetric. This means, in particular, that
the distribution $\xi$ is integrable and  thus the CR--structure
$(\xi,J_\xi)$ is Levi-flat.  }
\end{example}
The tangent bundle  $T(\R\times V)$ of a cylindrical manifold
canonically splits as  a direct
sum of two complex subbundles, namely  into the bundle $\xi$ and the trivial
$1$--dimensional complex bundle generated by  $\frac{\partial}{\partial t}$.
In this paper  the manifold $V$  is   mostly   assumed to be closed, though
in the last Section \ref{sec:noncompact} we will discuss some
results for   non-compact $V$.
  \subsection{Taming conditions}\label{sec:taming}
  It was first pointed out by M. Gromov (see \cite{Gr}) that
though the local theory of holomorphic curves in an almost complex manifold
is as rich as in the integrable case, the meaningful global theory does
not exist unless there is a symplectic form which  {\it tames}
  the almost
complex structure.
Let us recall that a linear  complex structure $J$ on a vector space
$E$ is called
  {\it tamed} by a  symplectic form
$\omega$ if the form $\omega$ is positive on complex lines. If, in
addition,  one adds the calibrating condition that $\omega$ is
$J$--invariant, then $J$ is   said to be {\it compatible} with $\omega$.
  In the latter case
  $$\omega(X,JY) -i\omega(X,Y),\;\;X,Y\in E,$$  is a Hermitian
  metric on $E$.
Let $(\R\times V,J)$ be a cylindrical almost complex manifold and let
$(\xi,J_\xi)$ and $\Reeb\in TV$ be the   CR--structure and the vector field
which determine   $J$, and $\l=\l_J$  the corresponding $1$--form.
Given a   maximal rank closed $2$--form
$\omega$ on $V$ we say,
to avoid an overused word ``compatible",    that $J$ is {\it adjusted}
  to  $ \omega $ if
    $\omega|_\xi$ is compatible with $J_\xi$ and,  in
    addition,
    $$
    L_{\Reeb}\omega=\Reeb\hook\omega=0.
    $$
    The latter condition means that the vector field $\Reeb$ is Hamiltonian with
    respect to $\omega$.
  Our prime interest in this paper are  {\it symmetric cylindrical
almost complex structures adjusted to a certain
  closed $2$--form}.
  Let us review the adjustment conditions in the Examples
\ref{ex:contact}--\ref{ex:mapping-torus}.
   In the contact case from Example \ref{ex:contact}
   $J$ is  adjusted to $ \omega=d\l $. When referring to the contact
case we will always assume that $d\l$ is chosen as
   the taming form.
   Suppose that
    $\pi\co V\to M$ ,   $\lambda$ and $ R$  are as in Example
\ref{ex:circle-bundle}.
     Suppose that the manifold $M$ is endowed with a compatible
       symplectic form $\underline\omega$ and an almost complex structure
                   $\underline{J}$.
        Set $\omega=\pi^*\underline\omega$ and lift $\underline J$ to
$J_\xi$ on $\xi$
        via the projection $\pi$.
       Then the  symmetric cylindrical almost complex structure
        $J$  on $\R\times V$  determined by the vector field $\Reeb$
and the CR--structure
        $(\xi,J_\xi)$ is  adjusted to the form  $\omega$.
   The adjustment condition in Example \ref{ex:mapping-torus}  requires
the existence of a closed $2$--form  $\omega$ on $V$ whose  restriction to
the leaves of the foliation $\xi$ is symplectic and
is compatible  with $J_\xi$.
An important special case of this construction, which appears in the
Floer homology theory,
is the case in which
the form $\lambda$ has an integral
cohomology class and hence  the manifold $V$ fibers over  the circle
$S^1$ with symplectic fibers. If $(M,\eta)$ is the fiber of this
fibration
then $V$ can be viewed
as  the mapping torus $$V=[0,1]\times M/\{(0,x)\sim(1,f(x)),\,x\in M\}$$
    of a symplectomorphism $f$ of a symplectic manifold $(M,\eta)$.
\subsection{Dynamics of the vector field $\Reeb$}\label{sec:dynamics}
Suppose that a symmetric cylindrical almost complex structure $J$ is
adjusted to a closed form $\omega$.
This implies that the vector field $\Reeb$ is Hamiltonian: its flow
preserves the form $\omega$.
Let us denote by $\Pc=\Pc_J$ the set of periodic trajectories,
counting their multiples, of the vector field $\Reeb$
  restricted to $V=\{0\}\times V$.
Generically, $\Pc$ consists of  only countably many periodic trajectories.
     Moreover, these trajectories   can be assumed to be  {\it non-degenerate}
in the sense  that the linearized Poincar\'e return map $A_\g$
  along any  closed
  trajectory $\g$, including multiples,
has no eigenvalues equal to $1$. We will refer to this generic case
as to the {\it Morse case}.
We will be also dealing in this paper  with a somewhat degenerate,
so-called {\it Morse--Bott
case.} Notice that any smooth family of periodic trajectories from
$\Pc$ has the same period.
Indeed, suppose we are given a map $\Phi\co S^1\times[0,1]\to V$
  satisfying  $\g_\tau=\Phi|_{S^1\times \tau}\in\Pc$ for all $\tau\in[0,1]$.
  Let us denote by $T_\tau $ the period of  $\g_\tau$.
  Then we have
  $$T_1-T_0=\int\limits_{\g_1}\l-\int\limits_{\g_0}\l=\int\limits_{S^1\times[0,1]}\Phi^*d\l=0\,,$$
in view of  $\Reeb\hook d\l=0$ so that  $\Phi^*d\l$ vanishes on
$S^1\times[0,1]$.
We say that a symmetric cylindrical almost complex structure $J$
is of  the  {\it Morse--Bott type} if, for every $T > 0$ the subset
$N_T\subset V$ formed by the closed trajectories from $\Pc$ of period $T$
is a smooth closed submanifold of $V$, such that
the rank of $\omega |_{N_T}$ is locally constant and $T_pN_T = \ker
(d(\varphi_T) - I)_p$, where $\varphi_t\co V\to V$ is the flow generated
by the vector field $\Reeb$.

{\em In this paper we will only consider symmetric cylindrical almost
complex structures $J$ for which either the Morse, or the Morse--Bott
condition is satisfied}.  

Given two homotopic periodic orbits
$\g,\g'\in\Pc$ and a homotopy $\Phi\co S^1\times[0,1] \to V$
connecting them, we define their {\it relative $\omega$--action} by the
formula\footnote{Sometimes, when it will be explicitly said so, we
consider the relative $\omega$--action between the $\Reeb$--orbits when
one, or both of the orbits come with the opposite orientation.}
\begin{equation}\label{eq:rel-action}
\Delta S_\omega(\g,\g';\Phi)=
\int\limits_{S^1\times[0,1]}\Phi^*\omega. 
\end{equation}
  If
  the taming form $\omega$ is exact,  ie, $\omega=d\theta$ for a
one-form $\theta$ on $V$, then
  $$\Delta
S_\omega(\g,\g';\Phi)=\int\limits_\g\theta-\int\limits_{\g'}\theta\,,$$
   so that the relative action is independent of the homotopy $\Phi$.

We will also introduce the {\it $\l$--action}, or simply the {\it
action} of a periodic $\Reeb$--orbit $\g$ by
$$S(\g)=\int\limits_{\g}\l\,.$$ Thus, in the
  contact case in which $\omega=d\lambda$,
  $$\Delta S_\omega(\g,\g';\Phi)=S(\g)-S(\g').$$

\section{Almost complex manifolds with cylindrical ends}\label{sec:ends}
\subsection{Remark about gluing two manifolds along their
boundary}\label{sec:gluing}
  We consider two manifolds $W$ and $W'$ with boundaries and let
$V$ and $V'$ be their boundary components. Given a diffeomorphism
$f\co V\to V'$, the manifold
$$\wt W=W\bigcup\limits_{V\mathop{\sim}\limits^f V'}W' $$ glued along
$V$ and $V'$ is defined as a piecewise
smooth manifold.   If $W$ and $W'$ are oriented, and $f\co V\to V'$
reverses the orientation then $\wt W$ inherits the orientation.
However, to define a smooth structure on $\wt W$ one needs to make
some additional choices (eg, one has to choose embeddings
$I\co (-\e,0]\times V\to W$ and $F\co [0,\e')\times V\to W'$ such that
$I|_{0\times V}$ is the inclusion $V\hookrightarrow W$ and
$F|_{0\times V}$ is the composition $V\mathop{\to}\limits^{f}
V'\hookrightarrow W'$).
Suppose that the manifolds
  $W$ and $W'$ are endowed with almost complex structures $J$ and $J'$. Then
  the tangent bundles
$TV$ and $TV'$  carry {\it CR--structures}, ie, complex subbundles
$$\xi=JTW\cap TW\quad\hbox{ and }\quad     \xi'=JTW'\cap TW'.$$
  Then, in order to define an almost complex structure
on $\wt W=W\bigcup\limits_{V'\mathop{\sim}\limits^f (V)}W' $  the orientation
reversing diffeomorphism
$f$ must preserve these structures. In other words, we should have
$df(\xi)=\xi'$ and
$df\co (\xi,J)\to(\xi',J')$ should be a  homomorphism of complex bundles.
In this case $df$ canonically extends to a complex bundle
homomorphism $$\wt{df}\co
(TW|_V,J)\to (TW'|_{V'},J') ,$$ and thus allows us to define a  $C^1$--smooth
structure on $\wt W$ and a {\it continuous}
almost complex structure $\wt J$ on $\wt W$.  To define a
$C^\infty$--smooth structure $\wt J$ on $\wt W$
we  may need not only  to choose some additional data, as in the case
of a smooth structure,
but also to perturb either $J$ or $J'$ near $V$.
  A particular choice of the perturbation  will  usually be irrelevant
for us, and  thus will  not be specified.

Suppose now that $(W,\Omega)$ and $(W',\Omega')$ are
symplectic manifolds with boundaries  $V$ and $V'$  and introduce
$\omega=\Omega|_V$, $\omega'=\Omega'|_{V'}$.  Suppose that
there exists a diffeomorphism $f\co V\to V'$ which reverses  the
orientations induced on $V$ and $V'$ by the symplectic
orientations of $W$ and $W'$ and which satisfies
$f^*\omega'=\omega$. Notice that for any $1$--form $\lambda$ on $V$
which does not vanish on the (1--dimensional) kernel of the form
$\omega$ we can form for a sufficiently small $\e>0$ a symplectic
manifold $\left((-\e,\e)\times V,\wt\omega+d(t\lambda)\right)$,
where $t\in(-\e,\e)$ and $\wt\omega$ is the pull-back of $\omega$
under the projection $(-\e,\e)\times V\to V$. According to a
version of Darboux' theorem any symplectic manifold containing
a hypersurface $(V,\omega)$,  is symplectomorphic near $V$ to
$\left((-\e,\e)\times V,\wt\omega+d(t\lambda)\right)$ via a
symplectomorphism fixed on $V$. In particular, the identity map
$V\to   0\times V$ extends to a symplectomorphism of a
neighborhood of $V$ in $W$ onto  the lower-half $ (-\e,0]\times V
)$, while the map $f^{-1}\co V'\to 0\times V$ extends to a
symplectomorphism of a neighborhood of $V'$ in $W'$ onto the upper
half $ [0,\e)\times V$. This allows us to glue $W$ and $W'$ into a
smooth symplectic manifold
$W\mathop{\cup}\limits_{V\mathop{\sim}\limits^{f}V'}W'$. The
symplectic structure on this manifold is independent of extra
choices up to symplectomorphisms  which are the identity
maps on $V=V'$ and also  outside of a neighborhood of this
hypersurface. Hence these extra choices will not be
usually specified.
  \subsection{Attaching a cylindrical end}\label{sec:attaching}
  
  Let $(\overline W, J)$ be a compact smooth manifold with boundary and let
  $J'$  a cylindrical almost complex structure on $\R\times V$ where
  $V=\partial\overline W$. If $J$ and $J'$ induce on $V$ the same CR--structure
then, depending on  whether  the orientation of $V$ determined by
$J'|_{\xi'}$ and $\Reeb$, is opposite  or coincides
  with the  orientation of the boundary of  $\overline W$ we can, as
it is described  in Section \ref{sec:gluing} above,
attach to $(\overline W,J)$ {\it the  positive cylindrical end}
  $(E_+=[0,\infty)\times V,J'|_{E_+})$, or the {\it negative} cylindrical end
  $(E_-=(-\infty,0]\times V,J'|_{E_-})$, ie, we
   consider the manifold
   \begin{equation}
   \label{eq:cyl-end}
   W=\overline W\mathop{\cup}\limits_{V=0\times V}E_\pm
   \end{equation}
   with the induced complex structure, still denoted by $J$.
  Alternatively, we say that an almost complex
  manifold $(X,J_X)$
  has a {\it cylindrical end}  (or {\it ends} if $X$ is disconnected
  at infinity) if it is biholomorphically equivalent to a manifold of the form
  (\ref{eq:cyl-end}).
   To get a concrete model of this kind (for a positive end), let us
choose a tubular neighborhood
   $U=[-1,0]\times V$ of $V=\partial \overline W$ in $\overline W$.
Let
$g^\d\co [-\d,\infty)\to[-\d,0)$,  with $0<\d<1$, be a monotone and (non-strictly)
  concave function which coincides with
$$t\mapsto -\frac{\d}2 e^{-t}$$ for $t\in [0,\infty)$ and which is
the identity map near $-\d$.
We define  a family of diffeomorphisms
$G^\d\co W\to\circW=\Int \overline W$  by means of the formula
\begin{equation*}
G^\d(w)=
\begin{cases}
(g^\d(t),x),&w=(t,x)\in U^\d\cup E_+\cr
  w,  &w\in W\setminus(  U^\d\cup E_+)\,,\cr
  \end{cases}
  \end{equation*}
   where $U^\d=[-\d,0]\times V\subset U$.
     The push-forward $(G^\d)^*J$ will be denoted by   $J^\d$.
     Thus $(\circW,J^\d)$ can be viewed as another model of an
     almost complex manifold with cylindrical end.
We say that $(X,J_X)$ has {\it asymptotically cylindrical}
  positive resp.  negative end if there exists a diffeomorphism
  $$ f\co W=\overline W\mathop{\cup}\limits_{V=0\times V}E_\pm\to X$$
  such that the families of mappings $f^s\co E_{\pm}\to X$, for $s\geq
0$, defined on the ends by the formulae $f^s(t, x)=f(t\pm s, x)$ have
the following properties,
\begin{itemize}
  \item $J^s\co =(f^s)^*J_X\longrightarrow J$ in $C^{\infty}_{\text{loc}}$.
  \item $J^s(\frac{\partial}{\partial t})=\Reeb$ for all $s\geq 0$.
\end{itemize}
  We say that $W$ has a {\it symmetric}  cylindrical end if     the
almost complex structure $J|_{E_\pm}$ is symmetric, ie,
$L_{\Reeb}\lambda=0$ where $\Reeb$ and $\lambda$ are the vector field
and the $1$--form
on $V$ introduced in  the definition of a cylindrical almost complex structure.
The notion of an almost complex manifold   with symmetric
asymptotically cylindrical end has an obvious meaning.
  \begin{examples}\label{ex:cyl-ends}
  {\rm
  \begin{itemize}
  \item[(1)] $\C^n$ has  a symmetric (and even contact type) cylindrical end.
  \item[(2)] For any complex manifold $W$ the punctured manifold
$W\setminus \{p\}$ has a
  cylindrical symmetric end.
  \item[(3)] More generally, let $(X,J)$ be an almost
    complex manifold,
   and $Y\subset X$   an almost  complex submanifold of  any real
codimension $2k$.
Then the manifold $(X\setminus Y,J|_{X\setminus Y})$ has a
symmetric asymptotically cylindrical end  corresponding  to the
cylindrical manifold $E=\R\times V$, where $V$ is the sphere
bundle over $Y$ associated with the complex normal bundle to $Y$
in $X$, and the vector field $\Reeb$ is tangent to the fibers. In
particular, there is another almost complex structure $J'$ on $X$
which is $C^1$--close to $J$ and coincides with $J$ on $T(X)|_Y$
such that the almost complex manifold $(X\setminus
Y,J'|_{X\setminus Y})$ has a symmetric cylindrical end. 
\item[(4)]
Let $(S\setminus Z,j)$ be a closed Riemann surface with punctures,
and $(M,J)$ be any almost complex manifold. Then
$\left((S\setminus Z)\times M, j\oplus J\right)$ has cylindrical
ends.
\end{itemize}
}
  \end{examples}
\subsection{Taming conditions for almost complex
   manifolds with cylindrical ends}\label{sec:taming-ends}
    Let  $(W=\overline W\mathop{\cup}\limits_{V} E  ,J)$
      be an almost complex  manifold with a cylindrical end $E\subset\R\times V$
        where $V=\partial\overline W$. Let
      $\Reeb$ be the associated vector field
      tangent to $ V $, and let
     $(\xi,J_\xi)$  be the CR--structure,  and $\lambda$ the $1$--form
       on $V$ associated
      with  the cylindrical structure $J|_E$.
     Suppose that the manifold $\overline W$ admits a symplectic form $\omega$.
   We say that the almost complex structure
  $J$  is {\it adjusted} to $\omega$
if
\begin{itemize}
\item the symplectic form $\omega$ is compatible with $J|_{\overline W}$,
\item $J|_{E}$ is adjusted to $\omega$ is the sense of the definition in
Section \ref{sec:taming}.
\end{itemize}
Notice that   the distinction between positive and negative ends
depends on the $1$--form $\lambda$.
In the contact case, ie, when the $1$--form $\l$ associated with  a
cylindrical end
$E$ is a contact form and $\omega=d\l$, this sign is beyond  our control;
the boundary component $V_0$ is positive  iff there exists an
outgoing  vector field $X$ transversal
to  $V_0$ which dilates  the symplectic form $\omega$,
ie, $L_X\omega=\omega$. On the other hand, in other cases (see, for
instance Examples \ref{ex:circle-bundle} and \ref{ex:mapping-torus})
  the sign of $\l$, and hence the sign
  of the corresponding end $E$ can be changed at our will. We  will
use  the notation
  $E_-$ respectively $E_+$ for the union of  the negative ends respectively
the positive ends.
    Similarly, the set $\Pc$ of periodic orbits of the  vector  field
$\Reeb$ also  splits into the disjoint
   sets
   $\Pc_-$ and $\Pc_+$ of orbits on $V_-$ and $V_+$.
Consider now an  almost complex structure $ J$ on $W=\overline
W\mathop{\cup}\limits_{V}E$
  with an {\it asymptotically} cylindrical end. The structure  $ J$ is
called  {\it adjusted} to
  $\omega$ if
  \begin{itemize}
\item  $\omega$ is compatible with $ J|_{\overline W}$,
\item  $\omega|_\xi$ is compatible with $J_t|_{\xi}$ for each $t\geq
0$ (or $t\leq 0$)
where $J_t$
is the pull-back of $J$ under the inclusion map
$V=V\times t \hookrightarrow E$.
\end{itemize}
         Note that in the situation of the last definition, the $2$--form
         $\omega$ can always be extended
         to $W=\overline W\mathop{\cup}\limits_{V}E$  as a symplectic
form $\wt \omega$
         taming $J$.
\subsection{Splitting}  \label{sec:split}
  The following splitting construction is an important source of
   manifolds with cylindrical ends.
   Let $W$ be a closed almost complex manifold, or a manifold with cylindrical
   ends, and $V\subset W$ a co-orientable compact real hypersurface.
Let $\xi=TV\cap JTV$
   be the CR--structure induced on $V$, ie, the distribution of maximal
   complex tangent
    subspaces  of $TV$.
    Let us cut $W$ open along $V$. The boundary of the newly created
    manifold $\circW$ (which is disconnected if $V$ divides $W$)
    consists of two  copies $V^{'}$, $V^{''}$
    of $V$.  Choose any vector field
    $\Reeb\in TV$ transverse to $\xi$ and  attach to $\circW$
    the ends $E_+=[0,\infty)\times V$
    and $E_-=(-\infty,0]\times V$  with   the unique cylindrical
almost complex structure
     determined by the CR--structure $(\xi, J_\xi)$ and the vector
field $\Reeb$. As it was pointed out in Section
     \ref{sec:gluing}
      the resulting  manifold
    $$\wt W=E_-\mathop{\cup}
    \limits_{V\times 0=V^{'}}
    \circW\mathop{\cup}\limits_{V^{''}=V\times 0} E_+\,.$$
     gets a canonical $C^1$--smooth structure and a continuous almost
complex structure
     $\wt J$. In order to make $\wt J$ smooth in the $C^{\infty}$
sense we may need to perturb $J$ on $\circW$ near $V$.
     A  specific choice of the perturbation will be irrelevant for our
purposes, and thus not specified.
         The manifold  $(\wt W,\wt J)$ has
      cylindrical ends and  is diffeomorphic to $W\setminus V$.
   Note that the above splitting construction can be viewed as the result of
   a ``stretching of the neck". Indeed one can consider a family of manifolds
   $$W^\tau=\circW\mathop{\cup}\limits_{V^{''}=(-\tau)\times
V,\,\tau\times V=V'}
   [-\tau,\tau]\times V \; , \:  \tau\in[0,\infty),$$
   with the almost complex structure $J^\tau$ on
    the insert  $[-\tau,\tau]\times V$  which is  uniquely determined
by the same
    translational invariance condition. Then we have
    $$(W^\tau,J^\tau)\xrightarrow[\tau \to \infty] {} (\wt W,
    \wt J)\,.$$
\begin{figure}[ht!]\small
\centering \psfrag{v}{$V$} \psfrag{w}{$W$} \psfrag{vt}{$[-\tau,
\tau]\times V$} \psfrag{wt}{$W^{\tau}$}
\psfraga <-10pt, -7pt> {em}{$E_{-}=(-\infty, 0]\times V$}
\psfraga <-10pt, 0pt> {ep}{$E_{+}=[0,\infty)\times V$} \psfrag{wi}{$W^{\infty}$}
\psfrag{wo}{$\circW$}
\includegraphics[width=0.7\textwidth]{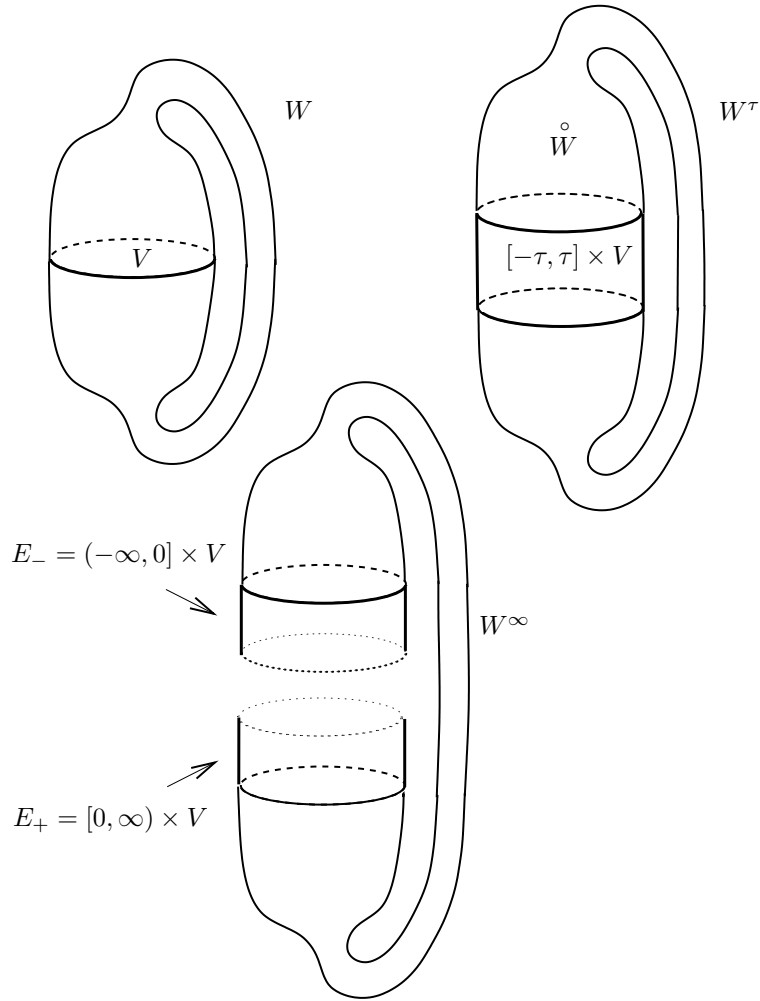}
\caption{Splitting}\label{fig1}
\end{figure}
Suppose that  $\omega$ is a symplectic form on $W$ compatible with $J$,
the vector field $\Reeb$ on the hypersurface $V$  generates the kernel of
$\omega|_V$, and the cylindrical structure on $\R\times V$ defined by
$J$ and $\Reeb$ is symmetric.
We say that in this case the splitting data are adjusted to the
symplectic form $\omega$.
This is, in particular,  the case when $V$ is a contact type
hypersurface in $W$, ie, a hypersurface
which admits in its neighborhood a transversal conformally symplectic
vector field.
In the adjusted case the splitting construction gives a manifold with
symmmetric cylindrical ends which is adjusted to $\omega$.
Note that the hypersurface $V$ is not assumed to divide $W$.
On the other hand, $V$ is allowed to be disconnected and to split
$W$ into several connected components.
Finally, the above splitting construction immediately generalizes to the
case  when the manifold $W$ itself has cylindrical ends.
  \section{Deligne--Mumford compactness revisited}\label{sec:DM}
  \subsection{Smooth stable Riemann surfaces}\label{sec:Riemann-smooth}
Let $\bS=(S,j,M)$ be a compact connected Riemann surface(without
boundary) with a
set $M$ of numbered disjoint marked points.
  Two Riemann surfaces $\bS=(S,j, M)$ and $\bS'=(S',j', M')$ are called
equivalent if there exists a diffeomorphism  $\varphi\co S\to S'$
such that $\varphi_*j=j'$ and $\varphi(M)=M'$ where we assume that
$\varphi$  preserves the ordering of the sets $M$ and $M'$.
  Let $\mu$
be the cardinality of ${ M}$, and
  $g$  the genus of $S$.
  The surface is called {\it stable} if
   \begin{equation}\label{eq:stability}
  2g+\mu\geq 3\,.
  \end{equation}
The stability condition is equivalent to the requirement that
  the group of conformal automorphisms of $\bS$, ie, biholomorphic maps
preserving the marked points, is finite.   The pair $(g,\mu)$ is
called the {\it signature} of the Riemann surface $\bS$.
  Given a  stable surface $\bS=(S,j, M)$, the
Uniformization Theorem asserts  the existence of a unique
   complete  hyperbolic metric of constant curvature
$-1$ of finite volume,  in the given conformal class $j$
on $\dot{S}=S\setminus{  M}$.   We will denote this metric by
  $h^{\bS}=h^{j,M}$.  Each puncture corresponds to  a cusp of the hyperbolic
metric $h^{j,M}$. In what follows we will always assume for a given
stable Riemann surface  $(S,j,M)$ that the punctured surface
$\dot S=S\setminus M$ is endowed with the uniformizing hyperbolic
metric $h^{j,M}$.
   Thus the {\it moduli space} $\Mc_{g,\mu}$
of  the stable Riemann surfaces of signature $(g,\mu)$ can be viewed
equivalently
  as the moduli space of  (equivalence classes of)
  hyperbolic metrics of finite volume and  of constant curvature
$-1$ on  the fixed surface $\dot{S}=S\setminus M$.
\subsection{Thick--thin decomposition}\label{sec:thick-thin}
Fix  an  $\e>0$. Given a stable Riemann surface $\bS=(S,j,M)$ we
denote by $\Tk_\e(\bS)$ and $\Tn_\e(\bS)$ its $\e$--thick and
$\e$--thin parts, ie,
\begin{equation}
\begin{split}
\Tk_\e(\bS)=&\{x\in \dot S|\;\rho (x)\geq\e\}\,\cr
\Tn_\e(\bS)=&\overline{\{x\in \dot S|\;\rho(x)<\e\}}\, ,\cr
\end{split}
\end{equation}
  where $\rho(x)$ denotes  the injectivity radius   of the metric
$h^{j,M}$ at the point $x\in \dot S$.
  It is a remarkable fact of  the
hyperbolic geometry that there exists a universal constant
$\e_0=\sinh^{-1} 1=\ln(1+\sqrt{2})$ such that for any
$\e<\e_0$ each component $C$ of $\Tn_\e(\bS)$ is conformally
equivalent either to a finite cylinder $[-L,L]\times S^1$ if the
component $C$ is not adjacent to a puncture,
  or to the punctured disc $D^2\setminus
0\cong[0,\infty)\times S^1$ otherwise, see for example \cite{Hummel}. Each
{\it compact}
component $C$ of the thin part contains a unique closed  geodesic
of length equal $2\,\rho(C)$, which will be denoted by $\Gamma_C$.
Here $\rho (C)=\inf_{x\in C}\rho (x)$.
When considering $\e$--thick--thin decompositions we will always
assume that $\e$ is chosen smaller than $\e_0$.

\begin{figure}[ht!]\small
\centering
\psfraga <40pt,0pt> {a}{$(a)$} 
\psfraga <40pt,0pt> {b}{$(b)$}
\psfraga <-10pt,3pt> {g}{$\Gamma_{C}$}
\includegraphics[width=0.8\textwidth]{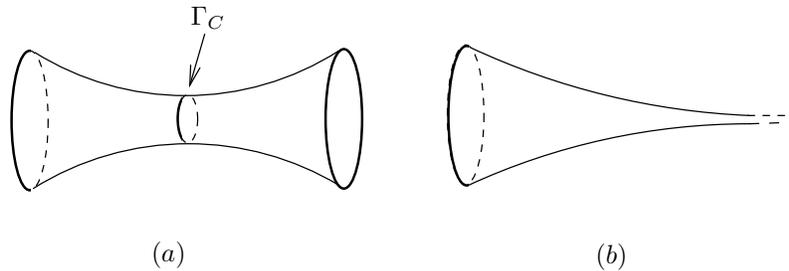}
\caption{ Thin parts non-adjacent   $(a)$ and adjacent $(b)$ to a
puncture}\label{fig2}
\end{figure}

\subsection{Compactification of a punctured Riemann
surface}\label{sec:blow-up} For each marked point $z\in M$ of a
   Riemann surface $\bS=(S,j,M)$ we define the surface $S^{z}$ with
   boundary as an {\it oriented blow-up} of $S$ at the point $z$.
   Thus, $S^{z}$ is the circle compactification of $S\setminus \{z\}$
   and it is a compact surface bounded by the circle
   {$\Gamma_z=(T_zS\setminus 0)/\R_+^{\ast}$, where
   $\R_+^{\ast}=(0,\infty)$.}  The conformal structure $j$ defines an
   action of the circle $ S^1=\R/\Z$ on $\Gamma_z$, and hence allows
   us to canonically metrize $\Gamma_z$.  The canonical projection
   $\pi\co S^{z}\to S$ sends the circle $\Gamma_z$ to the point $z$
   and maps $\Int\, S^{z}$ diffeomorphically to $S\setminus \{z\}$.
   Similarly, given a finite set $M=\{z_1,\dots, z_k\}$ of punctures
   we define a blow-up surface $S^{M}$ having $k$ boundary components
   $\Gamma_1,\dots,\Gamma_k$. It comes with the projection $\pi\co
   S^{M}\to S$ which collapses the boundary circles
   $\Gamma_1,\dots,\Gamma_k$ to the points $z_1,\dots,z_k$ and maps
   $\Int\, S^{M}$ diffeomorphically to $\dot{S}=S\setminus M$.  A
   Riemann surface $\bS$ is called {\it $\e$--thick} if $\Tn_\e(\bS)$
   consists only of non-compact (ie, adjacent to punctures)
   components.  It is not difficult to see that the subspace
   $\Mc^{\e}_{g,\mu}\subset \Mc_{g,\mu}$ of moduli of $\e$--thick
   Riemann surfaces of signature $(g,\mu)$ is compact with respect to
   its natural topology. However, to compactify the moduli space of
   all hyperbolic metrics on $S\setminus M$ one has to add {\it
   degenerate} metrics, or metrics with interior cusps, if the length
   of the closed geodesics $\Gamma_C$ in one or several components of
   the thin parts converge to $0$.
\subsection{Stable nodal Riemann surfaces}\label{sec:Riemann-stable}
We introduce a notion of a {\it
nodal}
  Riemann  surface.
Suppose we are given  a possibly disconnected  Riemann surface
$\bS=(S,j,M,D)$ whose set of marked points   is presented as a
disjoint union of
  sets $M$ and $D$, where the cardinality of the set $D$  is even.
   The marked points from $D$, which are   called {\it special}, are
  organized in pairs,
  $D=\{\od_1,\ud_1,\od_2,\ud_2,\dots,\od_k,\ud_k\}$.
The {\it nodal Riemann surface} is the equivalence class
  of surfaces $(S,j,M,D)$ under  the additional equivalence relations which make
  \begin{itemize}
  \item each pair $(\od_i,\ud_i)$, for
  $i=1,\dots,k$, and
  \item the set of all special pairs
  $\{(\od_1,\ud_1),(\od_2,\ud_2),\dots,(\od_k,\ud_k)\}$
  \end{itemize}
  unordered.
  For notational convenience  we will still denote  the nodal curve by
  $\bS=(S,j,M,D)$, but one should remember that the numeration of
pairs of of points in $D$, and the
  ordering of each pair  is not part of the structure.
  The nodal curve is called {\it stable} if the stability condition
  (\ref{eq:stability}) is
         satisfied for each component of the surface
  $S$ marked by the points from  $M\cup D$.
  With a  nodal surface $\bS$  we can associate the following singular
  surface with double points,
  \begin{equation}
  \label{eq:double}
  \widehat S_{D}=S/\{\od_i\sim\ud_i\,;i=1,\dots,k\}\,.
  \end{equation}
  We shall call the identified points  $\od_i\sim\ud_i$ a {\em node}.
The nodal surface $\bS$ is called {\it connected} if  the singular
surface
  $\widehat S_{D}$ is connected.
If the nodal surface $\bS$ is connected then its {\it arithmetic
genus} $g$
  (compatible with the definition of the deformation $S^{D,r}$ below)
  is defined as
  \begin{equation}\label{eq:arithm-genus}
  g=\frac12\#D-b_0+\sum\limits_1^{b_0}g_i+1\,,
  \end{equation}
   where
  $\#D=2k$ is the cardinality of $D$, and $b_0$ is the number of connected
  components of the surface
  $S$, and $\sum_1^{b_0}g_i$ is the sum of  the genera of the
connected components
  of $S$.  The {\it signature} of a nodal curve $\bS=(S,j,M,D)$ is the
pair $(g,\mu)$ where
  $g$ is the arithmetic genus and   $\mu=\#M$.
  A stable
  nodal Riemann surface $\bS=(S,j,M,D)$ is called {\it decorated} if
   for each special pair there is chosen
  an orientation reversing orthogonal map
  \begin{equation}\label{eq:decorat}
  r_i\co \oG_i=\left(T_{\od_i}(S)\setminus 0\right)/\R_+^{\ast}\to\uG_i
  =\left(T_{\ud_i}(S)\setminus 0\right)/\R_+^{\ast}\;.
  \end{equation}
  \begin{figure}[ht!]\small
\centering \psfrag{ogi}{$\oG_i$} \psfrag{ugi}{$\uG_i$}
\psfrag{ri}{$r_i$}
\includegraphics[width=0.7\textwidth]{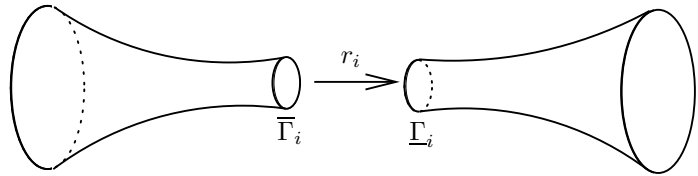}
\caption{The decoration}\label{fig3}
\end{figure}
Orthogonal orientation reversing requires $r(
e^{i\vartheta}z)=e^{-i\vartheta}r(z)$ for all
$z\in  \oG_i$.
We will also consider {\it partially decorated } surfaces if  the maps $r_i$
  are given only for a certain subset $D'$ of special double points.
  An equivalence of decorated
  nodal surfaces must   respect the decorating maps $r_i$.
  The moduli spaces of stable  connected Riemann
  nodal surfaces, and decorated stable connected nodal surfaces  of
signature $(g,\mu)$
  will be denoted by
  $\overline\Mc_{g,\mu}$ and  $\overline\Mc^{\$}_{g,\mu}$, respectively.
  Note that the moduli space $\Mc_{g,\mu}$ of {\it smooth} Riemann surfaces,
  ie, surfaces with the  empty set $D$ of double points,
  is contained in both spaces  $\overline\Mc_{g,\mu}$ and
  $\overline\Mc^{\$}_{g,\mu}$, so that the natural projection
  $\overline\Mc^{\$}_{g,\mu}\to \overline\Mc_{g,\mu}$ is the identity
  map on $\Mc_{g,\mu}\subset \overline\Mc^{\$}_{g,\mu}$.
  Let us consider the oriented blow-up $S^{D}$ at the points
  of $D$ as described in Section \ref{sec:blow-up} above. The circles
$\oG_i$ and $\uG_i$
  introduced in (\ref{eq:decorat}) serve as the boundary circles
   corresponding  to the points $\od_i,\ud_i\in D$.
    The canonical projection  $\pi\co  S^{D}\to S$,
  which collapses the circles $\oG_i$ and $\uG_i$ to the points
  $\od_i$ and $\ud_i$, induces on $\Int S^{D}$ a conformal
  structure. The   smooth
  structure of $\Int S^{D}=S\setminus D$   extends to $S^D$, while the
extended  conformal structure
  degenerates along the  boundary circles $\oG_i$ and $\uG_i$.
Given a {\it decorated}  nodal surface $(\bS,r)$,  where
$r=(r_1,\dots,r_k)$, we can glue  $\oG_i$ and $\uG_i$ by means of the
mappings  $r_i$, for
$i=1,\dots,k$,
  to obtain  a closed surface $S^{D,r}$. As it was pointed out in
Section \ref{sec:blow-up}, the special circles
  $\Gamma_i=\{\oG_i\sim\uG_i\}$ are endowed with the
  canonical metric.
  The genus of the surface $S^{D,r}$ is equal to the
  arithmetic genus of $\bS$. There exists a
  canonical projection
  $p\co  S^{D,r}\to \widehat S_{D}$
  which projects the circle  $\Gamma_i=\{\oG_i,\uG_i\}$  to the double point
  $d_i=\{\od_i,\ud_i\}$.
  The projection $p$ induces on the surface  $S^{D,r}$ a conformal
   structure, still denoted by $j$,  in the complement of the special circles
  $\Gamma_i$. The continuous extension of $j$ to $S^{D,r}$ degenerates
along the special  circles
  $\Gamma_i$.
   The uniformizing metric
  $h^{j,M\cup D}$ can also be lifted to a metric $h^{\bS}$ on
   $\dot {S}^{D,r}=S^{D,r}\setminus M$.
   The lifted metric
  degenerates along
  each  circle $\Gamma_i$, namely
  the length of  $\Gamma_i$ is  $0$,  and the
   distance of $\Gamma_i$ to any other point
  in   $\dot S^{D,r}$ is infinite. However,  we can still  speak about
geodesics on $\dot S^{D,r}$
   orthogonal to $\Gamma_i$.
Namely, two geodesic rays,  whose asymptotic directions at the cusps
$\od_i$ and $\ud_i$
  are related
  via the map $r_i$,  correspond to a compact geodesic
  interval in $S^{D,r}$,   which
   orthogonally intersects
   the circle $\Gamma_i$. Notice that     the smooth structure on the
oriented blow up
   surface $S^D$ with boundary   is compatible with some smooth structure
   on $S^{D,r}$. Moreover, using the  hyperbolic metric one can make this choice
   canonical.
   However, this smooth structure will be irrelevant for us
   and we will not discuss here the details of its construction.
   It will be convenient for us to view
   $\Tn_\e(\bS)$ and $\Tk_\e(\bS)$ as subsets
   of $\dot S^{D,r}$. This interpretation provides us with a
    compactification of non-compact  components  of $\Tn_\e(\bS)$ not
adjacent to points
    from $M$.
    Every compact component $C$ of $\overline{\Tn_\e(\bS)}\subset
S^{D,r}$   is a compact annulus.
    It contains either a closed
    geodesic $\Gamma_C$, or one of
     the special circles which will
     also be denoted by $\Gamma_C$. This special circle projects to a
node as described above.
The surface $S^{D,r}$ with all the
  endowed structures and the projection $S^{D,r}\to\wh S_D$
   is called the {\it deformation} of the decorated nodal
  surface $(\bS,r)$. One can also define a partial deformation
$S^{D',r'}$ of $(\bS,r)$
  by splitting the set $D$ into a disjoint union
  $D=D'\cup D^{''}$ which respects the pairing structure, and then
  applying the above construction to $D'$ while adjoining $D^{''}$
  to the set $M$ of the marked points.

\subsection{Topology of spaces $\overline{\Mc}_{g,\mu}$ and
$\overline{\Mc}^{\$}_{g,\mu}$ }
  \label{sec:DM-topology}

  \noindent In this section we define    the meaning of the convergence in the
spaces $\overline\Mc_{g,\mu}$
and $\overline\Mc^{\$}_{g,\mu}$. The introduced topologies are
compatible with certain metric structures on the spaces
$\overline\Mc_{g,\mu}$ and $\overline\Mc^{\$}_{g,\mu}$ which we
discuss in Appendix \ref{ap:metric-surfaces}. In particular,  the introduced
topologies are Hausdorff.
Consider a sequence of  decorated stable nodal marked Riemann
surfaces $$(\bS_n,r_n)=\{S_n,j_n,M_n,D_n,r_n\},\;\;n\geq 1\,.$$
The sequence $(\bS_n,r_n)$  is said to converge to a decorated stable
nodal surface $$(\bS,r)=(S,j,M,D,r)$$ if (for sufficiently large
$n$)
   there exists a sequence of diffeomorphisms
$$\varphi_n \co S^{D,r}\to S_n^{D_n,r_n}\;\;\hbox{
with}\;\;\varphi_n(M)=M_n$$ and such that
the following conditions are satisfied.
\begin{description}
\item{\bf CRS1}\qua For every $n\geq 1$,  the  images $\varphi_n(\Gamma_i)$  of the
special circles $\Gamma_i\subset S^{D,r}$  for $i=1,\ldots, k$, are
special circles or
closed geodesics of the metrics $h^{j_n,M_n\cup D_n}$ on $\dot
S^{D_n,r_n}$. Moreover, all special circles on $S^{D_n,r_n}$ are
among these images.
\item{\bf CRS2}\qua $h_n\to h^{\bS} \:\:
    \text{in $C^{\infty}_{\text{loc}}\bigl(S^{D, r}\setminus (M\cup
\bigcup\limits_1^k\Gamma_i)\bigr)$},$
where $h_n=\varphi^*_nh^{j_n, M_n\cup D_n}$.
  \item{\bf CRS3}\qua Given a   component $C$ of
   $\Tn_\e(\bS)\subset \dot S^{D,r}$ which contains a special circle
   $\Gamma_i$ and given a point  $c_i\in\Gamma_i$,
   we consider for every $n\geq 1$  the geodesic arc $\d^n_i$ for the
induced  metric
$h_n=\varphi^*_nh^{j_n, M_n\cup D_n}$ which  intersects $\Gamma_i$
orthogonally
  at the point $c_i$,
and whose  ends  are contained in the $\varepsilon$--thick part of
  the metric $h_n$. Then $(C\cap\d^n_i)\;$
  converges as $n\to\infty$ in  $C^0$  to a
{\it continuous} geodesic for the metric  $h^{\bS}$  which passes
through the point $c_i$.
    \end{description}
\begin{remark}{\rm
Let us point out  that  in view of the Uniformization Theorem, the
condition CRS2
  is equivalent to the
condition $$\varphi_n^* j_n\to j \;\; \text{in
$C^{\infty}_{\text{loc}}\bigl(S^{D, r}\setminus (M\cup
\bigcup\limits_1^k\Gamma_i)\bigr)$}.$$ Moreover,
the Removable Singularity Theorem  guarantees that the latter
condition is equivalent to
$$\varphi_n^* j_n\to j \;\; \text{in
$C^{\infty}_{\text{loc}}\bigl(S^{D, r}\setminus
\bigcup\limits_1^k\Gamma_i\bigr)$}.$$
  }
  \end{remark}
\begin{figure}[ht!]\small
\centering
  \psfrag{gi}{$\Gamma_i$}
  \psfrag{x}{$c_i$}
\psfrag{d}{$\delta_{i}^n$}
\includegraphics[width=0.7\textwidth]{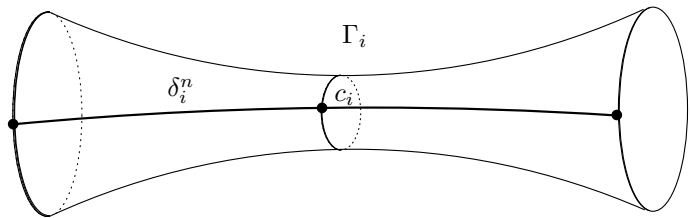}
\caption{Illustration for property CR3}\label{fig4}
\end{figure}
A sequence  $\bS_n\in\overline\Mc_{g,\mu}$ is said to converge to
$\bS\in\overline\Mc_{g,\mu}$ if
  there  exists a sequence    of decorations $r_n$ for $\bS_n$ and a decoration
  $r$ of $\bS$ such that
$ (\bS_n,r_n)$ converges to $(\bS,r)$ in
$\overline\Mc^{\$}_{g,\mu}$. In other words the topology on
$\overline\Mc_{g,\mu}$ is defined as the weakest topology on
$\overline\Mc_{g,\mu}$ for which the projection
$\overline\Mc^{\$}_{g,\mu}\to\overline\Mc_{g,\mu}$ is continuous.
Note that the convergence in the space $\overline\Mc_{g,\mu}$
can also be  defined by the  properties CRS1 and CRS2 above.

\begin{figure}[ht!]\small
\centering
\psfrag{i}{i)}
\psfrag{ii}{ii)}
\psfrag{iii}{iii)}
\psfraga <-3pt, 0pt> {y1}{$y^{(1)}$}
\psfraga <-3pt, 0pt> {y2}{$y^{(2)}$}
\psfrag{om}{$\text{old $m$}$}
\psfrag{m}{$m$}
\psfrag{z}{$z_0$}
\psfrag{S}{$S$}
\psfraga <-4pt, 2pt> {ox}{$\text{old $x$}$}
\psfraga <-2pt, 2pt> {oxp}{$\text{old $x'$}$}
\includegraphics[width=0.7\textwidth]{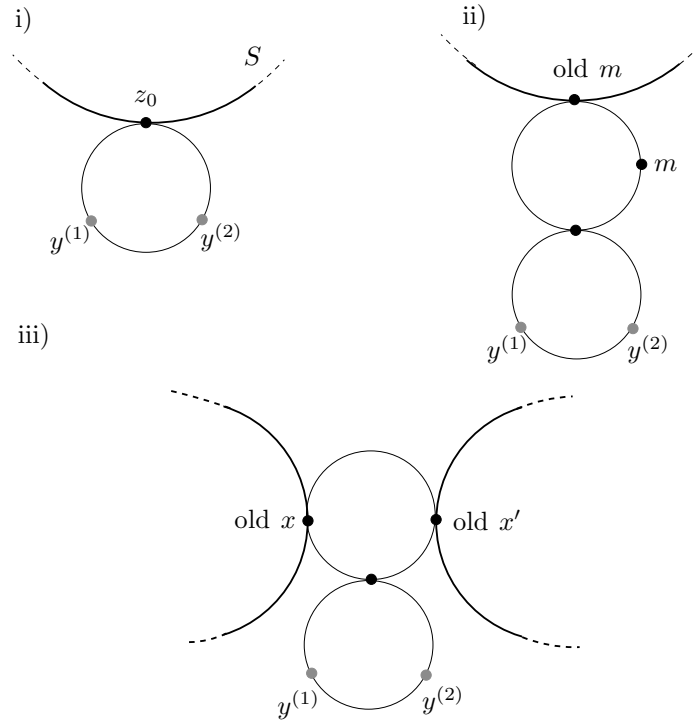}
\caption{The possible configuration  of  spherical bubbles  appearing
in Proposition
\ref{prop:2points}.}\label{Newfig1}
\end{figure}

   \begin{theorem}[Deligne--Mumford \cite{DM}, Wolpert \cite{Wolpert}]
\label{thm:DM} 
   
  The spaces $\overline\Mc_{g,\mu}$ and $\overline \Mc^{\$}_{g,\mu}$ are
  compact metric spaces, and serve as the compactifications  of the
space $\Mc_{g,\mu}$,
  ie, they coincide with the closure of  ${\Mc_{g,\mu}}$
  viewed as a subspace of  $\overline\Mc_{g,\mu}$ and of $\overline
\Mc^{\$}_{g,\mu}$,
  respectively.
  In particular, every sequence of smooth marked Riemann surfaces
$\bS_n=(S_n,j_n,M_n)$ of signature $(g,\mu)$ has
  a subsequence which  converges to a decorated nodal curve $\bS=(S,j,M,D,r)$
  of signature $(g,\mu)$.
  \end{theorem}
The next proposition illustrates the geometry  of the Deligne--Mumford
convergence in a special
case when one varies the configuration of  the marked points.
It follows from the definition of this convergence and from scaling
operations on the limit surface.
\begin{proposition}\label{prop:2points}
Let $\bS_n=(S_n,j_n,M_n,D_n)$ be a sequence of smooth marked nodal
Riemann surfaces
  of signature $(g,\mu)$ which converges to a  nodal curve $\bS=(S,j,M,D)$
  of signature $(g,\mu)$. Suppose that for each $n\geq 1,$ we are given
   a pair of points $Y_n=\{y^{(1)}_n,y^{(2)}_n\}\subset
  S_n\setminus (M_n\cup D_n)$
  such that  $$\mathrm{dist}_n(y^{(1)}_n,y^{(2)}_n)\xrightarrow[n\to
\infty] {}0.$$ Here $\mathrm{dist}_n$
  is the distance with respect to the hyperbolic metric
$h^{j_n,M_n\cup D_n}$ on $S_n\setminus (M_n\cup D_n)$.
  Suppose, in addition, that there is a  sequence
$R_n\mathop{\to}\limits_{n\to\infty} \infty$
  such that there exist injective  holomorphic maps
$\varphi_n\co D_{R_n}\to S_n\setminus (M_n\cup D_n)$,   where
   $D_{R_n}$ denotes the disc
  $\{|z|\leq R_n\}\subset\C$, satisfying  $\varphi_n(0)=y^{(1)}_n$ and
$\varphi_n(1)=y^{(2)}_n$.
  Then there exists a subsequence of the new  sequence
$\bS'_n=(S_n,j_n,M_n\cup Y_n,D_n)$ which converges to
  a nodal curve $\bS'=(S', j', M', D')$ of signature $(g,\mu +2)$,
which has one or two additional
   spherical  components. One of these components contains the marked
points $y^{(1)}$ and $y^{(2)}$
    which correspond to the sequences $y^{(1)}_n$ and $y^{(2)}_n$.
    The possible cases  are illustrated by Figure  \ref{Newfig1} and
described in detail below.
   \end{proposition}
    Let $r_n,r$ be some decorations of $\bS_n$  and $\bS$ such that
    $(\bS_n,r_n)\mathop{\to}\limits_{n\to\infty}(\bS,r)$.
    Let $\varphi_n\co S^{D,r}\rightarrow S^{D_n,r_n}_{n}$ be the sequence
of diffeomorphisms
  guaranteed by the  definition of  the convergence  of decorated
Riemann surfaces $\bS_n\rightarrow \bS$. Let  $\wh S_D$ be
  the singular surface with double points  as defined in
\eqref{eq:double}, and $\pi\co S^{D,r}\to\wh S_D$ the canonical
  projection.
  Set $Z_n=\pi(\varphi^{-1}(Y_n))\subset \wh S_D$.
Then the following scenarios are possible:
\begin{itemize}
\item[i)] The points $z^{(1)}_n,z_n^{(2)}\in Z_n$ converge to a
point $z_0$ which does not belong to $M$ nor  to $D$. Then the limit
{$\bS'$} of $\bS'_n$ has an extra sphere attached at $z_0$ on which
there are two
marked points $y^{(1)},y^{(2)}$, see Figure \ref{Newfig1}i).
\item[ii)] The points  $z^{(1)}_n,z_n^{(2)}\in Z_n$ converge to a
marked point $m\in M$.
In this  case the limit  $\bS'$ is  {$\bS$} with two extra spheres
$T_1$ and $T_2$ attached.
The sphere
$T_1$ is attached at its $\infty$ to the ``old" $m$ and has $m$ as
its $0$. $T_2$ is attached at its
$\infty$ to the point $1\in T_1$ and contains the marked points
$y^{(1)}$ and $y^{(2)}$, see Figure \ref{Newfig1}ii).
\item[iii)] The points in $Z_n$ converge to a double point $d$ which
corresponds to
  the node $\{x,x'\}\in D$. Then
between the nodal points $x$ and $x'$
a new sphere $T_1$ is inserted, say attached at  its $\infty$ to  $x$
and at $0$ to $x'$. At the point
  $1$ a second sphere $T_2$
carrying  two marked points  $y^{(1)},y^{(2)}$ is attached, see
Figure  \ref{Newfig1}iii).
\end{itemize}
  \section{Holomorphic curves in cylindrical almost complex
manifolds}\label{sec:holom-cylin}
\subsection{Gromov--Schwarz and monotonicity} \label{sec:Gromov-Schwarz}
    We  begin  by recalling two important  analytic results
    about  $J$--holomorphic maps (see \cite{Gr} and
\cite{Muller,Hummel} for the proofs).
   \begin{lemma}[Gromov--Schwarz]
   \label{lm:Gromov-Schwarz}
    Let $f \co  D^2(1) \rightarrow W$ be a
holomorphic disk in an almost complex manifold $W$ whose structure $J$ is
tamed by an exact symplectic form. If the image of $f$ is contained in a
compact set $K \subset W$, then
$$
\| \nabla^k f(x) \| \le C(K,k) \qquad \textrm{for all } x \in D^2(1/2),
$$
for every $k\geq 1$ where the constants do not depend on $f$.
\end{lemma}
  Note that locally  every  $J$ is tamed by some symplectic form,
  and hence Lemma \ref{lm:Gromov-Schwarz}
  holds for sufficiently small compact sets
  in every almost complex manifold.
\begin{lemma}[Monotonicity]\label{lm:monotonicity}
  Let $(W, J)$ be a compact almost complex manifold and suppose $J$ to
be tamed by $\omega$.
  Then there exists a positive constant      $C_0$ having the
following property.
  Assume that $f\co (S, j)\to (W, J)$ is a compact $J$--holomorphic
  curve with boundary and choose $s_0\in S\setminus \partial S$ and $r$ smaller
  that the injectivity radius of $W$. If the boundary $f(\partial S)$ is
   contained in the complement of the $r$--ball $B_r(f(s_0))\subset W$,
   then the area of $f$ inside of the ball $\overline{B}_r(f(s_0))$ satisfies
  $$
\int_{f(S)\cap B_r(f(s_0))} \omega \ge C_0r^2\,.
$$
\end{lemma}
\subsection{$\overline\partial$--equation on cylindrical manifolds}
\label{sec:equation}
Let $(W=\R\times V,J)$ be  a cylindrical almost complex manifold  and
let $(\xi,J_\xi)$,
$\Reeb\in TV$ and $\l$ be the corresponding CR--structure, the vector
field transversal to $\xi$ and the $1$--form
on $V$ determined by the conditions
$\xi=\{\l=0\}$ and $\l(R)=1$.
   We will denote by $p_{\R}$ and $p_V$ the projections to the first
and the second factor of $W$, and by
$\pi$  the projection $TV\to\xi$ along the direction of the  field $\Reeb $.
Let us agree on the following notational convention.
We shall use capital roman letters to denote maps to $W$, and corresponding
small letters for their projections to $V$, eg, if $F$ is a map to $W$, then
$p_V\circ F$ will be denoted by $f$. The $\R$--component $p_{\R}\circ F$ of $F$
  will usually be  denoted by $a$, or $a_F$ if necessary.
   Given a Riemann surface  $(\Sigma, j)$, the
$\overline\partial$--equation  defines the holomorphic maps (or curves)
$F=(a,f)\co (\Sigma,j)\to (W=\R\times V,J)$ as solutions of the equation
$$
TF\circ j=J\circ TF.
$$
For our distinguished structure $J$, the equation takes the form
\begin{equation}\label{main10}
\begin{gathered}
\pi\circ df\circ j=J\circ\pi\circ df \\
da=(f^*\l)\circ j\;.
\end{gathered}
\end{equation}
    Notice, that the second equation just means that the form $f^*\l\circ j$ is
exact on $\Sigma$ and that
the function $a$ is a primitive of the $1$--form ${f^*\l}\circ j$.
Thus the holomorphicity   condition
for $F=(a,f)$ is essentially just a condition on  its $V$--component $f$. If
$f$ satisfies the first of the equations (\ref{main10}) and the form
$(u^*\l)\circ j$ is exact then the coordinate $a$ can be reconstructed
uniquely up to an additive constant on each connected component of $\Sigma$.
  We will call the map $f\co \Sigma\to V$ a {\it holomorphic  curve in
  $V$} if $f$ satisfies the first of the equations (\ref{main10}) and the form
$(f^*\l)\circ j$ is exact.
   \subsection{Energy}\label{sec:energy}
   A crucial assumption in Gromov's compactness theorem for holomorphic curves
in a compact symplectic manifold is the finiteness of the area.
However, the area of a
non-compact
  proper holomorphic curve in  a
cylindrical manifold is never finite with respect to any complete metric.
  Moreover, in the contact case, or more generally in the case
  when the taming form $\omega$ is exact, there are no  non-constant compact
  holomorphic curves. We define below
  another quantity, called {\it energy}, which serves as a substitute
for the area in
  this case.
Suppose that the cylindrical almost complex structure $J$ is adjusted to a
closed form  $\omega$.
Let us  define the {\it $\omega$--energy} $E(F)$
of the holomorphic map $F=(a,f) \co  (S,j) \rightarrow
(\RR \times V, J)$ by the formula
\begin{equation}\label{eq:energy1}
E_\omega(F) =\int\limits_Sf^*\omega\,,
\end{equation}
and
the $\lambda$--energy
by the formula
\begin{equation}\label{eq:energy2}
E_\lambda(F)=\sup_{\phi\in\Cc  } \int\limits_S (\phi\circ a)da\wedge f^*\l \,,
\end{equation}
where the supremum is taken over the set $\Cc$ of all non-negative
$C^\infty$--functions
  $\phi\co \R\to \R$ having  compact support  and  satisfying  the condition
  $$ \int\limits_{\bR}\phi(s)ds=1\,. $$
  It is important to observe that the $\omega$--energy is additive,
while the $\l$--energy is not.
  Finally, the {\it energy} of $F$ is defined as the sum
  $$E(F)=E_\omega(F)+E_\lambda(F).$$
   Note that the $\omega$--energy $E_\omega(F)$ depends only on  the
$V$--component $f$ ofthe map $F$.
The following inequality is a straightforward consequence of the
definition of the $\lambda$--energy $E_\l(F)$.
\begin{lemma}\label{lm:action-bound}
For any holomorphic  curve $F=(a,f)\co (S,j)\to (\R\times V,J)$ and any
non-critical level  $\Gamma_c=\{a=c\}$,
$$\biggl|\int_{\Gamma_c}\lambda\biggr|\leq E_\l(F)\,.$$
\end{lemma}
The adjustment condition implies the following lemma.
\begin{lemma}\label{lm:positivity}
For any holomorphic curve $F$ we have $E_\omega(F),E_\lambda(F)\geq 0$.
  Moreover, we have
  $ \int_S (\phi\circ a)da\wedge f^*\l\geq 0$ for any $\phi\in\Cc$.
If $E_\omega(F)= 0$,  then the image $f(S)$ is contained in
a trajectory of the vector field $\Reeb$.
\end{lemma}
  \subsection{Properties of holomorphic cylinders}\label{sec:asymptotics}
We begin with
  a version of  the Removable Singularity Theorem in
our context.
\begin{lemma} \label{punctfin}
Let $F=(a,f) \co  (D^2 \setminus \{0 \}, j) \rightarrow (\RR \times
V, {J})$ be a
holomorphic map with $E(F) < \infty$. Suppose that the image
of $F$ is contained in a compact subset $K$ of $\RR \times V$.
Then $F$ extends continuously
to a holomorphic map $\widetilde F \co  (D^2,j) \rightarrow (\RR \times V, J)$
with $E(\widetilde{F}) = E(F)$.
\end{lemma}
\begin{proof}
Choose $\phi \in \Cc$ such that $\phi(t) > 0$ for all $(t,x) \in K$.
Then,  the (not necessarily closed) $2$--form
$\Omega=\pi^*\omega+\phi(t)dt\wedge\l$
  is non-degenerate on $K$ and tames $J$. Thus the energy bound for $F$
implies its area bound,
and hence, we can apply
the usual Removable Singularity Theorem (see for example \cite{MS2}). \end{proof}
The next proposition,  proven in \cite{HWZ1, HWZ6} in the
non-degenerate case and
in \cite{B} in the Morse--Bott case,
describes the behavior of finite energy holomorphic
curves near the punctures in the case of  non-removable singularities.
\begin{proposition} \label{lm:asymptot}
Suppose that the vector field $\Reeb$ is of Morse, or Morse--Bott type.
Let $F = (a,f) \co  \RR^+ \times \R/\Z \rightarrow (\RR \times V, J)$ be a
holomorphic map of finite energy  $E(F) < \infty$. Suppose that the
image of $F$ is
unbounded in $\RR \times V$. Then there exist a number $T\neq 0$ and
a periodic orbit $\gamma$
of $\Reeb$ of period $|T|$ such that $$\lim_{s \rightarrow \infty} f(s,t) =
\gamma( Tt)\quad\hbox{and}\quad
\lim_{s \rightarrow \infty}\frac {a(s,t)}s = T\quad\hbox{in}\quad
C^\infty(S^1)\,.$$
\end{proposition}
\begin{proposition}\label{prop:twist}
Given $E_0$, $\varepsilon>0$, there exist constants $\sigma, c>0$
such that for every $R>c$ and every holomorphic cylinder
$F=(a, f)\co [-R, R ]\times S^1\to {\mathbb R}\times U$ satisfying  the
inequalities
$$E_{\omega}(F)\leq \sigma\quad \text{and}\quad E(F)\leq E_0,$$
we have
$$f(s, t)\in B_{\varepsilon}(f(0, t))$$
for all $s\in [-R+c, R-c]$ and all $t\in \bR$.
\end{proposition}
  More precisely, there exists either a periodic orbit $\gamma$ of
the vector field $\Reeb$ on $V$ having period $T>0$ such that $f(s,
t)\in B_{\varepsilon}(\gamma (Tt))$ or a point $p^{\ast}\in \R\times
V$ such that $F(s, t)\in B_{\varepsilon}(p^*)$ for all $s\in [-R+c,
R-c]$ and all $t\in \R$. Hence a long holomorphic cylinder having
small $\omega$--energy is either close to a trivial cylinder over a
periodic orbit of the vector field $\Reeb$, or close to a constant
map. The proof in the non-degenerate case can be found in
in \cite{HWZ7}. In the Morse--Bott case the proof is given in
Appendix \ref{ap:estimates}.

\subsection{Proper holomorphic maps of punctured Riemann
surfaces}\label{sec:punctures}
Proposition  \ref{lm:asymptot} implies:
  \begin{proposition}\label{prop:punct}
  Let $(S,j)$ be a closed Riemann surface
  and let $$Z=\{z_1,\dots, z_k\}\subset S$$ be a set of punctures.
Every holomorphic map
$F=(a,f)\co (S\setminus Z,j)\to\R\times V$ of finite energy  and without
removable singularities
is asymptotically cylindrical near each puncture $z_i$ over
a periodic orbit $\g_i\in\Pc$.
\end{proposition}
  The puncture $z_i$
is called
{\it positive or negative} depending on the sign of the coordinate  function
  $a$
when approaching the puncture.   Notice that the change of the
holomorphic coordinates near the punctures
affects only the choice of the origin on the orbit $\g_i$;
  the parametrization of the asymptotic orbits induced by the
holomorphic polar coordinates
  remains  otherwise the same.
   Hence, the orientation induced on $\gamma_i$ by the holomorphic
    coordinates coincides with the orientation defined by the vector field
    $\Reeb$ if and only if the puncture is positive.
\begin{remark}
{\rm In the situation when $V$ is as in Example
\ref{ex:circle-bundle},  a holomorphic map
$F=(a,f)\co (S\setminus Z,j)\to\R\times V$  can be extended via the
Removable Singularity Theorem into a holomorphic map
$\overline F\co S\to \overline W$, where $\overline W$ is the
projectivization of the complex line bundle associated with the circle
bundle $V\to M$. Punctures $z_j$ are then mapped by $\overline F$ to
the divisors of $\overline W$ which correspond
to $0$-- and $\infty$--sections $M_0,M_\infty\subset \overline W$.
The multiplicity $k_j$ of the asymptotic orbit of $\Reeb$
corresponding to a puncture $z_j\in S$ equals  $t_j+1$
where $t_j$ is the order of tangency of $\overline F$ to $M_0\cup
M_\infty$ at $z_j$. This observation explains
why the compactness theorems proven in \cite{Ionel-Parker,Li,Ruan-Li}
follow from the results of this paper.}
\end{remark}
In  what follows we will only
consider  holomorphic maps of
Riemann surfaces  which are
conformally       equivalent to a compact Riemann surface $(S,j)$ with {
\it punctures} $Z=\{z_1,\dots,z_k\}$.
Let $S^Z$ be the oriented blow-up of $S$ at the points of $Z$, as it
is defined  in
Section  \ref{sec:blow-up} above. Thus $S^Z$ is a compact surface with
boundary consisting of circles $\Gamma_1,\dots,\Gamma_k$. Each of these circles
is endowed with a canonical $S^1$--action and we denote by
$\varphi_j\co S^1\to\Gamma_j$ the canonical
(up to a choice of the base point) parametrization of the boundary
circle $\Gamma_j$,  for $j=1,\dots,k$.
   Proposition  \ref{prop:punct} can be equivalently reformulated as follows.
    \begin{proposition}\label{prop:extension}
    Let $F=(a,f)\co (S\setminus  Z,j) \to (\R\times V,J)$ be a finite energy
    holomorphic  map without removable singularities.
         Then the map $f\co S\setminus Z\to V$ extends to a continuous map
    $\of \co S^Z\to V$ satisfying
    \begin{equation}\label{eq:param}
    \of\left(\varphi_j(e^{it})\right)=\gamma_j\left(\pm{Tt}\right),
    \end{equation} where
    $\gamma_j\co S^1=\R/\Z\to V$ is a periodic orbit of the vector field
$\Reeb$ of period
    $T$, parametrized by the vector field $\Reeb$. The sign  in the formula
    (\ref{eq:param}) coincides with the sign of  the puncture $z_j$.
    \end{proposition}
    We will call the map $\of\co S^Z\to V$  the {\it compactification }
of the map $f$.
\subsection{Bubbling lemma}
\begin{lemma} \label{lm:bubble}
There exists a constant $\hbar$ depending only on $(V, J)$ so that
the following holds true. Consider
  a sequence $F_n\co =(a_n, f_n):D\to (\bR\times V, J)$ of $J$--holomorphic maps
of  the  unit disc $D=\{|z|< 1\}\subset \C$ to $V$ satisfying
$E(F_n)\leq C$ for some constant $C$ and such that $a_n(0)=0$. Fix a
Riemannian metric on
  $V$. Suppose that  $\| \nabla F_n(0) \| \to \infty$ as $n \to
\infty$.  Then there exists a sequence of points $y_n\in D$
converging to $0$, and
sequences  of positive
  numbers $c_n, R_n\to \infty$ as $n\to \infty$ such  that
$||y_n||+c_n^{-1}R_n< 1$ and the
  rescaled maps
\begin{equation*}
\begin{split}
F^0_n  \co  D_{R_n}=\{|z|<R_n\}& \rightarrow (\R \times V, J),\\
z&\mapsto F_n(y_n+c^{-1}_nz),
\end{split}
\end{equation*}
converge  in
$C^{\infty}_{\text{loc}}(\C)$ to a holomorphic map $F^{0}\co\C\to
\bR\times V$ which satisfies the conditions
  $$E(F^0)\leq C\;\;\hbox{ and}
\;\;E_{\omega}(F^0)>\hbar.$$
  Moreover, this map is either a
holomorphic sphere or a holomorphic plane $\bC$ asymptotic as
$\abs{z}\to \infty$ to a periodic orbit of the vector field
$\Reeb$.  (To be precise we mean in the first case that $F^0$
smoothly extends to $S^2=\bC\cup \{\infty\}$).
\end{lemma}
  For the proof  of Lemma
\ref{lm:bubble} we shall need the following lemma from \cite{HV}.
\begin{lemma} \label{lm:Hofer}
Let $(X,d)$ be a complete metric space, $f \co  X \rightarrow \RR$ be
a nonnegative continuous function, $x \in X$, and $\delta > 0$.
Then there exist $y \in X$ and a positive number $\varepsilon \le
\delta$ such that
$$
d(x,y) < 2\delta, \qquad \sup_{B_\varepsilon(y)} f \le 2f(y), \qquad
\varepsilon f(y) \ge \delta f(x) .
$$
\end{lemma}
\begin{proof}[Proof of Lemma \ref{lm:bubble}]
Choose  $\delta_n >0$ such that $$\delta_n \rightarrow 0,\;\;\hbox{ and } \;\;
\delta_n \| \nabla F_n(0) \| \rightarrow \infty\,.$$ Applying
Lemma \ref{lm:Hofer}, we obtain new sequences $y_n \in S$ and $0 <
\varepsilon_n \le \delta_n$ such that $y_n \rightarrow x$ and
$$
\sup_{|z-y_n|\leq\varepsilon} \| \nabla F_n (z)\| \le 2 \| \nabla
F_n(y_n) \| , \qquad \varepsilon_n \| \nabla F_n(y_n) \| \rightarrow
\infty .
$$
Introduce $c_n = \| \nabla F_n(y_n) \|$ and $R_n = \varepsilon_n
c_n$.  Notice that for sufficiently large $n$ we have
$\|y_n\|+c_n^{-1}R_n<1 $. We  consider
  the rescaled maps $$F^0_n(z) = F_n(y_n + c_n^{-1}z).$$
  This sequence has  the following properties:
\begin{itemize}
\item $ \sup_{D_{R_n}} \| \nabla F^0_n(z) \| \le 2 , \qquad R_n
\rightarrow \infty $; \item  $E(F^0_n)\leq E(F_n)\leq C$; \item $
\| \nabla F^0_n(0) \|=1$.
\end{itemize}
Now, by Ascoli--Arzela's theorem, we can extract a converging
subsequence and thus we obtain a non-constant finite energy plane
$F^0$. If the  image of $F^0$ is contained in a compact subset of
$\RR \times V$, then by Lemma \ref{punctfin}, $F^0$ is a
holomorphic sphere. Otherwise,  we can apply Proposition
\ref{lm:asymptot} to
  deduce that $F^0$ is converging to an $\Reeb$--orbit $\gamma$ for large radius.
In both cases there is a constant $\hbar$ such that
$E_\omega(F^0)>\hbar>0$. Indeed, otherwise we could get a sequence
$F^{i}$ of holomorphic planes  satisfying $E_\omega(F^{i})\to
0$ as $i\to \infty$, having a uniform  gradient bound  and
normalized by  the condition $\|\nabla F^{i}(0)\|=1$. Then  Ascoli--Arzela's
theorem would imply the existence of a  non-constant limit
holomorphic plane $F^\infty\co \bC\to \bR\times V$ satisfying
$E_\omega(F^\infty)=0$ and $E(F^{\infty})<+\infty$. In view of Lemma
\ref{lm:positivity}
we then conclude that $F^{\infty}(\C)$  is contained in a cylinder
over an orbit of $\Reeb$, and hence
must coincide with the universal covering of this cylinder. But this
is impossible in view of finiteness
of the energy $E(F^{\infty})$. This finishes off the proof.
\end{proof}
\subsection{The symmetric case}
The next lemma shows that for a symmetric $J$, and in particular in
the contact case,
  the energies of holomorphic curves $F=(a,f)$
asymptotic to  prescribed periodic orbits from $\Pc$
  can  be uniformly bounded  in terms of the relative homology class represented
  by $f$.
  \begin{proposition}\label{prop:automatic-bound}
Let $(\R\times V,J)$ be a symmetric cylindrical almost complex structure
adjusted to a closed $2$--form $\omega$ on $V$.
  Suppose that the  holomorphic curve
$F\co (S\setminus Z,j)\to (\R\times V,J)$  of finite energy is
  asymptotic at the positive punctures to the
periodic orbits $\og_1,\dots,\og_k\in\Pc$ and to the periodic orbits
$\ug_1,\dots,\ug_l \in \Pc$ at  the negative punctures.
Then there exists a positive constant $C$ (which depends on $J,\lambda,\omega$
but not $F$)
such that
\begin{equation}
E(F)\leq C\int\limits_Sf^*\omega+\sum\limits_1^k S(\og_i).
\end{equation}
In particular,
the energies $E(F)$ are  uniformly bounded for all $F$ for which
$f$ represents a given homology class in $H_2(V,\bigcup\og_i\cup\ug_j)$.
\end{proposition}
\begin{proof}
First observe that
\begin{equation}
\label{eq:lambda-omega}
  |f^*d\lambda|\leq C f^*\omega\,.
\end{equation}
    Indeed, in the symmetric case $\Reeb\hook\  d\lambda=\Reeb\hook\ \omega=0$,
    and hence the value of both forms on any bivector $\sigma$ is equal
    to the value of these forms on the projection of $\sigma$ to $\xi$.
    But a complex direction of $J$ projects to a complex direction of $J_\xi$,
    and hence the inequality follows from
    the Wirtinger inequality due to the fact that $J_\xi$ is tamed by
$\omega|_\xi$.
Given any function $\phi\in \Cc$ we find, using Stokes' formula,
\begin{equation}\label{eq:Stokes2}
\begin{split}
\left|\int\limits_S (\phi\circ a)da\wedge f^*\l\right|
&\leq \left|\int\limits_SF^*d(\psi\l)\right|+\left|\int\limits_S
(\psi\circ a)f^*d\l\right|\cr
&\leq \sum\limits_1^k S(\og_i)+\int\limits_S|f^*d\l|,\cr
\end{split}
\end{equation}
where $\psi(s)=\int\limits_{-\infty}^s\phi(\sigma)d\sigma$.
Using  \eqref{eq:lambda-omega} and  \eqref{eq:Stokes2}   we obtain
     \begin{equation}
     E_{\lambda}(F) =\sup\limits_{\phi\in\Cc}\int\limits_S (\phi\circ
a)da\wedge f^*\l \leq
      \sum\limits_1^k S(\og_i)+C\int\limits_S  f^*\omega  \,.
    \end{equation}
\vspace{-30pt}

\end{proof}
We will also need the following property of holomorphic cylinders
with small boundary circles.
\begin{lemma}\label{lm:almost-sphere}
Let $F_n=(a_n,f_n)\co [-n,n]\times S^1\to\R\times V$ be a family of
holomorphic cylinders.
Suppose that $$E_\omega(F_n)\to 0\;\;\hbox{ and}\;\;
\lim\limits_{n\to\infty} F^\pm_n|_{\pm n\times S^1}=z_\pm\in\R\times
V\:  \text{in $C^1(S^1)$},$$
where $z_{\pm}$  are two points and where the maps $F^\pm_n$ differ
from $F_n$ by a translation with a sequence of constants,
$$F^\pm_n=(a_n-c_n^\pm ,f_n).$$
Then $$\diam F_n([-n,n]\times S^1))\to 0\:\: \text{ as $n\to \infty$}.$$
\end{lemma}
\begin{proof}
The cylindrical almost complex structure $J$  on $\R\times V$ is
tamed by an almost symplectic, ie, non-degenerate
but not necessarily closed,  differential $2$--form
$\Omega=\omega +dg(t)\wedge\lambda$, where $g\co \R\to(0,\varepsilon)$
is  $C^\infty$--function with a positive derivative.
The Monotonicity Lemma \ref{lm:monotonicity} implies the existence of
a constant $C>0$ such that  for a sufficiently small
  $\varepsilon>0$, for every holomorphic curve $S$ and every  point
$x\in S$ such that the ball $B_{\varepsilon}(x)$ does not intersect
the boundary
  of $S$,  we have (possibly after  translating $S$ along the
$\R$--direction) the inequality
  $$\int\limits_{S\cap B_{\varepsilon}(x)}\Omega\geq C\varepsilon^2\,.$$
   Suppose that $\diam F_n([-n,n]\times S^1)\geq\delta>0$. Then there
exists a point $y_n\in [-n,n]\times S^1$ such that
   $$\dist\left(F_n(y_n),F_n\left(\partial([-n,n]\times
S^1)\right)\right)\leq\frac\delta2\,.$$
   Choosing $\varepsilon<\delta/2$ we conclude, after possibly
translating the cylinder $F_n$, that
   \begin{equation}
   \label{eq:lower-bound}
   \int\limits_{[-n,n]\times S^1}F_n^*\Omega\geq
   \int\limits_{[-n,n]\times S^1\cap
F_n^{-1}(B_{\varepsilon}(F_n(y_n))}F_n^*\Omega\geq C\varepsilon^2\,.
   \end{equation}
   On the other hand,  by Stokes' theorem,
   \begin{equation}
   \label{eq:Stokes3}
   \int\limits_{[-n,n]\times S^1}F_n^*\Omega=\int\limits_{[-n,n]\times
   S^1}F_n^*(\omega-g(t)d\lambda)+
\int\limits_{F_n\left(-n\times S^1\cup n\times S^1\right)}g(t)\lambda\,.
\end{equation}
The second term on the right-hand side of this equality converges to
$0$. On the other hand,
the symmetry condition  and the Wirtinger inequality imply that
$$\int\limits_{[-n,n]\times S^1}|F_n^*d\lambda|\leq C\int\limits_{[-n,n]\times
S^1}F_n^*\omega\leq
CE_{\omega}(F_n)\to 0\,$$
as $n\to \infty$. Hence,
the right-hand side  of (\ref{eq:Stokes3}) converges to $0$, which
contradicts the positive lower bound (\ref{eq:lower-bound}) and
completes the proof of the lemma.
   \end{proof}
\subsection{The contact case: relation between energy and contact
area}\label{sec:relation}
We consider in this section the contact case, ie, we assume that
the $1$--form $\lambda$ implied by the definition of a cylindrical
structure $J$ is a contact form,  and  $J$ is adjusted to
$d\l=\omega$. We shall denote  the $d\l$--energy
$E_{d\l}(F)$ by $A(F)$ and call it the {\it contact area}.
\begin{lemma} \label{equiv}
Let $F=(a, f) \co  (S \setminus Z,j) \rightarrow (\RR \times V, J)$ be a
holomorphic map. Then the
following two statements are equivalent,
\begin{enumerate}
\item [\rm(i)] $A(F) < \infty$ and $F$ is a proper map;
\item[\rm(ii)] $E(F) < \infty$ and $S$ has no removable punctures.
\end{enumerate}
\end{lemma}
\begin{proof} (i) $\Rightarrow$ (ii).
By properness of $F$, the limit of the coordinate function $a$ near
each puncture
is either $+\infty $ or $-\infty$, and hence all  the   punctures can
be divided into positive and negative punctures,
according to the particular  end of $\bR \times V$ which the
holomorphic curves approaches near the puncture.
In a neighborhood $U$ of a puncture $p$, let $z$ be a complex
coordinate vanishing at $p$. Let $D_r(p) = \{ q \in U \, | \;
|z(q)| \le r\}$ and $C_r(p) = \partial D_r(p)$ oriented
counter-clockwise  for a positive
   puncture $p$, and clockwise for a negative one. Consider
$\int\limits_{C_r(p)} f^*\lambda$ as a function
of $r$. It is increasing and bounded above (resp. decreasing and
bounded below) if the puncture
is positive (resp. negative), since $d\lambda \ge 0$ on complex lines and
$\int\limits_{D_r(p)}d\l<C$.
Hence $C_r(p)$ has a finite limit for $r \rightarrow 0$  for all,
positive and negative
punctures.
Now, let $\phi \in \mathcal{C}$ and let $\phi_n \in \mathcal{C}$ such
that $\phi_n \circ a=0$   in $D_{\frac1n}(p)$ for all punctures $p$.
Such functions exist, by properness of $F$. Moreover,
we can choose $\phi_n$ so that $\| \phi - \phi_n\|_{C^0} < \varepsilon_n$,
  with $\varepsilon_n \rightarrow 0$ for $n \rightarrow \infty$.
We have
\begin{equation}\label{eq:new-old-energy}
\int\limits_S (\phi_n\circ a)da\wedge  \lambda  =
\int\limits_S F^* d(\psi_n\lambda)-\int\limits_S(\psi_n\circ a)f^*d\l \,,
\end{equation}
where $\psi_n(s)=\int\limits_{-\infty}^s\phi_n(\sigma)d\sigma$.
Notice that $\psi_n\circ a=1$ in
$D_{\frac1n}(p)$ when $p$ is a positive puncture and $\psi_n\circ a=0$ in
$D_{\frac1n}(p)$ when $p$ is a negative one.
By Stokes theorem,
\begin{equation}
\int\limits_S F^* d(\psi_n\lambda)=\lim\limits_{r \rightarrow 0}
\sum\limits_{p } \int\limits_{C_r(p)} f^*\lambda\,,
\end{equation}
  where the sum is taken over all positive punctures $p$.
  Therefore,
\begin{equation}\label{eq:new-energy}
\begin{split}
\int\limits_S (\phi_n\circ a)da\wedge  \lambda &=\lim\limits_{r
\rightarrow 0}\sum\limits_{p } \int\limits_{C_r(p)} f^*\lambda -
  \int\limits_S(\psi_n\circ a)f^*d\l\cr
  &\leq
  \lim\limits_{r \rightarrow 0}\sum\limits_{p } \int\limits_{C_r(p)}
f^*\lambda < C'<+\infty \, .
  \cr
  \end{split}
  \end{equation}
Moreover,
$$\int\limits_S (\phi_n\circ a)da\wedge  \lambda\to \int\limits_S
(\phi \circ a)da\wedge  \lambda \, $$
as $n\to \infty$. Hence, $$\int\limits_S (\phi \circ a)da\wedge
\lambda\leq C',$$
  and thus $E(F) \leq\int\limits_Sd\l+ C' < \infty$.
\noindent (ii) $\Rightarrow$ (i). First, we obviously have
  $A(F) \le E(F) < \infty$. Moreover, $F$ has only positive and
negative punctures by assumption,
and thus the map $F$ is proper.
\end{proof}
The energy and the contact area of
  a holomorphic map $F\co (S\setminus Z,j)\to(\R\times V,J)$ of finite energy
    are easily computable in view of Stokes' formula and formula
(\ref{eq:new-energy}). The result is as follows.
  \begin{lemma} \label{value}
Under the condition (i) or (ii) of Lemma \ref{equiv}, we denote by
$\og_1, \ldots,
\og_k$ (resp. $\ug_1, \ldots, \ug_l$) the periodic orbits
of $\Reeb$ asymptotic to the positive (resp. negative) punctures of
$S$. Then, with $d\lambda=\omega$,
\begin{equation}\label{eq:energy-area}
\begin{split}
E_{\omega}(F)&=A(F) = \sum_{j=1}^k S(\og_j) - \sum_{j=1}^l S(\ug_j)\, \cr
E_{\l}(F) &= \sum_{j=1}^k S(\og_j)\cr
E(F) &= 2\sum_{j=1}^k S(\og_j)-\sum_{j=1}^l S(\ug_j). \cr
  \end{split}
\end{equation}
\end{lemma}
\section{Holomorphic curves in almost complex manifolds
with cylindrical ends}  \label{sec:curves-ends}
\subsection{Energy}
Let  $(W,J)$ be  an almost complex manifold with
cylindrical ends, ${W}=
E_-\cup \overline W\cup E_+$,  adjusted to a
symplectic form $\omega$ on $\overline W$.
The diffeomorphism
  $G^\d\co W\to\circW=\Int\overline W$
   defined in Section \ref{sec:attaching}
  allows us to identify $J$ with the almost complex structure
  $J^\d=G^\d_*J$ on $\circW$.
  Both
(equivalent) points
of view will be useful
for us.    We will denote by $J_\pm$  the cylindrical almost complex
structures which are restrictions of $J$
to the ends $E_\pm$, where
$E_-=(-\infty,0]\times V_-$ and $ E_+=[0,\infty)\times V_+$.
We need to generalize the definition of energy for a holomorphic
curve into an almost
complex manifold $(W,J)$ with cylindrical ends, when the structure
$J$ is adjusted
to a symplectic form $\omega$ on $\overline W$.
First, we define the {\it $\omega$--energy} by the formula
\begin{equation}\label{eq:energy21}
  E_\omega(F)=\int\limits_{F^{-1}(\overline W)}F^*\omega+
   \int\limits_{F^{-1}(E_-)}f_-^*\omega +
   \int\limits_{F^{-1}(E_+)}f_+^*\omega\,,
   \end{equation}
  where $F|_{E_\pm}=(a_\pm,f_\pm)$.
  Next, we define the $\l$--energy $E_\l(F)$ in a way similar to
formula (\ref{eq:energy2}) in the case of cylindrical manifolds.
\begin{equation}\label{eq:energy22}
E_\lambda(F)=\sup_{\phi_\pm\in\Cc   }
\left(\int\limits_{f^{-1}(E_+)} (\phi_+\circ a_+)da_+\wedge f^*\l+
\int\limits_{f^{-1}(E_-)} (\phi_-\circ a_-)da_-\wedge f^*\l\right) \,,
\end{equation}
where the supremum is taken over all pairs $(\phi_-,\phi_+)$ from
the set $\Cc$ of all   $C^\infty$--functions
$\phi_\pm\co \R_\pm\to\R_+$ such that
{
$$
\int\limits_{0}^{\infty}\phi_+(s)ds=
\int\limits_{-\infty}^{0}\phi_-(s)ds= 1.
$$
}
  Finally, the {\it energy} of $F$ is defined as the sum
  $$E(F)=E_\omega(F)+E_\lambda(F).$$
  Similarly to Lemma \ref{lm:positivity} for cylindrical manifolds we have
  \begin{lemma}\label{lm:positivity-ends}
For any holomorphic curve $F$ in $(W,J)$ we have
$$E_\omega(F),E_\lambda(F)\geq 0.$$
  Moreover,
$$ \int\limits_{f^{-1}(E_+)} (\phi_+\circ a_+)da_+\wedge f^*\l\geq 0$$ and
$$\int\limits_{f^{-1}(E_-) }(\phi_-\circ a_-)da_-\wedge f^*\l\geq 0$$
for all  admissible functions $\phi_\pm$.
  \end{lemma}
\subsection{Asymptotic properties of holomorphic curves}
Straightforward extensions of
  Lemmas \ref{punctfin} and \ref{lm:asymptot}
  allow us to describe the asymptotic behavior of $f$ near the
punctures  of $S$ as follows.
\begin{proposition} \label{prop:punct-cob}
Let $(S,j)$ be a closed Riemann surface,  $Z=\{z_1,\dots, z_k\}\subset S$
   a set of punctures, and $S^Z$  the oriented blow-up at the points of $Z$.
Any holomorphic map
$F\co (S\setminus Z,j)\to {W}$ of finite energy and without removable
singularities
is asymptotically cylindrical near each puncture $z_i$ over
a periodic orbit $\g_i\in\Pc=\Pc_-\cup\Pc_+$.
  The map  $F^\d=G^\d\circ F\co (S\setminus Z,j)\to \circW$
  extends to a smooth map $\overline F^\d\co S^{Z}\to \overline W$, so
that the boundary
    circles  $\Gamma_i$ are mapped to orbits $\gamma_i\in\Pc$
equivariantly  with respect
    to the canonical action of the circle $\R/\Z$ on $\Gamma_i$ and
the action of
    $\R/\Z$   on
    $\gamma_i$ which is  generated by the time $1$ map of the vector
field $T_i\Reeb$
    if
    $\gamma_i\in\Pc_-\;$ and of  $-T_i\Reeb$ if $\gamma_i\in\Pc_+$, where
      $T_i=S(\gamma_i)=\int\limits_{\gamma_i}\lambda $.\footnote{A
possibly confusing difference in signs
       here and in Lemma
      \ref{lm:asymptot} in the cylindrical case is
      caused by the fact that the action $\R/\Z$ on $\Gamma_i$ is
defined by the linear complex
      structure on the plane $T_{z_i}S$ tangent to $S$ at the
corresponding puncture.}
\end{proposition}
  Punctures associated with  orbits from $\Pc_+$ are called {\it
positive}, while punctures associated with  orbits from $\Pc_-$ are
called {\it negative}.
As in the contact cylindrical  case,  the conditions of finiteness of
the full-energy and the $\omega$--energy
are essentially equivalent.
The following statement is similar   to Proposition \ref{prop:automatic-bound}.
   \begin{proposition}\label{prop:automatic-bound-ends}
Let $(W,J)$ be an almost complex manifold
with  symmetric cylindrical ends  adjusted to a symplectic form
$\omega$.
   Suppose that a  holomorphic curve with punctures
$F\co (S\setminus Z,j)\to (W,J)$ is
  asymptotic at the positive punctures to the
periodic orbits $\og_1,\dots,\og_k\in\Pc_+$ and to the periodic orbits
$\ug_1,\dots,\ug_l \in \Pc_-$ at  the negative punctures.
Then there exists a positive constant $C$ (which depends on $J,\lambda,\omega$
but not $F$)
such that
\begin{equation}
E(F)\leq C\left(
\int\limits_{F^{-1}(\overline W)}F^*\omega+
\int\limits_{F^{-1}(E_+)}f_+^*\omega
+\int\limits_{F^{-1}(E_-)}f_-^*\omega\right)+
\sum\limits_1^k S(\og_i) + \sum\limits_1^l S(\ug_i).
\end{equation}
In particular, the energies
$E(F)$ are  uniformly bounded for all $F$ for which
$f$ represents a given homology class in
$H_2(V,\bigcup\limits_i\og_i\cup\bigcup\limits_j\ug_j)$.
\end{proposition}
\section{Holomorphic buildings in cylindrical manifolds\break
$W=\boldsymbol{\bR}\times V$}\label{buildings-cyl}

  \subsection{Holomorphic buildings of height $1$}\label{ssect2.3}
We first introduce in a more systematic way the types of
holomorphic curves needed for the  compactification of the moduli
spaces of holomorphic curves in a cylindrical manifold.
Let $(S,j, M\cup Z,D)$ be a nodal Riemann surface, such that
the set of its marked points is presented as a disjoint
union of two ordered sets $M$ and $Z$. The points
from  $Z$ are called {\it punctures}, the points from $M$  are
called {\it marked points}.  The
set
$$
D=\{\od_1,\ud_1,\od_2,\ud_2,\dots,\od_s, \ud_s\}.
$$
  of {\it special marked} points is viewed as an unordered set of unordered
  pairs.
  The
surface $S$ may be disconnected, and the points of any given
special pair $(\od_i,\ud_i)$ may belong to the same component, or  to different
components of $S$.
A holomorphic curve $F=(a, f)$ is called the {\it trivial} or {\it
vertical cylinder } if
  $S$ is the Riemann sphere, the
sets $M$ and $D$ are empty, the set $Z$ consists of exactly $2$ points
and $f$ maps $S$ onto a periodic orbit $\g$.
A {\it nodal holomorphic curve  (or building) of height  $1$} is  a
proper holomorphic map
$$F=(a,f)\co (S\setminus Z,D,M,j)\to (\bR\times V,J)$$
  of finite energy  which
  sends elements of each special pair
to one point: $$F(\od_i)=F(\ud_i)\quad\hbox{ for each}\quad
i=1,\dots, s.$$ The curve $F$ is called {\it stable}  in $\R\times V$, if the
following conditions are satisfied:
\begin{description}
\item{\bf Stab 1}\qua at least one connected component of the curve
is not a trivial
cylinder,
  \item{\bf Stab 2}\qua if $C$ is a connected component of $S$ and    the
map $F|_{C}$ is constant then the Riemann surface $C$ together
with all its marked, special marked points and punctures
   is stable in the sense
of  Section \ref{sec:Riemann-stable} above.
   \end{description}
         Lemma \ref{lm:asymptot} and Proposition \ref{prop:punct-cob}
describe the behavior of a holomorphic curve of height $1$ near each
puncture.
  In particular, one can associate a  periodic orbit  $\g_i\in\Pc$
   to each
  puncture $z_i\in Z$. The coordinate $a$ of the map $F$
  tends near each puncture either to $+\infty$ or $-\infty$.
   Respectively, we call the punctures positive or negative, and denote the
   set of positive resp. negative punctures by $\oZ$ resp. $\uZ$.
  The signature of a
   holomorphic curve of height $1$ is the quadruple of integers
$(g,\mu,p^+,p^-)$, where $g$
  is the arithmetic genus \eqref{eq:arithm-genus}  of $S$, where
$\mu=\#M$ is the number of marked points, and
  where $p^\pm$ are the numbers of positive respectively negative punctures.
     As in the case  of Riemann surfaces,  a holomorphic curve $F$ of
height $1$ is called
{\it  connected} if the singular Riemann surface
  $\widehat S_{D}$ is connected.
   Two nodal holomorphic curves,
    $$(F,S,j,M,Z,D)\;\;\hbox{ and}\;\;(F',S',j',M',Z',D'),$$ of height $1$
   are called {\it equivalent} if there exists a diffeomorphism
   $\varphi\co S\to S'$ such
   that
   \begin{itemize}
  \item  $\varphi_*j=j'$
  \item $f'\circ\varphi=f$, $a'\circ\varphi=a+\const$
  \item $\varphi$  sends  the ordered sets $M$ and $Z$ isomorphically
  to $M'$ and $Z'$.
  \item $\varphi|_D$ is an isomorphism $D\to D'$ of unordered sets of
unordered pairs.
  \end{itemize}
    In particular, {\it we identify curves which  differ by a
translation along the $\R$--factor.}
     If the
    curves $F$
    and $F'$ are connected, then  we can say equivalently  that we identify
the curves
    which have
    the same projections $f$ and $f'$ to the contact manifold $V$.
  The moduli space of  stable {\it connected }  smooth (ie, without double
  points)
holomorphic curves   of signature $(g,\mu,p^+,p^-)$
  is denoted by $\Mc_{g,\mu,p^+,p^-}(V)$. The bigger moduli space of
  stable {\it connected }  nodal
  holomorphic curves  of height $1$ and   of signature $(g,\mu,p^+,p^-)$
  will be  denoted by ${}^1\Mc_{g,\mu,p^+,p^-}(V)$.
  Unlike the case of Riemann surfaces,  the
  space ${}^1\Mc_{g,\mu,p^+,p^-}(V)$ is not large enough
  to compactify $\Mc_{g,\mu,p^+,p^-}(V)$.
   For this purpose  we need  holomorphic curves (buildings) of height $>1$
discussed in the next section.

\subsection{Holomorphic buildings of height $k$} \label{sec:height-k}
\begin{figure}[ht!]\small
\centering
\psfrag{z1}{$z_1$}
\psfrag{z2}{$z_2$}
\psfrag{z3}{$z_3$}
\psfrag{z4}{$z_4$}
\psfrag{z5}{$z_5$}
\psfrag{z6}{$z_6$}
\psfrag{z7}{$z_7$}
\psfrag{m1}{$M_1=\{z_3, z_4\}$}
\psfrag{m2}{$M_2=\{z_1, z_7\}$}
\psfrag{m3}{$M_3=\{z_2, z_5, z_6\}$}
\psfrag{m}{$M=M_1\cup M_2\cup M_3=\{z_1, z_2, \ldots, z_7\}$}
\psfrag{ph1}{$\Phi_1$}
\psfrag{ph2}{$\Phi_2$}
\includegraphics[width=0.5\textwidth]{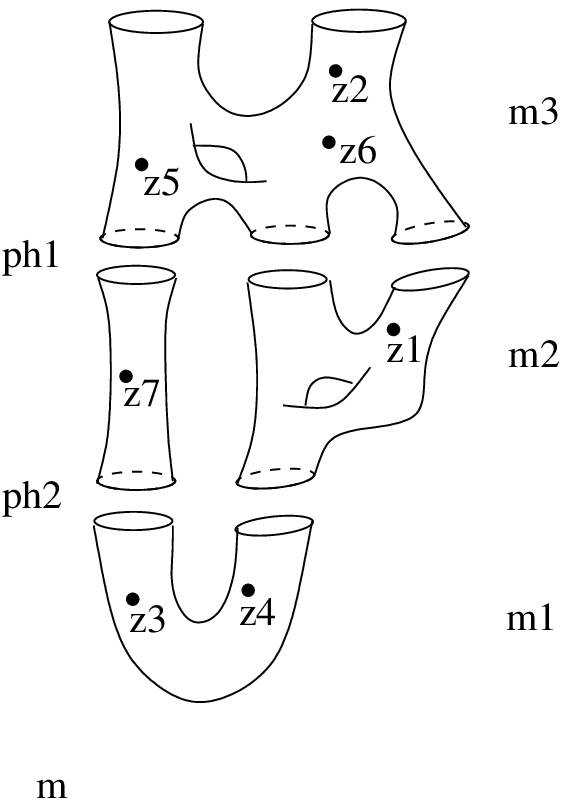}
\caption{Holomorphic building of height three with an ordered set of
marked points}\label{fig6}
\end{figure}

Suppose we are given $k$ stable, possibly disconnected nodal
curves of height $1$,
$$
F_m=(a_m,f_m;S_m,j_m,M_m,D_m, Z_m=\oZ_m\cup \uZ_m)\,,\; m=1,\dots,
k\,.
$$
Suppose,  in addition, that we are given a cross-ordering $\sigma$ of
$M=\mathop{\bigcup}\limits_{m=1}^{k} M_m$, which is compatible
with the ordering of each individual $M_i$, but may mix the points
of different $M_i$ in an arbitrary way.  See Figure \ref{fig6}, where
the  ordering of each $M_i,\,i=1,2,3,$
is  induced by the natural ordering of the index set.
Let
  $\wh S^{Z_m}_m$ be   the circle compactification  of  the Riemann
surface $(S_m,j_m)$ at punctures $Z_m$, as
  described in Sections
\ref{sec:blow-up} and \ref{sec:Riemann-stable} above. We denote by
$\Gamma_m^+$ and $\Gamma_m^-$ the sets of boundary circles which correspond to
the sets $\oZ_m$ and $\uZ_m$ of  punctures.
   Suppose  that for  each $m=1,\dots, k-1$
    the number $p_m^+$ of positive punctures of $F_m$
   is equal to the number $p_{m+1}^-$ of negative  punctures of $F_{m+1}$,
   and that there is given an orientation  reversing
   diffeomorphism  $\Phi_m\co \Gamma_m^+\to \Gamma_{m+1}^-$ which is an
   orthogonal map on each boundary component. Using these maps we can,
    similarly to the construction
   of the surface $S^{D,r}$ in Section \ref{sec:Riemann-stable},
form a piecewise
   smooth surface
   $$ S^{Z,\Phi}=S_1^{Z_1}\mathop{\cup}\limits_{\Phi_1} S_2^{Z_2}\mathop{\cup}
   \limits_{\Phi_2}\dots\mathop{\cup}\limits_{\Phi_{k-1}}S_k^{Z_k}\;.$$
The sequence  $F=\{F_1,\dots,F_k\}$ of holomorphic curves
    of height $1$, together with the
    decoration maps
    $\Phi=\{\Phi_1,\dots,\Phi_{k-1}\}$ and the  cross-ordering
$\sigma$ is called a {\it holomorphic building of height (or level)
$k$},  if the  compactified maps $\of_m\co S_m^{Z_m}\to V$ fit together
into a continuous map
   $\of\co S^{Z,\Phi}\to V$. This property implies, in particular, that
for $m=1,\dots, k-1$,
    the  curve $F_m$ at its positive punctures
   is asymptotic to the same periodic orbits as the curve $F_{m+1}$
    at its corresponding negative punctures.
   Two holomorphic buildings of height $k$, namely $(F,\Phi,\sigma)$
and $(F',\Phi',\sigma')$, where
   \begin{align*}
(F,\Phi,\sigma)&=\left(\{F_1,F_2,\dots,F_k\},\{\Phi_1,\dots,\Phi_{k-1}\},\sigma\right),\\
   F_i&=(a_i,f_i; S_i,,j_i,M_i,D_i,Z_i),\; \text{for $i=1,\dots, k$},\\
\intertext{and}
(F',\Phi',\sigma')&=\left(\{F'_1,F'_2,\dots,F'_k\},\{\Phi'_1,\dots,\Phi'_{k-1},\sigma'\}\right),\\
  F'_i&=(a'_i,f'_i; S'_i,,j_i',M'_i,D'_i,Z'_i),\;\text{for $i=1,\dots, k$},
\end{align*}
are called
   {\it equivalent}
   if there exists a sequence
   $\varphi=\{\varphi_1,\dots,\varphi_k\}$ of diffeomorphisms
   having the following properties,
   \begin{itemize}
   \item $\varphi_m$,  for $m=1,\dots, k$, is an equivalence between
the height $1$ holomorphic buildings
   \begin{eqnarray*}F_m&=&(a_m,f_m; S_m,j_m,M_m,D_m,Z_m)\;\;\hbox{ and }\\
   F'_m&=&(a'_m,f'_m; S'_m,j'_m,M'_m,D'_m,Z'_m),
     \end{eqnarray*}
   \item  $\varphi$ commutes with   the sequences $\Phi$ and $\Phi'$
   of attaching maps, ie,
   $$\Phi'_{m+1}\circ\varphi_m=\varphi_{m+1}\circ\Phi_{m},\; \text{for
$m=1,\dots,k-1$}\, ,$$
   \item $\varphi_*\sigma=\sigma'$.
   \end{itemize}
   Additionally, we identify holomorphic buildings which differ by a
synchronized re-ordering
   of  the pair of the sets  $\oZ_m$  and $\uZ_{m+1}$ for $1< k\leq m$.
To keep the notation simple we will usually drop the cross-ordering
$\sigma$
from the notation and will write
$(F,\Phi)$ for a holomorphic building of height $k$. The main
points to remember are the following. First of all the union of
marked points coming from the various levels is ordered. Secondly,
the union of the negative punctures on the first level and the
positive punctures on the  highest level are ordered. Thirdly, there
is a compatibility  between two consecutive levels in the sense
that asymptotic limits match (as specified by the decoration map
$\Phi$).
The genus $g$ of a height $k$ building $(F,\Phi,\sigma)$
    is by definition the arithmetic genus of $S^{Z,\Phi}$. Its signature
    is defined as the quadruple
    $(g,\mu,p^-,p^+)$, where  $\mu$ is the total cardinality of the set
    $M=\bigcup\limits_1^k M_i$ and $p^+=p_k^+$ and $p^-=p_1^-$.
\begin{figure}[ht!]\small
\centering
\psfrag{og1}{$\og_1$}
\psfrag{og2}{$\og_2$}
\psfrag{og4}{$\og_4$}
\psfrag{og7}{$\og_7$}
\psfrag{ug3}{$\ug_3$}
\psfrag{ug5}{$\ug_5$}
\psfrag{ug6}{$\ug_6$}
\psfrag{m1}{$m_1$}
\psfrag{m2}{$m_2$}
\psfrag{m3}{$m_3$}
\psfrag{m4}{$m_4$}
\psfrag{m5}{$m_5$}
\psfrag{m6}{$m_6$}
\psfrag{ph1}{$\Phi_1$}
\psfrag{ph2}{$\Phi_2$}
\psfrag{ph3}{$\Phi_3$}
\includegraphics[width=0.65\textwidth]{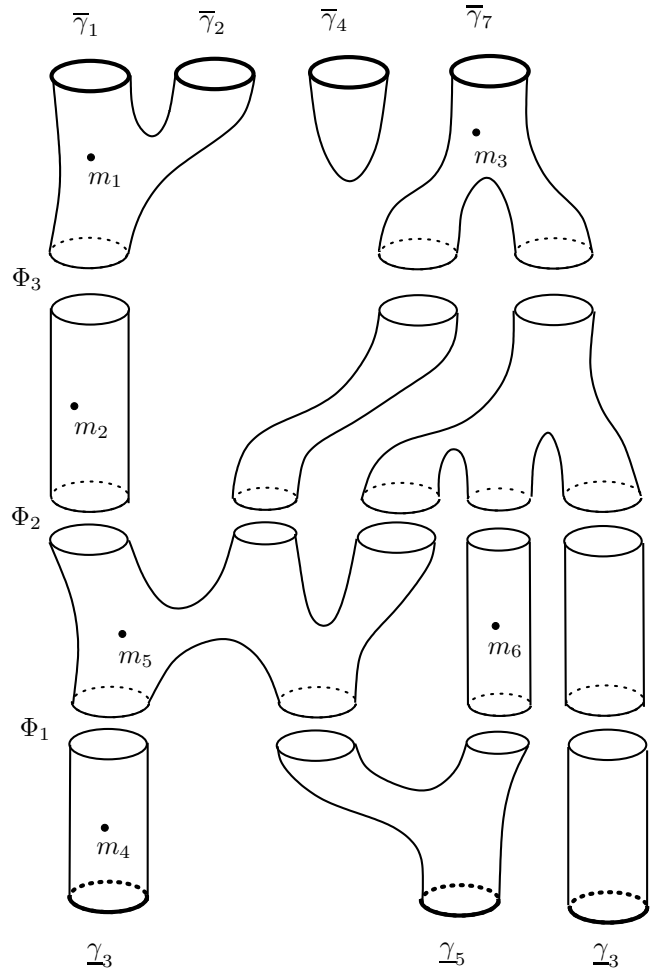}
\caption{Holomorphic building of height  $4$}\label{fig7}
\end{figure}
  Note that if
    $\varphi=\{\varphi_1,\dots,\varphi_k\}$ is an equivalence  between
    $(F,\Phi)$ and $(F',\Phi')$ then the homeomorphisms
         $\varphi_m\co S_m\to S'_m$, for $m=1,\dots,k-1$,  fit together
into a homeomorphism
   $\overline \varphi\co S^{Z,\Phi}\to( S')^{Z',\Phi'}$ between the two surfaces\footnote
   {The converse, however, is not true unless the asymptotic orbits
   associated with the punctures are {\it simple}.}
       $$S^{Z,\Phi}=S_1^{Z_1}\mathop{\cup}\limits_{\Phi_1}
S_2^{Z_2}\mathop{\cup}
   \limits_{\Phi_2}\dots\mathop{\cup}\limits_{\Phi_{k-1}}S_k^{Z_k}\;$$
   and $$(S')^{Z',\Phi'}=(S'_1)^{Z'_1}\mathop{\cup}\limits_{\Phi'_1}
(S'_2)^{Z'_2}\mathop{\cup}
   \limits_{\Phi'_2}\dots\mathop{\cup}\limits_{\Phi'_{k-1}}(S'_k)^{Z'_k}\;.$$
It is also useful to note that in the cases in which  all the curves
$F_m$, for $m=1,\dots, k$, are connected,
    or in which they are disconnected but
     have exactly one component different from a trivial cylinder,
     one can define the holomorphic building of height $k$ purely in
terms of their
    $V$--components $f_1,\dots, f_m$ of the maps $F_1,\dots,F_k$.
   The stability condition for $F$ means the stability of all its components
   $F_1,\dots,F_k$.
    The moduli space of equivalence classes of stable holomorphic
buildings of height $k$ and
   signature $(g,\mu,p^+,p^-)$ is denoted
   by ${}^k\Mc_{g,\mu,p^-,p^+}(V)$. We
   set $$\overline \Mc_{g,\mu,p^-,p^+}(V)=
   \bigcup\limits_{k=1}^\infty{}^k\Mc_{g,\mu,p^-,p^+}(V)$$ and
  $$\overline \Mc_{g,\mu}(V)=\bigcup\limits_{p^-,p^+\geq
0}\overline\Mc_{g,\mu,p^-,p^+}(V)\,.$$
   With each holomorphic building
    \begin{gather*}\label{eq:height-k}
(F,\Phi,\sigma)=\left(\{F_1,F_2,\dots,F_k\},\{\Phi_1,\dots,\Phi_{k-1}\},\sigma\right),\\
     F_m=(a_m,f_m;S_m,j_m,M_m,D_m, Z_m=\oZ_m\cup \uZ_m)\,,\; \text{for
$m=1,\dots, k$},
       \end{gather*}
       of height $k$
                we can associate the {\it underlying nodal Riemann  surface}
                 $$\bS_F=
    \left(\bigcup\limits_1^k (S_m,j_m), M=\bigcup
    \limits_1^kM_m\cup \uZ_1\cup \oZ_k,D=\bigcup
    \limits_1^k
    D_m\cup\bigcup\limits_1^{k-1}(\oZ_m\cup \uZ_{m+1})\right). $$
    Here we treat the punctures from
     $\uZ_1$ and $\oZ_k$ as extra marked points, while the set
$\oZ_1\cup \uZ_2\cup\dots\cup\oZ_{k-1}\cup \uZ_k$ as
    additional {\it special} marked points where the puncture $\oz\in
\oZ_i$ is coupled with the
    puncture $\uz\in  \uZ_{i+1}$, $i=1,\dots,k-1$, if the map $\Phi_i$ maps
    the compactifying circle $\Gamma_{\oz}$ associated with  the puncture
    $\oz$ onto the circle
    $\Gamma_{\uz}$ associated with $\uz$.  The ordering of $M$ is
given by $\sigma$, rather than
     the natural ordering of the  union of ordered
    sets $M_1,\dots,M_k$. The maps $\Phi_i$ define the decorations
     at these double points in the sense of Section
\ref{sec:Riemann-stable} above. Hence,
     $\bS_F$ is partially decorated, and
     in the case when each of the
     height  $1$ nodal curves $F_m$, for $m=1,\dots, k$, forming $F$,
is  equipped  with its
     own decoration $r_m$,
       the Riemann surface $\bS_F$ gets a full  decoration
$r_{F,\Phi}=\{r_1,\dots,r_k,\Phi_1,
     \dots,\Phi_k\}$.
     It is important to realize that the stability of the curve $F$
does not guarantee
     the stability of the Riemann surface $\bS_F$.
     However one can always add a few marked points to some of the sets
     $M_i$ in order to stabilize
     the Riemann nodal surface   $\bS'=\bS_{F'}$ which underlies
     the new holomorphic building $F'$.

    \subsection{Topology of  $\overline{\Mc}_{g,\mu,p_-,p_+}(V)$}
    \label{sec:topology-curves}

\noindent The notion of  convergence  in  $\overline\Mc_{g,\mu,p_-,p_+}(V)$
which we define below is compatible with the metric space structure on
    $\overline\Mc_{g,\mu,p_-,p_+}(V)$ defined in Appendix
\ref{ap:metric-cylindrical}.
    In particular, the topology introduced here is Hausdorff.
    Suppose that we are given a sequence
    $$(F_i,\Phi_i)\in \overline\Mc_{g,\mu,p_-,p_+}(V),\;\;
  \text{for $  i \geq 1$}, $$  of   holomorphic buildings of height $\leq k$.
      The sequence  $(F_i,\Phi_i)$
         converges to a building\ $(F,\Phi)\in
\overline\Mc_{g,\mu,p_-,p_+}(V)$
    of height $k$
         if  there exist sequences $M'_i$ of  extra sets of marked points
         for the buildings  $(F_i,\Phi_i)$  and a set $M'$ of extra
         marked points  for  the building $(F,\Phi)$,  which have the
same cardinality $N$ and which
          stabilize the
         corresponding underlying Riemann surfaces, and
         such that the following conditions  are satisfied.
          Let $$(\bS_{F_i},r_{F_i,\Phi_i})   =(S_i,j_i,M_i\cup
M'_i,D_i,r_i)$$ and
    $$(\bS_F,r_{F,\Phi})=(S,j,M\cup M',D,r)$$
    be the decorated stable nodal Riemann surfaces underlying
     the curves $$(F_m,\Phi_m)\;\;\hbox{ and}\;\;(F,\Phi)$$ with extra
marked points.
      Then there exist
     {diffeomorphisms} $\varphi_i\co  S^{D,r}\to
    S^{D_i,r_i}$ with $\varphi_i(M)=M_i$ and  $\varphi_i(M')=M'_i$
    which satisfy the conditions CRS1--CRS3 in the  definition
    of convergence of Riemann surfaces and, in addition,  the
following conditions.
  \begin{description}
      \item{\bf CHC1}\qua  The sequence of the compactified projections
    $\overline f_i\circ \varphi_i\co  S^{D,r}\to V$ converges to
$\overline f\co S^{D,r}\to V$  uniformly.
    \item{\bf CHC2}\qua Let us denote by $C_l$ the union of components
    of $S^{D,r}\setminus\bigcup\Gamma_m$ which correspond
       to the same level $l=1,\dots, k$ of the building $F$.
    Then there exist sequences of real numbers $c_i^l$, for
$l=1,\dots, k$ and  $i\geq 1,$ such that
      $(a_i\circ\varphi_i-a-c^l_i)|_{C_l}\to 0$ in the
$C^0_{\mathrm{loc}}$--topology.
    \end{description}
The $C^0_{\text{loc}}$--convergence in  CHC2 can be equivalently replaced  by
   the $C^{\infty}_{\text{loc}}$--conver\-g\-ence in view of the elliptic
regularity theory.
\begin{figure}[ht!]\small
\centering
\includegraphics[width=0.6\textwidth]{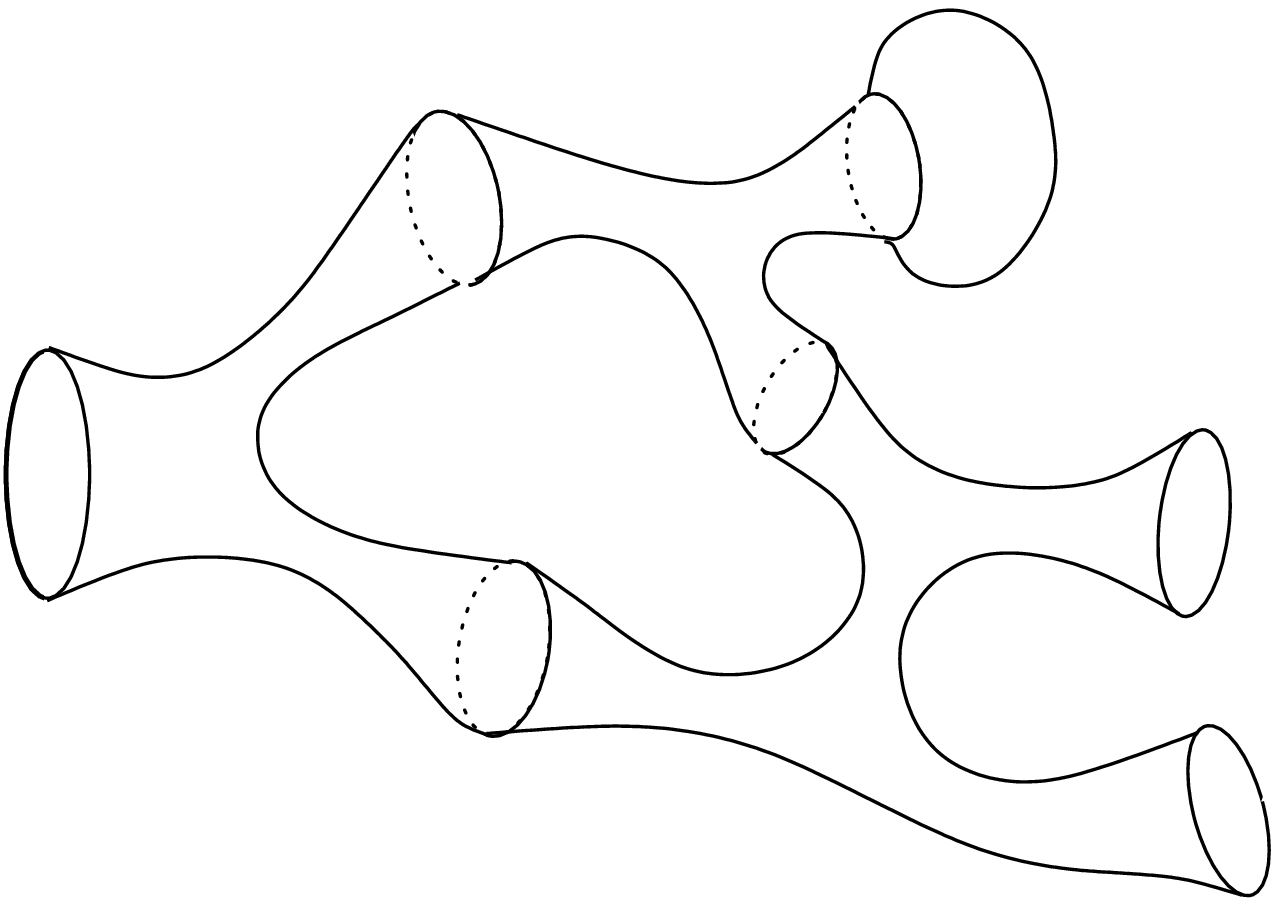}
\caption{The map $\overline{f}\co S^{D,r}\to V$}\label{fig8}
\end{figure}
  \section{Holomorphic buildings   in
manifolds with cylindrical ends}\label{sec:buildings-split}
\subsection{Holomorphic buildings  of height $k_-|1|k_+$}
We now generalize the notion of a holomorphic building to the
case,  in which  the target manifold has cylindrical ends, rather
than being cylindrical.
A nodal holomorphic curve in $W$, or a {\it holomorphic building of
height $1$} is  defined similarly to a nodal
curve of height $1$ in a cylindrical manifold,
  ie, it is a  proper holomorphic map
$$F=(S\setminus Z,j,D,M)\to ({W},J)$$
  of finite energy  which
  sends elements of every  special pair
to one point, ie, $F(\od_i)=F(\ud_i)$ for $i=1,\dots, s$.
Suppose that we are given:
\begin{enumerate}
\item[(i)] A holomorphic building of  height  $k_+$  in the
cylindrical manifold
$\R\times V_+$:
  $$(F^+,\Phi^+)=\left(\{F_1,F_2,\dots,F_{k_+}\},\{\Phi_1,\dots,
  \Phi_{k_+-1}\}\right)$$
  $$F_i=(a_i,f_i;S_i,j_i,M_i,D_i,\uZ_i\cup \oZ_i)\,, \;  \text{for
$i=1,\dots,k_+$}.$$
  \item[(ii)] A  holomorphic building of   height $k_-$
  in the cylindrical manifold
$\R\times V_-$:
  $$(F^-,\Phi^-)=\left(\{F_{-k_-},F_{-k_-+1},\dots,F_{-1}\},\{\Phi_{-k_-},\dots,
  \Phi_{-2}\}\right)$$
   $$F_i=(a_i,f_i;S_i,j_i,M_i,D_i,\uZ_i\cup \oZ_i),\; \text{for $
i=-k_-,\dots,-1$}.$$
\item[(iii)] A
nodal holomorphic  curve $(F_0,S_0,D_0,M_0, \uZ_0\cup\oZ_0,j_0)$ in
$({W},J)$. We denote by
$\Gamma_0^\pm$ the sets of boundary circles which correspond to the
punctures $\uZ_0$ and $\oZ_0$.
\item[(iv)]  An ordering of $\bigcup\limits_{k_-}^{k_+}M_i$ which is
compatible with the ordering of each individual $M_i$ but not
necessarily respecting the numbering
of  the sets $M_{k_-},\dots,M_{k_+}$.
\end{enumerate}
  Suppose
   that
   \begin{itemize}
   \item
   the number $p_0^+$ of positive punctures of $F_0$
   is equal to the number $p_{1}^-$ of negative  punctures of $F_{1}$,
   \item  the number $p_{-1}^+$ of positive punctures of $F_{-1}$
   is equal to the number $p_{0}^-$ of negative  punctures of $F_{0}$,
   \item for $m=-1,0$
   there is given an orientation  reversing
   diffeomorphism  $\Phi_m\co \Gamma_m^+\to \Gamma_{m+1}^-$ which is
   orthogonal  on each boundary component.
   \end{itemize}
\begin{figure}[ht!]\small
\centering \psfrag{rvp}{${\mathbb R}\times V_{+}$}
\psfrag{rvm}{${\mathbb R}\times V_{-}$} \psfrag{w}{$W$}
\includegraphics[width=0.35\textwidth]{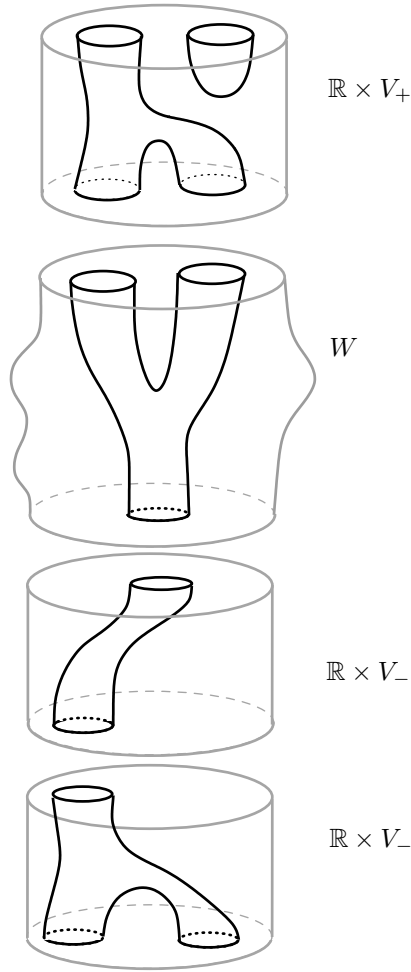}
\caption{Holomorphic building of height  $2\vert 1\vert
1$}\label{fig9}
\end{figure}
   Let $S^0=S_0^{Z_0}$ be the oriented blow-up of $S_0$ at the punctures
   from $Z_0$, and let the surfaces
      \begin{equation}
        \begin{split}
        S^+=S^{Z^+,\Phi^+}=&S_1^{Z_1}\mathop{\cup}\limits_{\Phi_1}
S_2^{Z_2}\mathop{\cup}
\limits_{\Phi_2}\dots\mathop{\cup}\limits_{\Phi_{k_+-1}}S_{k_+}^{Z_{k_+}}\;
\cr
S^-=S^{Z^-,\Phi^-}=&S_{-k_-}^{Z_{-k_-}}\mathop{\cup}\limits_{\Phi_{-k_-}}
S_{-k_-+1}^{Z_{-k_-+1}}
   \mathop{\cup}
\limits_{\Phi_{-k_-+1}}\dots\mathop{\cup}\limits_{\Phi_{-2}}S_{-1}^{Z_{-1}}\;\cr
   \end{split}
   \end{equation}
be defined as in Section \ref{sec:height-k}.
   Gluing $S^-$ and $S^0$ by means of  $\Phi_{-1}$, and $S^0$ and
$S^+$ by means of  $\Phi_0$
   we obtain the  piecewise-smooth surface
    \begin{equation}\label{eq:underlying}
    \overline S= S^-\mathop{\cup}\limits_{\Phi_{-1}} S^0\mathop{\cup}
   \limits_{\Phi_0}S^+\,.
   \end{equation}
  The last condition in the definition
  of a building in $ W$ of height $k_-|1|k_+$ can now be formulated as follows:
     \begin{enumerate}
\item[(v)] for a sufficiently small $\d>0$ the
  maps $$\overline f^-\co S^-\to V_-\,,\;\;G^\d\circ F^0\co S^0\to \circW,\;\;\hbox{
   and  }\;\;
   \overline f^+\co S^+\to V_+$$ fit together into a continuous map
   ${\overline F}\co \overline S\to \overline W$.
\end{enumerate}
We will also say sometimes that a holomorphic building of height
$k_-|1|k_+$ consists of
$3$ {\it layers}, namely the {\it lower} layer is a  holomorphic
building $F_-$ of height $k_-$,
  the {\it main} layer is a holomorphic curve $F_0$
( of height $1$), and the {\it upper}  layer is a holomorphic
building  $F_+$ of height $k_+$.
   The equivalence relation for holomorphic buildings of height
$k_-|1|k_+$ is defined similarly to buildings in cylindrical
manifolds except that  there is no translation to be quotient out  in
the central layer.
\begin{figure}[ht!]\small
\centering \psfrag{fp}{$\overline{f}^{+}$}
\psfrag{fm}{$\overline{f}^{-}$} \psfrag{ow}{$\overline{W}$}
\psfrag{gf}{$G^{\delta}\circ F^{\circ}$}
\includegraphics[width=0.5\textwidth]{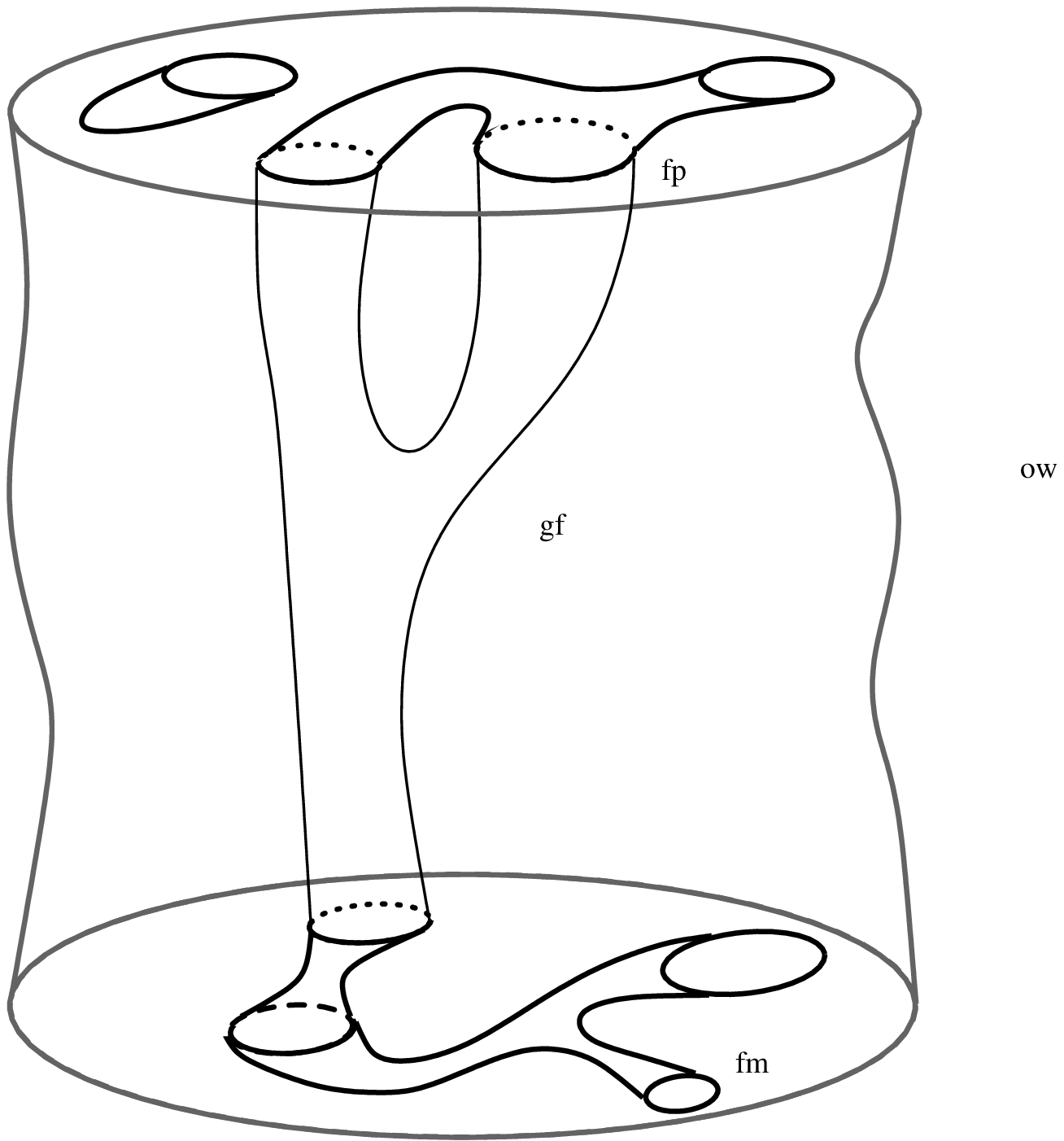}
\caption{The maps $\overline{f}^{-}$, $G^{\delta}\circ F^{\circ}$
and $\overline{f}^{+}$ for a continuous map
$\overline{F}\co \overline{S}\to \overline{W}$}\label{fig10}
\end{figure}
   As in the cylindrical case the genus $g$ of the building
$(F,\varphi)$ of height  $k_-|1|k_+$
    is by definition the arithmetic genus of $\overline S$. Its signature
    is defined as the quadruple
    $(g,\mu,p^-,p^+)$, where  $\mu$ is the total cardinality of the set
    $M=\bigcup_{k_-}^{k_+} M_i$,  and where $p^+=p^+_{k_+}$ and
$p^-=p^-_{-k_-}$.
The energies  $E(F)$, $E_\l(F)$ and $E_\omega(F)$  of a curve
  $F$ of height  $k_-|1|k_+$ are defined by the formulas
\begin{equation}
\begin{split}
E_\omega(F) =& \sum\limits_{-k_-}^{k_+} E_\omega (F_i)\,\cr
  E_\l(F) =&
\max\limits_{-k_-\leq i\leq k_+} E_\l(F_i)\, \cr
E(F)=&E_\l(F)+E_\omega(F)\,.\cr
\end{split}
  \end{equation}
A holomorphic building  $(F^-,F^0,F^+)$  of height  $k_+|1|k_-$
  is called {\it stable} if all its three  layers,
  \begin{itemize}
  \item the height  $k_+$  building
in $\R\times V_+$,
\item  the height
$k_-$ building $F^-$ in $\R\times V_-$  and
\item the curve $F^0$ in $ W$,
\end{itemize}
  are stable. Here the map $F_0$ into $W$ is called stable if for
every component
   $C$ of the underlying Riemann surface $S_0$ either the restriction
   $F_0\vert_{C}$ is non constant or  in case $F_0\vert_{C}$ is constant,
    then $C$, equipped with all its marked points from $D_0\cup M_0$
and all its punctures
    from $Z_0$, is stable. Equivalently, this requires that the
automorphism group of every
    component of $S_0$ equipped with all its distinguished points, is finite.
  The moduli spaces of holomorphic buildings of
  signature $(g,\mu,p^-,p^+)$ and  height $k_-|1|k_+$ in $( W,J)$
  is denoted   by ${}^{k_-,k_+}\Mc_{g,\mu,p^-,p^+}(W,J)$.
   We set
     \begin{equation}
     \begin{split}
     {}^{k_-,k_+}\Mc_{g,\mu}(W,J)&=\bigcup\limits_{p^-,p^+\geq 0}
      {}^{k_-,k_+}\Mc_{g,\mu,p^-,p^+}(W,J)\cr
      \overline\Mc_{g,\mu,p^-,p^+}(W,J)&=\bigcup\limits_{k_-,k_+\geq 0}
      {}^{k_-,k_+}\Mc_{g,\mu,p^-,p^+}(W,J)\cr
      \overline \Mc_{g,\mu}(W,J)&=\bigcup\limits_{p^-,p^+\geq 0}
      \overline\Mc_{g,\mu,p^-,p^+}(W,J)\ .\cr
      \end{split}
      \end{equation}

\subsection{Topology of  $\overline{\Mc}_{g,\mu}(W,J)$}

\noindent         In  this section we
              will   spell out the meaning of convergence
              of a sequence of  smooth  curves   $$F^{(k)}\in \Mc_{g,\mu}(W,J)=
                {}^{k_{-},k_{+}}\Mc_{g,\mu}(W,J),\quad\hbox{
with}\quad k_{-}=0=k_{+},$$
  for $k\geq 1$,  to
              a building
$$F=(\{F_{-k_-},\dots,F_0,\dots,F_{k_+}\};\{\Phi_{-k_-},\dots,\Phi_0,\dots,\Phi_{k_+}\})
$$
              from
              ${}^{k_-,k_+}\Mc_{g,\mu,p^-,p^+}(W,J).$
              A more general definition of convergence in the case in
which  $F^{(k)}$ is a
              sequence of holomorphic buildings from
               $${}^{k_-,k_+}\overline{\Mc}_{g,\mu,p^-,p^+}(W,J)=
               \mathop{\bigcup\limits_{0\leq i\leq k_-}}\limits_{0\leq
j\leq k_+}
              {}^{i,j}\Mc_{g,\mu,p^-,p^+}(W,J)$$ is similar and left
to the reader.
   The sequence $F^{(k)}$ converges
   to $F$
  if  there exist sequences $M^{(k)}$ of  extra sets of marked points
         for the curves  $F^{(k)}$  and a set $M$ of extra
         marked points  for  the building $F$  which have the same
cardinality $N$ and which
          stabilize the
         corresponding underlying Riemann surfaces
         such that the following conditions  are satisfied.
           Let $\bS_k=(S^{(k)},j^{(k)}, M^{(k)})$
              be Riemann surfaces underlying $F^{(k)}$  with the extra
marked points $M^{(k)}$, and
              $(\bS,\Phi)$ be the decorated Riemann surface underlying
the building $F$ with the extra set
              $M$ of marked points.
              We consider, as in (\ref{eq:underlying}), the surface
              $$ \overline S= S^-\mathop{\cup}\limits_{\Phi_{-1}}
S^0\mathop{\cup}
   \limits_{\Phi_0}S^+=
   S_{-k_-}^{Z_{-k_-}}\mathop{\cup}\limits_{\Phi_{-k_-}} S_{-k_-+1}^{Z_{-k_-+1}}
\mathop{\cup}\limits_{\Phi_{-k_-+1}}\dots\mathop{\cup}\limits_{\Phi_{k_+-1}}S_{k_+}^{Z_{k_+}}\,.$$
       with a  conformal structure
   $j$ which is degenerate along the union $\Gamma$ of special circles.
   Let $\overline F\co \overline S\to\overline W$ be the map described in
part (v) of the definition of a
    holomorphic building
   of height $k_-|1|k_+$.
   We also abbreviate  $$\dot{S}_i:= (\overline S\setminus\Gamma)\cap
S_i\; \text{for  $ i=-k_-,\dots,k_+\,.$}$$
          Suppose   that there exists a sequence of
         diffeomorphisms
   $\varphi_k\co \overline S\to S^{(k)}$ which satisfies the conditions
CRS1--CRS3 in the  definition
    of the convergence of decorated Riemann surfaces in section
\ref{sec:DM-topology}, and require,
     in addition, the following conditions.
   \begin{description}
    \item {\bf CHCE1}\qua For sufficiently large $k\geq K$,  the images
$F^{(k)}\circ\varphi_k|_{\dot{S}_i}$ for
     $i=-k_-,\dots, -1$, are
   contained in the cylindrical end $E_-$, and the
   the images $F^{(k)}\circ\varphi_k|_{\dot{S}_i}$ for  $i=1,\dots, k_+$, are
   contained in the cylindrical end $E_+$ of the manifold $W$.
   \item{\bf CHCE2}\qua There exist constants $c^{(k)}_i$ for
   $i=-k_-,\dots, -1,1,\dots, k_+$ and $k\geq K$, such that $\wt
   F^{(k)}_i\circ\varphi_k|_{\dot{S_i}}$ converge to $F_i$ uniformly on compact
   sets, where $F^{(k)}_i=(a^{(k)}_i,f^{(k)}_i)$ and $\wt
   F^{(k)}_i=(a^{(k)}_i+c^{(k)}_i,f^{(k)}_i)$.
   \item{\bf CHCE3}\qua The sequence $G^\d\circ
F^{(k)}\circ\varphi_k\co \overline S\to\overline W$ converges
   uniformly  to $\overline F$.
   \end{description}
    Note that the space $\overline \Mc_{g,\mu}(W,J)$ can be
             metrized  similar to the way it is done
               in Appendix
              \ref{ap:metric-cylindrical} below for the moduli spaces
of holomorphic
              buildings in cylindrical manifolds.
              For different values of $p^\pm$  the      spaces $
\overline\Mc_{g,\mu,p^-,p^+}(W,J)$
              are disjoint open--closed subsets of $\overline \Mc_{g,\mu}(W,J)$.

\section{Holomorphic buildings in split almost complex\break manifolds}
\subsection{Holomorphic buildings of height $\mathop{\vee}\limits_1^{k_0}$}

\noindent Let us recall the splitting construction from Section
\ref{sec:split}. We begin with a closed almost complex manifold
$(W,J)$, cut it open along a co-oriented  hypersurface $V$ to get a
manifold $\circW$ with two new  boundary
components $V',V''$ diffeomorphic to $V$, and attach to $V'$ and
$V''$ cylindrical ends,
thus obtaining a manifold  with cylindrical ends
of the
form $$\wt W=(-\infty,0] \times V \mathop{\cup}\limits_{0\times V=V'}
\circW \mathop{\cup}\limits_{V''=0\times V} [0,+\infty)
  \times V\,.$$
Of course, the manifold $\wt W$ is diffeomorphic to $W\setminus V$.
  The almost complex structure $J$  canonically extends
  to $\wt W$ as an almost complex structure $\wt J$ which is translation
  invariant on the ends.
A  holomorphic building $(F,\Phi)$  of height
$\mathop{\vee}\limits_1^{k_0}$ in the split manifold $(\wt W, \wt J)$
  is determined by the following data:
\begin{enumerate}
\item[(i)] a height 1  holomorphic curve (building) $F_0$ in $(\wt W, \wt J)$;
\item[(ii)] a height $k_0$ holomorphic building
$$(F',\Phi')=(\{F_1,\ldots,F_{k_0}\},
\{\Phi_1, \ldots,\Phi_{k_0-1}\})\,,$$
$$F_i=(a_i,f_i;S_i,j_i,M_i,D_i,\uZ_i\cup \oZ_i)\,, \; \text{for $
i=1,\dots,k_0,$}$$
in $(\RR \times V,J)$;
\item[(iii)] orientation reversing diffeomorphisms $\Phi_0 \co  \Gamma_0^+
\to \Gamma_1^-$ and $\Phi_{k_0} \co  \Gamma_{k_0}^+ \to \Gamma_0^-$,
orthogonal on each
boundary component;
\item[(iv)] an ordering of $\bigcup\limits_{0}^{k_0}M_i$ which is
compatible with the ordering of each individual $M_i$ but not
necessarily respecting the numbering
of   the sets $M_{0},\dots,M_{k_0}$.
\end{enumerate}
\begin{figure}[ht!]\small
\centering
\psfrag{wo}{$\circW$}
  \psfrag{rv}{${\mathbb R}\times V$}
   \psfrag{ph0}{$\Phi_0$}
   \psfrag{ph1}{$\Phi_1$}
  \psfrag{ph2}{$\Phi_2$}
\includegraphics[width=0.35\textwidth]{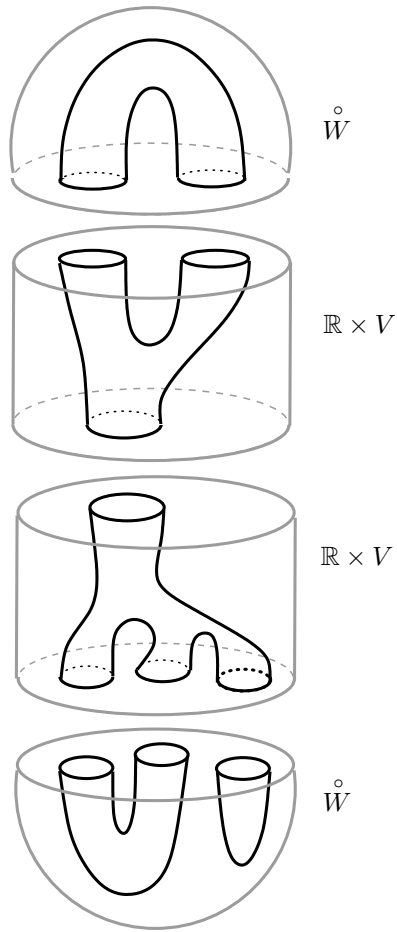}
\caption{Holomorphic buildings of height
$\mathop{\vee}\limits_1^{2}$.   The parts in the two pieces of  the
main layer are parameterized by $S_0$,
the  two levels  of the insert layer are parameterized by $S_1$ and
$S_2$} \label{fig11}
\end{figure}
Using blow-up at the punctures and
identification of the boundary components by means of the mappings $\Phi$,
we define the following surfaces:
\begin{eqnarray*}
S^0 &=& S^{Z_0}_0 \\
S' &=& S_1^{Z_1} \mathop{\cup}\limits_{\Phi_1} S_2^{Z_2} \mathop{\cup}
\limits_{\Phi_2} \ldots \mathop{\cup}\limits_{\Phi_{k_0-1}} S_{k_0}^{Z_{k_0}}
\\
\overline{S} &=& S' \mathop{\cup}\limits_{\Phi_0,\Phi_{k_0}} S^0.
\end{eqnarray*}
The last requirement in the definition of the holomorphic
building of height $\mvee\limits_1^{k_0}$ is  the following:
\begin{enumerate}
\item[(v)] the maps $\overline{f'} \co  S' \to V$ and $G^\delta \circ
F_0 \co  S_0 \to \circW$   fit together into a continuous map $\overline{F} \co 
\overline{S} \to W$.
\end{enumerate}
The holomorphic building  has  two layers: the main layer $F_0$, and
the {\it insert}
  layer $(F',\Phi')$.
  The equivalence relation for holomorphic buildings of height
$k_-|1|k_+$ is defined similarly to buildings in cylindrical
manifolds except that  there is no translation to be quotient out in
the insert layer.
The holomorphic building $(F,\Phi)$ is said to be stable if  all its
layers are stable.
The moduli spaces of holomorphic buildings  of
signature $(g,\mu)$ and  height $\mvee\limits_1^{k_0}$ in the split
manifold $(\wt W,\wt J)$
is denoted by ${}^{k_0}\Mc_{g,\mu}(\wt W,\wt J)$.
We set, recalling section \ref{sec:split},
\begin{equation}
\begin{split}
\overline\Mc_{g,\mu}(\wt W,\wt J)&=\bigcup\limits_{k_0 \geq 0}
{}^{k_0}\Mc_{g,\mu}(\wt W,\wt J)\cr
\Mc_{g,\mu}(W^{[0,\infty)},J^{[0,\infty)}) &= \bigcup\limits_{\tau
\in [0,\infty)} \overline \Mc_{g,\mu}(W^\tau,J^\tau)\cr
\overline \Mc_{g,\mu}(W^{[0,\infty]},J^{[0,\infty]})&=
\Mc_{g,\mu}(W^{[0,\infty)},J^{[0,\infty)}) \cup
\overline\Mc_{g,\mu}(\wt W,\wt J).
\end{split}
\end{equation}
The space $\overline \Mc_{g,\mu}(W^{[0,\infty]},J^{[0,\infty]})$ can be
topologized by introducing a metric similar to the way it is  done in
Appendix \ref{ap:metric-cylindrical}
  for the case of holomorphic buildings in cylindrical manifolds.
   The formula for the distance between
a holomorphic curve $(F,\Phi)$ in $(W^\tau,J^\tau)$ and a holomorphic curve
$(F',\Phi')$ in $(W^{\tau'},J^{\tau'})$ must contain the additional term
$\Abs{\frac1{1+\tau} - \frac1{1+\tau'}}$.
Let us  spell out the meaning of the convergence in
this topology,  of a sequence $F^{(k)}$ of stable
holomorphic curves into a  sequence of almost complex manifolds $(W^k,J^k)$
degenerating into the split almost complex manifold $(\wt W,\wt J)$.
We say that a sequence of stable holomorphic curves
$F^{(k)}$ into  $(W^k,J^k)$ converges to a stable level $k_0$
holomorphic building $(F,\Phi)$ of height $\mvee\limits_1^{k_0}$
in the split manifold $(\wt W,\wt J)$ if there exist
\begin{itemize}
\item extra sets of marked points $M^{(k)}$ and $M$ of the same
cardinality, which stabilize the
Riemann surfaces which underly the holomorphic curves $F^{(k)}$ and
the holomorphic building $(F,\Phi)$,
\item a sequence of diffeomorphisms
$\varphi_k \co  \overline{S} \rightarrow S^{(k)}$,
\item  sequences
$c_i^{(k)} \in \RR$,  for $i = 1, \ldots, k_0$,
\end{itemize}
  such that the   conditions CRS1--CRS3 of the Deligne--Mumford
convergence in  Section \ref{sec:DM-topology} are satisfied and such
that, in addition, the following conditions are met,
  \begin{description}
\item{\bf CHCS1}\qua $F_0^{(k)} \circ \varphi_k |_{\dot{S}_0}$ converges to $F_0$
uniformly on compact sets,
\item{\bf CHCS2}\qua for sufficiently large $k \ge K$, the images $F^{(k)} \circ
\varphi_k |_{\dot{S}_i}$  for $i = 1, \ldots, k_0$, are contained in the
cylindrical portion $[-k,k] \times V$ of $W^k$,
\item{\bf CHCS3}\qua  $(c_i^{(k)} + a_k \circ \varphi_k,
f_k \circ \varphi_k)|_{\dot{S}_i}$ converge uniformly on compact sets to
$F_i$ for $i = 1, \ldots, k_0$.
\end{description}
\subsection{Energy bounds for holomorphic curves in the process of splitting}
Suppose now  that $W$ is endowed with a symplectic structure $\omega$
compatible with $J$, and that the splitting
along $V$ is adjusted to $\omega$.
  The {\it energy}, and $\omega$--{\it energy} of a  holomorphic
building of height $\mvee\limits_1^{k_0}$ in  a split almost complex manifold
is naturally defined by the formulas
\begin{equation}
\begin{split}
E_\omega(F) =&  E_\omega (F_0)+E_\omega (F') \cr
  E_\l(F) =&
\max( E_\l(F_0),E_\l(F')) \cr
E(F)=&E_\l(F)+E_\omega(F)\,.\cr
\end{split}
  \end{equation}
So far we never specified the symplectic forms on  the family of
almost complex manifolds $(W^\tau,J^\tau)$ converging
to the split manifold $(\wt W,\wt J)$. This can be done but only in
such a way that in the limit the
  symplectic structure either
degenerates  or blows up. Instead, we slightly modify the notion of
energies for holomorphic curves in
  $(W^\tau,J^\tau)$.
  Let us recall that
  $$W^\tau=\circW\cup I_\tau,$$ where
  $I_\tau=[-\tau,\tau]\times V\,.$ Given a holomorphic curve
  $F\co (S,j)\to(W^\tau,J^\tau)$ we define its $\omega$--energy as
$$E_\omega(F)=\int\limits_{F^{-1}(\circW)}F^*\omega+\int\limits_{F^{-1}(I_\tau)}F^*p_V^*\omega\,,$$
  where $p_V$ is the projection $I_\tau=[-\tau,\tau]\times V\to V$.
  We also define
  $$E_\l(F)=\sup\int\limits_{F(S)\cap I_\tau}(\phi\circ p_{\R}\circ
F)dt\wedge \l,$$
  where $p_{\R}$ is the projection $I_\tau=[-\tau,\tau]\times V\to
[-\tau,\tau]$,
  and the supremum is taken over all function $\phi\co [-\tau,\tau]\to \R_+$ with
  $\int_{[-\tau ,\tau]}\phi(t)dt=1$. Finally, we set
  $$E(F)=E_\l(F)+E_\omega(F)\,.$$
With these definitions we immediately get
\begin{lemma}
Given a sequence of holomorphic curves $F^{(k)}$ in $(W^k,J^k)$ which
converges to a holomorphic  building
$F$ in the split manifold $(\wt W,\wt J)$, then
$$\lim\limits_{k\to \infty} E_\omega(F^{(k)})=E_\omega(F)\,.$$
\end{lemma}
  It turns out that a  uniform bound on the $\omega$--energy
automatically implies a uniform bound on the full energy.
  \begin{lemma}\label{lm:split-bound}
  There exists a constant $C$ which depends only on $(W,J)$, $V$ and
$\l$ such that for
   every $\tau>0$ and every  holomorphic curve $F\co (S,j)\to(W^\tau,J^\tau)$,
    $$E(F)\leq CE_\omega(F)\,.$$
  \end{lemma}
\proof
  Let us  denote  $V_+=V\times\tau\subset
  W^\tau=\circW\cup [-\tau,\tau]\times V$.
We will show that there exists a constant $K$ such that for any
$\tau>0$ and any
holomorphic curve $F\co (S,j)\to(W^\tau,J^\tau)$ we have
\begin{equation}\label{eq:bound-lambda}
\Big|\int\limits_{F^{-1}(V_+)}\l\Big|\leq KE_\omega(F)\,.
\end{equation}
Take a {sufficiently small} tubular neighborhood
$U_\e=V_+\times[0,1]$ of $V_+$ inside $\circW$, so that
$V_+=V_+\times 0$, and pull-back $\l$ to $U_\e$ via the projection
$V_+\times[0,\e]\to V_+=V$. Let $t\in[0,1]$ denote a coordinate in
$U_\e$ which corresponds to the second factor and $V_+'=V_+\times
1$. Applying Stokes' theorem we find
\begin{equation*}
\begin{split}
\int\limits_{F^{-1}(U_\e)}F^*d(t\l)&=\int\limits_{F^{-1}(V_+')}F^*\l,\cr
\int\limits_{F^{-1}(U_\e)}F^*d\l&=\int\limits_{F^{-1}(V'_+)}F^*\l-\int\limits_{F^{-1}(V_+)}F^*\l\,.\cr
\end{split}
\end{equation*}
On the other hand, taking into account that $F$ is $J$--holomorphic
and that $\omega$ is compatible with $J$
on $\circW$,  we have
\begin{equation*}
\begin{split}
\Big|\int\limits_{F^{-1}(U_\e)}F^*d(t\l)\Big|&\leq
K_1\int\limits_{F^{-1}(U_\e)}F^*\omega\,,\cr
\Big|\int\limits_{F^{-1}(U_\e)}F^*d\l\Big|&\leq
K_2\int\limits_{F^{-1}(U_\e)}F^*\omega\, . \cr
\end{split}
\end{equation*}
Consequently,
$$
\Big|\int\limits_{F^{-1}(V_+)}F^*\l\Big|
\leq (K_1+K_2)\int\limits_{F^{-1}(U_\e)}F^*\omega\leq KE_\omega(F)\,.
$$
Similarly, with the obvious notation,
$$
\Big|\int\limits_{F^{-1}(V_-)}F^*\l\Big|\leq KE_\omega(F).
$$
{The rest of the proof follows the lines of
Proposition \ref{prop:automatic-bound}.} \qed
\section{Compactness theorems}\label{sec:compact}
\subsection{Statement of main theorems}
In this section we  prove  the  main  results of the paper.
\begin{theorem}\label{thm:compact-cyl}
Let $(\RR \times V,J)$ be a symmetric cylindrical almost complex manifold.
Suppose that  the almost complex structure $J$ is adjusted to
the taming symplectic form $\omega$. Then
for every $E>0$, the space
$\overline\Mc_{g,\mu}(V)\cap \{E(F)\leq E\}$ is compact.
\end{theorem}
\begin{theorem}\label{thm:compact-ends}
Let $(W=E_-\cup\overline W\cup E_+,J)$
be an almost complex manifold with symmetric cylindrical ends. Suppose that
$J$ is adjusted to a symplectic form $\omega$ on $W$. Then
for every $E>0$, the space
$\overline\Mc_{g,\mu}(W,J)\cap \{E(F)\leq E\}$ is compact.
\end{theorem}
\begin{theorem}\label{thm:compact-split}
Let $(\wt W, \wt J)$  be a split almost complex manifold which is
obtained, as in Section
\ref{sec:split} above, by splitting a closed almost complex manifold
$(W,J)$ along a co-oriented hypersurface $V$. Suppose that
$J $  is compatible with  a symplectic form $\omega$ on $W$, and
$(\wt J|_{\RR \times V},\omega|_{\RR \times V})$ satisfy  the
symmetry condition from Section
\ref{sec:cyl-structures}. Then for every  $E>0$, the space
$\overline \Mc_{g,\mu}(W^{[0,\infty]},J^{[0,\infty]}) \cap \{E(F) \leq E\}$
is compact.
\end{theorem}

  First note that   in all three cases it is enough to prove
    the sequential compactness because  the  corresponding moduli
spaces $\overline\Mc_{g,\mu}(V)$
    are metric spaces.
   Next, it is enough to
consider sequences of curves of height $k=1$. Indeed,  we can handle
each level separately.
{\sl  Moreover, we can assume   all these curves to be smooth, ie, having no
  double points  $D$,  because the double points can be treated as
extra marked points.}
  Finally,
the energy bound and the Morse--Bott condition guarantee that there
are only finitely many
possibilities for the asymptotics at the punctures,
Hence it is enough to prove the following three theorems.
\begin{theorem} \label{thm:compact-cyl2}
Let $$\bF_n=\left(F_n=(a_n,f_n);S_n,j_n,M_n,
\uZ_n\cup \oZ_n\right) ,\;$$
$$ \oZ_n=
\left\{\left(\oz_1\right)_n, \dots, \left(\oz_{p^+}\right)_n\right\}\;\;,
  \uZ_n=
\left\{\left(\uz_1\right)_n, \dots, \left(\uz_{p^-}\right)_n\right\},$$
be a sequence of  smooth holomorphic   curves
in $(W=\RR \times V,J)$ of the
same signature $(g,\mu, p^-,p^+)$ and which are asymptotic
at the  corresponding punctures to orbits from the same component of
the space of periodic orbits $\Pc$.
Then there exists a subsequence
that converges to a stable holomorphic building $F$ of height $k$.
\end{theorem}
\begin{theorem} \label{thm:compact-ends2}
Let $$\bF_n=\left(F_n;S_n,j_n,M_n,
\uZ_n\cup \oZ_n\right) ,\;$$
$$ \oZ_n=
\left\{\left(\oz_1\right)_n, \dots, \left(\oz_{p^+}\right)_n\right\}\;\;,
  \uZ_n=
\left\{\left(\uz_1\right)_n, \dots, \left(\uz_{p^-}\right)_n\right\},$$
be a sequence of  smooth holomorphic   curves
in $(W,J)$  of the
same signature $(g,\mu, p^-,p^+)$ and which are asymptotic
at the  corresponding punctures to orbits from  the same component of
the space of periodic orbits $\Pc$.
Then there exists a subsequence
that converges to a stable holomorphic building $\bF$ of height $k_-|1|k_+$.
\end{theorem}
\begin{theorem} \label{thm:compact-split2}
Let $$\bF_n=\left(F_n;S_n,j_n,M_n,
\uZ_n\cup \oZ_n\right) ,\;$$
$$ \oZ_n=
\left\{\left(\oz_1\right)_n, \dots, \left(\oz_{p^+}\right)_n\right\}\;\;,
  \uZ_n=
\left\{\left(\uz_1\right)_n, \dots, \left(\uz_{p^-}\right)_n\right\},$$
be a sequence of  smooth holomorphic   curves
in  manifolds $(W^n,J^n)$  converging  to a split manifold $(\wt
W,\wt J)$. Suppose all the curves have the
same signature $(g,\mu, p^-,p^+)$  and are asymptotic
at the  corresponding punctures to orbits from the same components of
the space of periodic orbits $\Pc$.
Then there exists a subsequence
that converges to a stable holomorphic building $\bF$ of height
$\mvee\limits_1^{k_0}$.
\end{theorem}
\subsection{Proof of Theorems \ref{thm:compact-cyl2},
\ref{thm:compact-ends2} and \ref{thm:compact-split2}}
The  proofs of all three theorems are very similar, though differ in
some details.
  In each of the four  steps of the
proof we first discuss in detail the cylindrical case, and then indicate
  which changes, if any, are necessary for the two other  theorems.
\subsubsection{Step 1: Gradient bounds}
\noindent  First, we observe that there is a bound $N=N(E)$,
depending only on the energy
    of the curve,
on the number of marked points which should be added to stabilize
the underlying surfaces $\bS^F_n$. Hence, by adding $N$ extra marked
points we can assume that
all the surfaces $\bS_n=\bS^F_n=(S_n,j_n,M_n,\uZ_n\cup\oZ_n)$ are stable.
Using Theorem \ref{thm:DM}  we may assume, by passing to a
subsequence, that the sequence
of Riemann surfaces
  $\bS_n$   converges to a decorated Riemann nodal surface
   $$\bS=\{S,j,M,D,\oZ\cup  \uZ,r\}.$$
    Next, we are going to add more marked points to obtain gradient bounds
on the resulting  punctured surfaces. First, we do it in the
cylindrical case. We continue to write $\rho (x)=\inj (x)$ for the
injectivity radius.
\begin{lemma}\label{lm:bubble-analysis}
  There exists an integer $K=K(E)$ which depends only on the energy bound
  $E$ such that,  by adding to   each marked point set $M_n$  a disjoint set
   $$Y_n=\{y^{(1)}_n,u^{(1)}_n,\dots, y^{(K)}_n,u^{(K)}_n\}\subset
\dot{S}_n=S_n\setminus (M_n\cup\uZ_n\cup\oZ_n)$$
  of cardinality $2K$,
  we can arrange   a uniform gradient bound
  \begin{equation}\label{grad-bound}   
  \|\nabla F_n(x)\|\leq \frac{C}{\rho (x)},\;x\in \dot{S}_n\setminus Y_n
  \end{equation}
  where the gradients are computed with respect to the cylindrical
metric on $\R\times V$
  associated with a fixed Riemannian
  metric on $V$, and the hyperbolic metric on $\dot{S}_n\setminus
Y_n$, and where $\rho (x)$
  is the injectivity radius of this hyperbolic metric at the  point
$x\in \dot{S}_n\setminus Y_n$.
  \end{lemma}
\begin{proof}
Suppose  we are given a sequence of points  $x^{(1)}_n\in \dot{S}_n$
which   satisfies
the property
\begin{equation*}
\lim\limits_{n\to\infty}\rho (x^{(1)}_n)\|\nabla F_{n}(x^{(1)}_n) \|
\rightarrow \infty.
\end{equation*}
By translating the maps $F_n=(a_n,f_n)$ along the $\R$--factor of
$\bR\times V$ we can arrange that
$a_n  (x^{(1)}_n)=0$ for all $n$. There exist (injective) holomorphic
charts $\psi_n\co D\to \cD_n\subset
\dot{S}_n$
with $\psi_n(0)=x^{(1)}_n$
and with $$C_1\rho \bigl(x^{(1)}_n\bigr)\leq\|\nabla\psi_n\|\leq
C_2\rho  \bigl(x^{(1)}_n\bigr)$$ for two positive constants $C_1,C_2$.
This can easily been seen by taking fundamental domains in the
hyperbolic upper half-plane
  uniformizing components of the thin part of $\dot S_n$.
Then  we have $\|\nabla(F_n\circ\psi_n)(0)\|\to\infty$ as
$n\to\infty$ and hence, using  Lemma \ref{lm:bubble}
    we conclude
  that   there exist sequences $y^{(1)}_n \rightarrow 0$ and $c_n, R_n
\rightarrow \infty$ as
$n \rightarrow \infty$ such that the rescaled maps
\begin{equation*}
{\wt {F}_n} \co  D_{R_n} \rightarrow (W, J): z
\mapsto F_{n}\circ\tilde{\psi}_n(z),
  \end{equation*}
where
  $$
  \tilde{\psi}_n(z):=\psi_n \biggl (y^{(1)}_n + c_n^{-1} z\biggr),
  $$
converge to a holomorphic map $\wt F_\infty$ satisfying $E(\wt
F_\infty) \leq C$ and $E_\omega(\wt F_\infty) > \hbar$.
Here   $D_{R}$ denotes the disc $\{|z|<R\}\subset \C$.
Moreover, this map is either a holomorphic sphere or
  a holomorphic plane $\bC$ asymptotic as $\abs{z}\to \infty$ to a
closed $\Reeb$--orbit. Let us
  choose a sequence $u_n^{(1)}=y^{(1)}_n+c_n^{-1}\in D$ and set $\bar
y^{(1)}_n=\psi_n(y^{(1)}_n)$,
  $\bar u^{(1)}_n=\psi_n(u^{(1)}_n)$.
  Then  $\bar y^{(1)}_n$ and $\bar u^{(1)}_n$ are distinct points in
  $\psi_n(D_{R_n})\subset\cD_n$,
and $\dist_n(\bar y^{(1)}_n,\bar
u^{(1)}_n)\mathop{\to}\limits_{n\to\infty} 0$, where $\dist_n$ is
the distance function on $\dot{S}_n$ defined by the hyperbolic
metric $h^{j_n,M_n\cup\uZ_n\cup\oZ_n}$. Thus according to
Proposition \ref{prop:2points} (a subsequence of) the sequence of
marked Riemann surfaces $\bS^{(1)}_n=(S_n,j_n,M_n\cup\{\bar
y^{(1)}_n,\bar u^{(1)}_n\},\uZ_n\cup\oZ_n)$ converges to a nodal
decorated Riemann surface $\bS^{(2)}$ obtained
from $\bS^{(1)}=\bS$ by adding one or two spherical components.
Figure \ref{Newfig12}
illustrates a possible scenario, while all possible cases are
illustrated by Figure \ref{Newfig1}.
  Note that by construction
  exactly one  of these components contain  the  marked points $\bar
y^{(1)}$ and $\bar u^{(1)}$,
   which correspond  to
  the sequences $\bar y^{(1)}_n$. This bubble serves
   as the domain   of the map $\wt F_\infty$, and thus have the
   $\omega$--energy concentration
for large $n$  exceeding
$\hbar$.
\begin{figure}[ht!]\small
\centering
\psfrag{x}{$x_i$}
  \psfraga <-1pt, -1pt> {y}{$\bar{y}^{(1)}$}
\psfrag{u}{$\bar{u}^{(1)}$}
\includegraphics[width=0.63\textwidth]{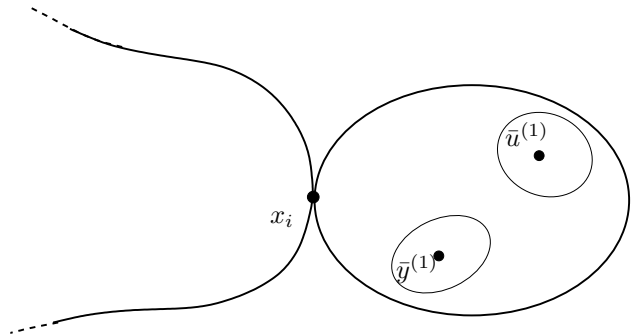}
\caption{The surface $\bS^{(2)}$ a
spherical bubble containing  the points $\bar y^{(1)}$ and $\bar
u^{(1)}$}\label{Newfig12}
\end{figure}
Set now $M_n^{(1)}=M_n\cup\{y^{(1)},u^{(1)}\}$ and repeat the above
analysis for the sequence
of holomorphic curves
$$\bF_n^{(2)}=\left(F_n=(a_n,f_n);S_n,j_n,M^{(1)},\uZ_n\cup\oZ_n\right)\,.$$
If the inequality  \eqref{grad-bound}  does not hold for the
hyperbolic metric on
$$\dot{S}^{(2)}=\dot{S}_n\setminus\{y_n^{(1)},
u_n^{(1)}\}$$ then  we repeat the above bubbling-off analysis, thus
constructing:
\begin{itemize}
\item  a sequence
of points  $x^{(2)}_n\in \dot{S}^{(2)}_n$ having
the property
\begin{equation*}
\lim\limits_{n\to\infty}\rho  \bigl(x^{(2)}_n\bigr)\| \nabla
F_{n}(x^{(2)}_n) \| \rightarrow \infty.
\end{equation*}
\item holomorphic charts
$\psi^{(2)}_n\co D\to \cD_n^{(2)}\subset
\dot{S}_n$ satisfying  $\psi^{(2)}_n(0)=x^{(2)}_n$
and $$C^{(2)}_1\rho
\bigl(x^{(2)}_n\bigr)\leq\|\nabla\psi^{(2)}_n\|\leq C^{(2)}_2\rho
\bigl(x^{(2)}_n\bigr),$$
  for two positive constants $C^{(2)}_1,C^{(2)}_2$.
  \item sequences $y^{(2)}_n \rightarrow 0$ and $c^{(2)}_n, R^{(2)}_n
\rightarrow \infty$ as
$n \rightarrow \infty$ such that the  rescaled maps
\begin{equation*}
\begin{split}
  \wt F^{(2)}_n  \co  D_{R^{(2)}_n} &\rightarrow (W, J),\cr
    z &\mapsto F^{(2)}_{n}\circ\wt\psi^{(2)}_n(z),\cr
  \end{split}
  \end{equation*}
  where $$\wt\psi^{(2)}_n(z)=\psi^{(2)}_n\biggl(y^{(2)}_n +
\frac{z}{c^{(2)}_n}\biggr),$$
converge to a holomorphic map $\wt F^{(2)}_\infty$ satisfying $$E(\wt
F^{(2)}_\infty) \leq C\;\;
\hbox{ and}\;\;
  E_\omega(\wt F^{(2)}_\infty) > \hbar.$$
  \end{itemize}
Notice that there exist sequences   $K^{(1)}_n,
K^{(2)}_n\mathop{\to}\limits_{n\to\infty}\infty$
  satisfying {
  $$K^{(1)}_n<R_n \quad \text{and}\quad K^{(2)}_n<R^{(2)}_n,$$  such
that  the discs
  $\tilde{\psi}_n(D_{K^{(1)}_n})$ and  $\tilde{\psi}^{(2)}_n(D_{K^{(2)}_n})$
  do not intersect.}  Indeed, for any fixed $K>0$ the gradients
  $\|\nabla F_n\|$ are uniformly bounded on $\psi_n(D_{K})$ and
  go to $\infty$ as $n\to\infty$ on $\psi^{(2)}_n(D_{K})$.
  Set $$u^{(2)}_n=y_n^{(2)}+\frac1{c_n^{(2)}},\;\;
  \bar y^{(2)}_n=\psi^{(2)}_n(y^{(2)}_n),\;\;
  \bar u^{(2)}_n=\psi^{(2)}_n(u^{(2)}_n).$$
  Then the sequence
$$\bF_n^{(3)}=\left(F_n=(a_n,f_n);S_n,j_n,M^{(2)}=M_n^{(1)}\cup\{y^{(2)}_n,u^{(2)}_n\},\uZ_n\cup\oZ_n\right)\,.$$
  has in the limit an extra bubble $\wt F^{(2)}_\infty$   disjoint
from $\wt F^{(1)}_\infty=\wt F_\infty$ whose energy satisfies
$E_\omega(\wt F^{(2)}_\infty)>\hbar$.
Hence, the uniform bound on the $\omega$--energy guarantees
  that,  after adding finitely many pairs of marked points
$$y^{(1)}_n,u^{(1)}_n,\dots, y^{(K)}_n,u^{(K)}_n, $$
the inequality \eqref{grad-bound} holds true  for the hyperbolic metric
  $h^{j_n,M^{(K)}_n\cup\uZ_n\cup\oZ_n}$,  where
$M_n^{(K)}=M_n\cup\{y^{(1)}_n,u^{(1)}_n,\dots,
y^{(K)}_n,u_n^{(K)}\}$. The proof of Lemma 
is complete.
\end{proof}

\noindent{\bf Case of manifolds with cylindrical ends and the
splitting case}\nl  Lemma \ref{lm:bubble-analysis} has obvious
analogues in these two cases. In the case of mappings to a
manifold $W$ with cylindrical ends  the gradients should be
computed with respect to a fixed metric on $W$ which is
cylindrical at the ends. In the splitting case the gradients are
computed with respect to a sequence of metrics on the target
manifold  which arise in the process of splitting. These
metrics have  longer and longer
cylindrical inserts.  The proof works  without any serious
changes, except that one needs to analyze separately the case in which
the sequence  $F_n(x_n)$
   stays in a compact subset
of $W$ and the case in which  (a subsequence of) it is contained in
the cylindrical part.
In the second case the proof is identical while in
the first case    it does not make sense, and there is no need  to
shift the map in order to fix its
  $\R$--component.
  \subsubsection{Step 2: Convergence of Riemann surfaces and
convergence away from nodes}
This step is common in  all three theorems.

\noindent We will assume from now on that we added enough extra
marked points  to stabilize
the underlying surfaces $\bS^F_n$ and to ensure the gradient bounds
\eqref{grad-bound}.  Then using Theorem \ref{thm:DM}  we may assume,
by passing to a subsequence, that the sequence
of Riemann surfaces
  $\bS_{F_n}$   converges to a decorated Riemann nodal surface
   $$\bS=\{S,j,M,D,\oZ\cup  \uZ,r\}.$$
 From section \ref{sec:DM-topology} we  recall that in our situation,
where $D_n=\emptyset$ and where $M_n\cup Z_n$  and $M\cup Z$  in
abuse of notation are again denoted by $M_n$ and $M$, this means the
following. There exists a sequence of diffeomorphisms
$$\varphi_n\co S^{D, r}\to S_n\:\: \text{with $\varphi_n (M)=M_n$}$$
having the following properties.
\begin{itemize}
\item[$\bullet$] There exist disjoint closed geodesics $\Gamma^n_i$
for $i=1,\ldots, k$ on $S_n\setminus M_n$ with respect to the
hyperbolic metrics $h^{j_n, M_n}$, for all $n\geq 1$, such that
$\Gamma_i:=\varphi^{-1}_n (\Gamma_i^n)$ are special circles on $S^{D,
r}\setminus M$.
\item[$\bullet$]
$\varphi_n^* j_n\to j $ in $C^{\infty}_{\text{loc}}\bigl(S^{D, r}\setminus
\bigcup\limits_1^k\Gamma_i\bigr).$
\item[$\bullet$]  Given a   component $C$ of
   $\Tn_\e(\bS)\subset \dot S^{D,r}$ which contains a special circle
   $\Gamma_i$ and given a point  $c_i\in\Gamma_i$,
   we consider for every $n\geq 1$  the geodesic arc $\d^n_i$ for the
induced  metric
$h_n=\varphi^*_nh^{j_n, M_n}$ which  intersects $\Gamma_i$
orthogonally
  at the point $c_i$,
and whose  ends  are contained in the $\varepsilon$--thick part of
  the metric $h_n$. Then $(C\cap\d^n_i)\;$
  converges as $n\to\infty$ in  $C^0$  to a
continuous geodesic for the metric  $h^{\bS}$  which passes
through the point $c_i$.
    \end{itemize}
In addition,  according to Lemma \ref{lm:bubble-analysis} we may assume that
  \begin{equation*}\label{eq:grad-bound}
  \|\nabla (F_n\circ\varphi_n)(x)\|\leq \frac{C}{\rho (x)},\;x\in
{S}\setminus\bigcup\Gamma_i\,.
  \end{equation*}
These gradient bounds allow
to apply locally the Gromov--Schwarz Lemma \ref{lm:Gromov-Schwarz}
  to conclude uniform bounds for {\it all} derivatives of $F_n\circ
\varphi_n$ on the $\e$--thick part
  on $\Tk_\e(\bS)$ for every  $\e>0$, and therefore
  Ascoli--Arzela's  theorem  allows us
   to extract a subsequence converging in $C^{\infty}_{\text{loc}}$ on
   $S\setminus \bigcup\Gamma_i=\bigcup\limits_\e\Tk_\e(\bS).$\footnote
   {Let us recall that the convergence in cylindrical manifolds is
defined up to translation
   along the $\R$--factor. In particular, when the surface
    $S\setminus \bigcup\Gamma_i$ is disconnected one may need to shift
   the maps of the sequence restricted to different components by
different constants.}
\subsubsection{Step 3: Convergence in the thin part}
{\bf Cylindrical case}\qua
Let us denote by $C_1,\dots, C_N$ the connected components of
$S\setminus \bigcup\Gamma_i$.
We already may assume that the holomorphic maps $F_n\circ\varphi_n$
converge on each component $C_i$ for
$i=1,\dots, N$.
Our next goal is  to understand the asymptotic behavior of the  limit
map $F=(a, f)$ on
the component $C_i$ near a node. First, if $F$ is bounded near the node, then,
by the removable singularity theorem, Lemma  \ref{punctfin}, the map
$F=(a, f)$ extends continuously on
$C_i$ across the node. On the other hand, if $F$ is unbounded near
the node, the behavior
of $F$ is described by Proposition  \ref{lm:asymptot}.
Namely, there exists a closed R--orbit $\gamma\in\Pc$ such that the map $F$
is asymptotic to $\gamma$ near the node either at the positive or  at
the negative
end. Moreover, the map $f$
extends continuously to the circle at infinity which compactifies
the puncture.
\medskip

\noindent{\bf Behavior near a node adjacent to two components}\qua
Given a node of $S$ adjacent to the two  components $C_i$ and
$C_j$, the asymptotic behavior of $F$ on the two components might
be different at first sight. For example, $F$ could be asymptotic to different
$\Reeb$--orbits, or $F$ could be asymptotic to an $\Reeb$--orbit on $C_i$ and
could converge to a point  on $C_j$, or  it could converge  to
different points.   Even if $F$ is asymptotic to
the same orbit on $C_i$ and $C_j$ we still have to worry about  a
loss of $\omega$--energy (which may happen only in the non-contact
case), and a possible shift in the asymptotic
parameterizations of the orbits.
To each node  adjacent to two components, we can associate
two asymptotic limits $\gamma^+$ and $\gamma^-$,   one for each
component of $S\setminus\bigcup\Gamma_i$ adjacent to the node. Each
$\gamma^{\pm}$ is either
a point or a periodic orbit from $\Pc$.
\begin{figure}[ht!]\small
\centering
\psfrag{st}{$S^2$}
\includegraphics[width=0.7\textwidth]{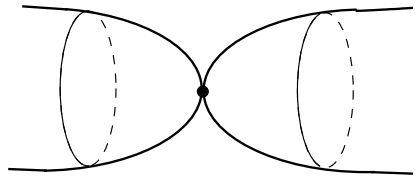}
\caption{Node adjacent to two components}\label{fig13}
\end{figure}
  The node in question appeared as a result of the degeneration of
a  component of the $\e$--thin part of $S_n$. In other words,
there exists a component $T^\e_n$  of the $\e$--thin part
  of the hyperbolic  metric $ h^n =\varphi_n^*h^{j_n,M_n}$ on
$S=S^{D,r}$,  with  conformal
  parametrization
$$
  g^\e_n\co 
A^{\e}_n=[-N^\e_n,N^\e_n]\times S^1\to (T^\e_n,j_n),
$$
such that in the $C^{\infty} (S^1)$--sense,
  \begin{equation}\label{eq:thin-limit}
  \lim\limits_{\e\to 0}
\lim\limits_{n\to\infty}f_n\circ\varphi_n \circ  g^\e_n|_{(\pm
N^\e_n)\times S^1}=\gamma^\pm\, ,
\end{equation}
where $\gamma^{\pm}$ are either two periodic orbits of $\Reeb$ or two
points of $V$.
Moreover, for any constant $K$ we have
\begin{equation}\label{eq:minus-constant}
\lim\limits_{n\to\infty}f_n\circ\varphi_n \circ  g^\e_n|_{(\pm
N^\e_n\mp K)\times S^1}=\gamma^\pm\,.
\end{equation}
\begin{remark}\label{rm:orientation}
{\rm
It is possible that both orbits $\g^+$ and $\g^-$ may appear at the
same end of the cylindrical manifold.
In this  case one of the orbits must have an  opposite orientation.
For the following discussion the orientation
of $\g^\pm$ will be irrelevant, and thus it will not  be specified it
in our notation.}
\end{remark}
Note that the parameterizations $g^{\e}_n$ can be chosen so that they
satisfy the gradient
bounds
$$\|\nabla g^{\e}_n(x)\|\leq C\rho  \bigl(g^{\e}_n(x)\bigr),$$
where the gradients are computed with respect
to the flat metric  in the source and the hyperbolic metric in the
target, while the
  injectivity radius is computed with respect to the hyperbolic metric.
   Together with the estimate \eqref{grad-bound} this implies a  uniform
   (ie, independent of $n$ and $\e$) gradient bound
   \begin{equation}\label{eq:grad-bound2}
   \sup\limits_{x\in A_n^{\e}}\|\nabla\left(F_n\circ\varphi_n \circ
g^\e_n(x)\right)\|
   \leq C\,.
   \end{equation}
Given  a sequence $\e_k\to 0$  let us choose a subsequence
$\e_{k_n}\to 0$ such that
$$\lim\limits_{n\to\infty}F_{k_n}\circ\varphi_{k_n}\circ
g^{\e_{k_n}}_{k_n}|_{\pm N^{\e_{k_n}}_{k_n}\times S^1}
=\gamma^\pm\, ,$$
and introduce  the abbreviated notation $\wh N_{n}=N^{\e_{k_n}}_{k_n}$,
$\wh g_n=\varphi_{k_n}\circ g^{\e_{k_n}}_{k_n}$,
$\wh f_n=f_{k_n}\circ \wh g_{k_n}$ and      $\wh F_n=F_{k_n}\circ \wh
g_{k_n}$, so that   we get
$$\lim\limits_{n\to\infty}\wh f_n(\pm\wh N_n\times S^1)=\g^\pm .$$
For large $n$,
  the loops $\wh f_n|_{(\pm \wh N_n)\times S^1}$ are sufficiently
$C^{\infty}(S^1)$-=close to $\g^\pm$ and hence the cylinder $\wh
f_n|_{[-  \wh N_n,\wh N_n]\times S^1}$
defines a   homotopically unique map $\Phi\co S^1\times[0,1]\to V$
satisfying $\Phi|_{S^1\times 0}=\gamma^-$ and
$\Phi|_{S^1\times 1}=\gamma^+$.
    We can assume that the homotopy class  of $\Phi$ is independent of $n$.
In  the notation of Section \ref{sec:dynamics} we
distinguish  the following two cases,
\begin{description}
\item{\bf C1}\qua $\Delta S_\omega(\g^+,\g^-;\Phi)=0$
\item{\bf C2}\qua $\Delta S_\omega(\g^+,\g^-;\Phi)>0$.
\end{description}
\medskip

\noindent {\bf Case C1}\qua We shall show in this case that   $\gamma^+$ and
$\gamma^-$  are geometrically the
same, and that also  their parameterizations are the same.
Moreover,  we shall show that in this case the limit map $f$
extends continuously to the circle $\Gamma$ which is associated to this node.
Assume first that one of the asymptotic limits, say $\gamma^-$, is an
$\Reeb$--orbit.
Let us  show that $\gamma^+$ is  also   an $\Reeb$--orbit and
  $\gamma^-(t)=\gamma^+(t)$ for all $t$. Indeed, by assumption,
   $E_{\omega}(\wh F_n\vert_{[-\wh N_n, \wh N_n]})\to 0$ and $E(\wh
F_n)\leq E_0$. Thus we may apply Proposition
    \ref{prop:twist}
  and find for every $\sigma>0$ a constant $c>0$ so that $\wh f_n(s,
t)\in B_{\sigma}(\wh f_n (0, t))$ for all
  $(s, t)\in [-\wh N_n+c, \wh N_n-c]\times S^1$ and $n$ large enough.
  Since    $\wh f_n(\pm \wh N_n, t)\mathop{\to}\limits_{n\to\infty}
\gamma^{\pm}(t)$
   we conclude  that
  $\gamma^+(t)=\gamma^-(t)$ for all $t$. This also proves that  the
limit map $f$ continuously
extends to the circle $\Gamma$ associated to this node.
Similarly, this time using Lemma \ref{lm:almost-sphere},
  one shows that if $\gamma^-$ is a point, then also $\gamma^+$ is a
point and $\gamma^-=\gamma^+=p$.
   Moreover, for large $n$ the image $\wh F_n([-\wh N_n, \wh
N_n]\times S^1)$ is contained in an
    arbitrary small neighborhood of a point $p^{\ast}=(a, p)\in \bR\times V$.
\medskip

\noindent{\bf Case C2}\qua  Assume now $\Delta S_\omega(\g^+,\g^-;\Phi)
=\delta>0$. We point out
that the uniform energy bound implies a uniform upper bound for the
periods of $\gamma^{\pm}$.
  This implies, as we shall first show, an a priori lower bound for $\delta$.
   Namely, we have
   \begin{lemma}\label{lm:lower-bound}
   There exists a ``quantum constant" $\hbar=\hbar (E)>0$ such that
\begin{equation}
\label{eq:hbar}
\Delta S_\omega(\g^+,\g^-;\Phi) >\hbar\, .\footnote{Let us recall
that we allow the orbits $\g^\pm$ be oriented by the vector
field $-\Reeb$ (see Remark \ref{rm:orientation}). However, we always
assume, that their orientation is
  chosen in such a way
that $\Delta S_\omega(\g^+,\g^-;\Phi)>0$.}
\end{equation}
\end{lemma}
\begin{proof}
To prove the claim we argue indirectly and assume that
we have sequences of orbits (or points) $\gamma^{\pm}_k$ with
  bounded periods and holomorphic cylinders $\wh F_{n, k}\co [-\wh
N_{n,k}, \wh N_{n,k}]\times S^1\to \bR \times V$
  satisfying
\begin{equation*}
\begin{gathered}
\lim\limits_{n\to\infty}\wh N_{n,k}=\infty\\
\lim\limits_{n\to\infty}\wh F_{n,k}(\pm\wh N_{n,k}\times S^1)=
\g^{\pm}_k\:\: \text{in $C^{\infty}(S^1)$}\\
\Delta S_\omega(\g_k^+,\g_k^-;\Phi_k)=\delta_k>0\\
\end{gathered}
\end{equation*}
and
\begin{equation*}
\delta_k\to 0\quad \text{as $k\to \infty$.}
\end{equation*}
Here $\Phi_k$ is the relative homotopy class of $\wh f_{n, k}$ which
is  well defined for large $n$.
  Now arguing as in C1 we can choose a diagonal subsequence
   $$\wh F_{n,k_{n}}\co [-\wh N_{n,k_{n}}, \wh N_{n, k_{n}}]\times S^1\to
\bR\times V$$
   which converges uniformly to a trivial cylinder over a periodic
orbit $\gamma$,
    or to a point $\gamma \in \bR\times V$, in the following sense.
There exists a constant $c$ such that
$$
\sup_{[-\wh N_{n,k_{n}}+c, \wh N_{n, k_{n}}-c]\times S^1}d\left(\wh
f_{n,k_{n}}(s, t), \gamma (t)\right)
\xrightarrow[\tau \to \infty] {}  0\,$$ if $\g$ is a periodic orbit, and
$$\sup_{[-\wh N_{n,k_{n}}+c, \wh N_{n, k_{n}}-c]\times S^1}d\left(\wh
F_{n,k_{n}}(s, t), (0,\gamma)\right)
\xrightarrow[\tau \to \infty] {}  0\,$$ if $\gamma (t)=\gamma$  is a point.
In the latter case we assume the $\R$--component
of the maps  $\wh F_{n,k_{n}}$  to be fixed by the condition $\wh
F_{n,k_{n}}(0,0)=0$.
In the situation where $\gamma$ is a periodic orbit we recall
  that in view of compactness of $V$ and
   the Morse--Bott condition,   the periods of $\gamma^+_{k_n}$ and
$\gamma^-_{k_n}$ are equal to that of $\gamma$ for large $n$.
Hence $\delta_{k_n}=\Delta S_{\omega}(\gamma^+_{k_n}, \gamma^-_{k_n};
\Phi_{k_n})=0$
for large $n$, in contradiction to the assumption. If $\gamma$ is a
point we arrive at the same contradiction,
  hence proving the claim \eqref{eq:hbar} above.
\end{proof}
In view of Lemma \ref{lm:bubble} we may assume that the same $\hbar$
also  serves as
a lower bound of the $\omega$--energies of all holomorphic planes and
spheres which appear
as a  result of the bubbling off analysis.
Let us recall that   the inequality \eqref{eq:grad-bound2}  provides
us with the  uniform gradient bounds for the maps
$\wh F_n$.  So, no bubbling off can  occur anymore.
   Analogously to the broken trajectories in Morse theory and Floer
   theory we shall see that the worst which can happen in the limit is
   the splitting of our
long cylinder into a finite sequence of cylinders, which, in the
image under the map meet  at their ends along periodic orbits
$\gamma$ of the vector field $\Reeb$.
\begin{figure}[ht!]\small
\centering
\psfrag{g}{$\gamma$}
\psfrag{gp}{$\gamma^+$}
\includegraphics[width=0.9\textwidth]{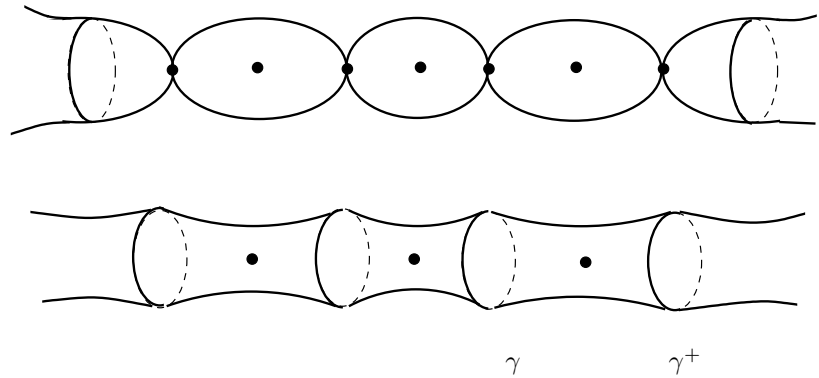}
\caption{Splitting of a long cylinder}\label{fig15}
\end{figure}
To see this we consider as in Case C1
  a sequence $\wh F_n \co [-\wh N_n, \wh N_n]\times S^1\to \bR\times V$
   of cylinders having uniform  gradient bounds and satisfying
   $$E_\omega(\wh F_n)>\hbar,\; \lim\wh f_n|_{-\wh N_n\times S^1}=\gamma^-,\;
   \lim\wh f_n|_{\wh N_n\times S^1}=\gamma^+\,.$$
   Then there exists $K_n<\wh N_n$ such that
$$
\int\nolimits_{[K_n, \wh N_n]}\wh F_n^{\ast}\omega =\hbar.
$$
Then it follows from \eqref{eq:minus-constant}
that $\wh N_n-K_n\to \infty$ as $n\to \infty$.  Choose a point $p_n$
on the circle $K_n\times S^1$ and translate the
$\bR$--coordinate $\wh a_n$ of the map $\wh F_n=(\wh a_n,  \wh f_n)$
in such a way that $\wh a_n(p_n)=0$.
By  means of Ascoli--Arzela theorem (a subsequence of) the sequence $\wh F_n$
converges to a holomorphic cylinder $\wh F_{\infty}\co \bR \times S^1\to
\bR\times V$.
If $\wh F_{\infty}$ is a trivial cylinder over an orbit $\gamma$, or
a constant map to a point, also denoted by $\gamma$, then $\Delta
S_{\omega} (\gamma^+, \gamma; \Phi)=\hbar$,
  where the homotopy class $\Phi$ is determined by the maps $\wh
F_n\vert_{[K_n, \wh N_n]\times S^1}$
   for sufficiently large $n$.  But this contradicts the a priori
estimate \eqref{eq:hbar}.
     Hence the cylinder $\wh F_{\infty}$ is not constant and different
from the trivial cylinder over a periodic orbit,
      and hence has energy $E_{\omega}(\wh F_{\infty})>\hbar$.
      Consequently, by adding the extra marked points $\wh g(p_n)\in
S_n$ we create an
      additional  spherical  component  $C$ of the limit Riemann
surface $\bS$. The  limit map $F_\infty|_C$
       uses up more than $\hbar$ of the  $\omega$--energy. Hence,
iterating this analysis
        of long cylinders  we arrive at the situation where
         the case C2 cannot occur anymore.
\medskip

\noindent{\bf Behavior near a puncture}\qua We should
also analyze the behavior of $F$ on a component of the thin part,
    adjacent to a marked point (puncture).  Fix a positive  $\e<\e_0$.
    The corresponding component $T_n$ of the $\e$--thin part of $S_n$
admits  in this case a holomorphic parametrization
$$g_n\co A=[0,\infty)\times S_1\to(T_n,j_n)\,$$
such that for every $n$ there exists a periodic orbit (or a point)
$\g_n$ satisfying,
in the $C^\infty(S^1)$--sense,
$$
\lim\limits_{s\to\infty}\wh F_n|_{(\pm s)\times S^1}=\gamma_n\, $$
where $\wh F_n=(\wh a_n,\wh f_n)=F_n\circ\varphi_n \circ  g_n$.
Moreover, due to the compactness of $V$ and of the space of periodic
orbits of bounded period we
can assume
that there exists a periodic orbit or a point
$\og=\lim\limits_{n\to\infty} \g_n$.
We can also assume an analogue to  the inequality
\eqref{eq:grad-bound2}, namely
\begin{equation}\label{eq:grad-bound3}
   \sup\limits_{x\in A}\|\nabla \wh F_n(x) \|
   \leq C\,,
   \end{equation}
so that (a subsequence of) the sequence $\wh F_n\co A\to\R\times V$, normalized
by the condition $\wh a_n(0,0)=0$   converges in  the
$C^\infty_\loc$--topology to a holomorphic
map $\wh F\co A\to\R\times V$ asymptotic to a closed orbit, or a point $\ug$.
Choose sequences $\uN_n,\oN_n\mathop{\to}\limits_{n\to\infty}\infty$
and $\uN_n<\oN_n$,
satisfying
$$\lim\limits_{n\to\infty}\wh f_n|_{\uN_n\times S^1}=\ug,\;\;
\lim\limits_{n\to\infty}\wh f_n|_{\oN_n\times S^1}=\og\,.$$
Notice that,  given  two sequences
$\uN_n',\oN_n'\mathop{\to}\limits_{n\to\infty}\infty$ and a constant
$c>0$
satisfying
  $\uN_n\leq \uN_n'\leq \uN_n+c$ and $\oN_n-c\leq\oN_n'\leq\oN_n$, we conclude
\begin{equation}\label{eq:2limits}\lim\limits_{n\to\infty}\wh
f_n|_{\uN_n'\times S^1}=\ug,\;\;
\lim\limits_{n\to\infty}\wh f_n|_{\oN_n'\times S^1}=\og\,.
\end{equation}
For large $n$,
  the cylinder $\wh f_n|_{[  \uN_n,\oN_n]\times S^1}$
defines a   homotopically unique map $\Phi\co S^1\times[0,1]\to V$
satisfying $\Phi|_{S^1\times 0}=\ug$ and
$\Phi|_{S^1\times 1}=\og$. As in the above  analysis  of a node
adjacent to two components
we distinguish  the following two cases:
\begin{description}
\item{\bf C1$'$}\qua $\Delta S_\omega(\ug,\og;\Phi)=0$
\item{\bf C2$'$}\qua $\Delta S_\omega(\ug,\og;\Phi)>0$.
\end{description}
Because these  cases are essentially similar to C1 and C2 we will
restrict ourselves to the situation in which
both $\ug$ and $\og$ are $\Reeb$--orbits.
\medskip

\noindent {\bf Case C1$'$}\qua
Let us show that in this case  $\ug(t)=\og(t),\;t\in S^1$. Indeed,
applying  Proposition
    \ref{prop:twist}
  we find for every $\e>0$ a constant $c>0$ such that
  \begin{equation}\label{eq:twist2}
  \wh f_n(s, t)\in B_{{\e}}(\wh f_n ( N_n, t)),
  \end{equation}
  where $N_n=\frac{\uN+\oN_n}2$, for all
  $(s, t)\in [\uN_n+c, \oN_n-c]\times S^1$ and $n$ large enough.
Clearly, (see \eqref{eq:2limits} above)
   $\lim\limits_{n\to\infty}
  \wh f_n(\uN_n+c,t)=\lim\limits_{n\to\infty}
  \wh f_n(\uN_n,t)=\ug(t)$ and $\lim\limits_{n\to\infty}
  \wh f_n(\oN_n-c,t)=\lim\limits_{n\to\infty}
  \wh f_n(\oN_n,t)=\og(t)$, $t\in S^1$. On the other hand,
\eqref{eq:twist2} implies that
$\lim\limits_{n\to\infty}
  \wh f_n(\uN_n+c,t)=\lim\limits_{n\to\infty}
  \wh f_n(\oN_n-c,t).$
  Hence, $\ug(t)=\og(t),\;t\in S^1$. In fact, a similar argument
implies a stronger statement. Namely,
  the sequence $\wh F_n\co A\to\R\times V$ converges to the holomorphic map
$\wh F$ {\it uniformly}, ie, for any $\e>0$ there exists $N,K>0$ such that
$\wh f_n(s,t)\in B_\e\left(\wh f(s,t)\right)$ for all $t\in S^1$,
$n\geq K$ and $s\geq N$.
\medskip

\noindent {\bf Case C2$'$}\qua The analysis of this case is identical to
the case C2.
\medskip

\noindent{\bf Manifolds with cylindrical ends and splitting}\qua
The proof is essentially the same as for  cylindrical manifolds. The
only  difference  arises when,
   all  the images $\wh F_n([-\wh N_n,\wh N_n]\times S^1)$ for large
$n$  intersect the non-cylindrical part
of $W$. In the analysis of Cases C1 and C1$'$, when we need to show
that if $\g^-$is  a point then $\g^+$ is the same point,
we may  use  the Monotonicity Lemma \ref{lm:monotonicity} instead of
Lemma \ref{lm:almost-sphere}. Similarly,
the Monotonicity Lemma implies Lemma \ref{lm:lower-bound} for the
case when either $\g^-$ or $\g^+$ is a point in the
non-cylindrical part. The rest of analysis is the same with the only
difference,
that the notion of convergence in the non-cylindrical
part does not involve any  freedom of translation.
\subsubsection{Step 4: Level structure}
{\bf Cylindrical case}\qua Let us introduce an ordering in the set of
components of $S\setminus\bigcup\Gamma_i$.
Given  two components $C_i$ and $C_j$,  we choose  two points $x_i \in C_i$ and
$x_j \in C_j$, and define  $C_i \le C_j$ if
  $$\limsup\nolimits_{n \rightarrow \infty}[ a_n\circ \varphi_n (x_i)
- a_n\circ \varphi_n (x_j)] < \infty.$$
If $C_i \le C_j$ and $C_j \le C_i$, then we  write  $C_i \sim C_j$.
Clearly, this ordering is independent of the choice of the points $x_i$ and
$x_j$.
Now, we can label the components $C_i$ with their level number as follows.
The set of  components minimal with respect to  the above ordering
will constitute  level 1. Then,
after removing these components, the set of minimal components will be of
level 2, etc.
Clearly, this labelling is constant across nodes that are mapped at finite
distance. However, it may happen that the level number jumps by an integer
$N > 1$ across a node in the limit.
   In that case, we have to insert $N-1$
additional components between these two components, each of them a vertical
cylinder over the  orbit corresponding to the above node. We can do
this by adding additional points so that the images under the maps
have the appropriate behavior.
Finally, we remove the marked points that we added in all the
previous steps of the proof.
If level $i$ becomes unstable because of this, we remove it and decrease
by $1$ the labelling of higher levels.
Hence, we obtain a level structure that satisfies all the conditions
for a stable  holomorphic building of height $k$ in the cylindrical
almost complex manifold $(\R\times V,J)$.
\medskip

\noindent{\bf Manifolds with cylindrical ends}\qua
  In every component $C_i$ we  pick a point $x_i\in C_i$.
We assign to a component $C_i$ the level number $0$  if $F_n\circ
\varphi_n(x_i)$ is contained in a compact part of $W$ for all $n$.
This property is independent of the choice of the point $x_i$.
  For any other component $C_i$ the sequence $F_n\circ \varphi_n
(x_i)$ is contained in one of the end components of $W$  for $n$
sufficiently large, and,  in particular, $F_n\circ \varphi_n (x_i)$
can be written as $(a_n\circ \varphi_n(x_i),f_n\circ \varphi_n
(x_i))$.
  Let us  introduce an ordering on the set of components $C_i$ associated to an
  end component $E$.
We write  $C_i \le C_j$ if
$\limsup_{n \rightarrow \infty} [a_n\circ \varphi_n (x_i) - a_n\circ
\varphi_n (x_j)] < \infty$.
If $C_i \le C_j$ and $C_j \le C_i$, then we  write  $C_i \sim C_j$.
Clearly, this ordering is independent of the choice of the points $x_i$ and
$x_j$.
Now, we can label the components $C_i$ with their level number as
follows. If $E$ is a positive end then
the set of minimal components for the above ordering will be of level 1. Then,
after removing these components, the set of minimal components will be of
level 2, etc. If $E$ is a negative end then the set of maximal
components for the above ordering will be of level -1.
Then,
after removing these components, the set of maximal components will be of
level -2, etc.
Clearly, this labelling is constant across nodes that are mapped at finite
distance. However, it may happen that the level number jumps by an integer
$N > 1$ across a node at infinity. In that case, we have to insert $N-1$
additional components between these two components, each of them a vertical
cylinder over the  orbit corresponding to the above node.
Finally, remove the marked points that we added in all the previous
steps of the proof.
If level $i$ becomes unstable because of this, we remove it and for a
positive (resp. negative) end
decrease (resp. increase) by $1$ the labelling of higher (resp. lower) levels.
Hence, we obtain a level structure that satisfies all the necessary conditions
for a stable  holomorphic building of height $k_-|1|k_+$ in the
almost complex   manifold
$W$ with cylindrical ends.
\medskip

\noindent{\bf Splitting}\qua
Again choosing in each component $C_i$ a point $x_i\in C_i$ we assign
to a component $C_i$ the main layer if
  there exists $N_i>0$ such that for  sufficiently large  $n$ we have
  $$F_n\circ \varphi_n (x_i)\in\circW\cup\left([-n,n]\setminus
(-N_i,N_i)\right)\times V\subset W^n\,.$$
  For all the other components $F_n\circ \varphi_n (x_i)$ is contained
inside the cylindrical part $[-n,n]\times V\subset
  W^n$for large $n$, and hence one can define the  partial order, and
after  that the labelling of the levels exactly as it was done  in the
  cylindrical case above. The proofs of the Theorems
\ref{thm:compact-cyl2}--\ref{thm:compact-split2} are complete.
\noindent \qed
  \section{Other compactness results}\label{sec:other}
  \subsection{Cylindrical structure over a non-compact manifold}
\label{sec:noncompact}
  It is clear that the compactness theorem for holomorphic curves
  in a cylindrical
  almost complex manifold $(\R\times V, J)$ fails if $V$ is
non-compact and if one does
  not impose some extra conditions on the behavior of $J$. In this
section we discuss
one type of conditions under which  the compactness theorem  can
still  be proven.
A cooriented hypersurface  $\Sigma\subset V$ is called {\it
pseudo-convex} if  $\wt\Sigma=\R\times \Sigma$ is
a (non-strictly) pseudo-convex hypersurface in the almost complex
manifold $\R\times V$.
\begin{theorem}\label{thm:compact-convex}
Let $(\R\times V,J)$ be a cylindrical almost complex manifold as in
Theorem \ref{thm:compact-cyl}, except that
$V$   is  either
\begin{itemize}
\item{ } compact manifold with a pseudo-convex boundary, or
\item{ } can be exhausted by compact domains with pseudo-convex boundaries.
\end{itemize}
Then for every $E>0$,  the space
$\overline\Mc_{g,\mu}(V)\cup\{E(F)\leq E\}$ is compact.
\end{theorem}
The proof of Theorem \ref{thm:compact-cyl} obviously remains valid
here because the pseudo-convex surfaces
  serve as barriers
through which holomorphic curves cannot  escape.

Here are some situations for which  the pseudo-convexity condition
for $\partial V$ is satisfied.
\begin{examples}
\begin{itemize}
\item[\rm a)] If $\dim\,V=3$ and  if the vector field $\Reeb$ is tangent
to the  hypersurface $\Sigma$, then $\Sigma$
   is pseudo-convex.
{\rm In fact, $\wt\Sigma$ is}   Levi-flat  {\rm in this case.  It is
foliated by holomorphic cylinders
$\R\times \gamma$ over  trajectories $\g$ of the vector field
$\Reeb\in  T\Sigma$.}
\item[\rm b)] If $\l_{T\Sigma}=0$,  then $\Sigma$ is pseudo-convex.
{\rm Again, $\wt\Sigma$ is Levi-flat  in this case. This situation
cannot, of course, appear in the contact case.}
\item[\rm c)] If $V= \R^{2n-1}$ and
$\lambda=dz+\frac12\sum\limits_1^{n-1}x_idy_i-y_idx_i$ is the
standard contact form,
then every  geometrically convex (with respect to the standard affine
structure)  hypersurface in $V$ is pseudo-convex.
\end{itemize}
\end{examples}
  \subsection{Degeneration to the Morse--Bott case}
The Compactness Theorem \ref{thm:compact-cyl} is limited to cylindrical almost
complex manifolds with fixed $\omega$, $J$ and $\l$, satisfying
either the Morse
or the Morse--Bott condition of Section \ref{sec:dynamics}. In this section,
we will state a compactness theorem providing a transition from the Morse
to the Morse--Bott case.
We consider a symmetric cylindrical almost complex manifold
$(\bR\times V, J)$ satisfying
the Morse--Bott condition adjusted to a closed form $\omega$ on $V$.
As in Section \ref{sec:dynamics}, we denote by $N_T$  the
submanifolds of $V$  foliated
  by the $T$--periodic trajectories of the vector field $\Reeb$. Let us
fix $T_0>0$. Below we construct a
  special perturbation of $\Reeb$ so that
  all but a finite number of closed orbits on $N_T$ for the periods
$T\leq T_0$ are destroyed while
  the remaining closed orbits become  non-degenerate.
  Let us choose a smooth function $G\co  V \to \bR$ having the following
properties for all periods
  $T\leq T_0$:
  \begin{itemize}
  \item[(1)] along the submanifolds $N_T$, we have $dG(\Reeb )=0$ and
$dG(v)=0$ for every
  vector $v$ normal to $N_T$ with  respect to the natural Riemannian metric
  $g(A, B)=\lambda (A)\cdot \lambda (B)+\omega (A, JB)$, where $A, B\in TV$,
  \item[(2)] the restriction $G\vert_{N_T}$ satisfies the Morse--Bott
condition and the critical
   submanifolds of $G\vert_{N_T}$ consist of finitely many closed
$\Reeb$--orbits.
  \end{itemize}
  For $\e$ small, consider perturbations $\omega_{\varepsilon}$,
   $\Reeb_{\e}$, $J_{\e}$ of $\omega$, $\Reeb$ and $J$, which are
determined by the
   following properties:
  \begin{itemize}
  \item{ } $\omega_{\e}=\omega+\e d(G\l )$
  \item{ } $\Reeb_{\e}\hook\omega_{\e}=0$ and $\lambda (\Reeb_{\e})=1$
  \item{ } $J_{\e}=J$ on $\xi$ and $J_{\e}\frac{\partial }{\partial
t}=\Reeb_{\e}$.
  \end{itemize}
  Writing
  $$\Reeb_{\e}=\Reeb +\e X_{\e},$$
  the vector field $X_{\e}$ satisfies $\l (X_{\e})=0$ so that $X_{\e}\in \xi$.
   Moreover, we have
   \begin{equation}\label{new}
   \begin{split}
   X_\e\hook\omega_\e&=\frac1{\e}\left(
   \Reeb_\e\hook\omega_\e-\Reeb\hook\omega_\e\right)
   \cr
   &=-\frac1{\e}\left(
   \Reeb\hook(\omega+\e dG\wedge\lambda-G\e d\lambda)\right)=dG\,.\cr
   \end{split}
\end{equation}
In particular, the trajectories of $\Reeb$ which
belong to the critical point locus of $G|_{N_T}$,  remain periodic trajectories
of $\Reeb_\e$, while all other closed trajectories of period $\leq
T_0$ got destroyed, as the next lemma
states.
\begin{lemma}
For every $T_0 > 0$, there exists $\e_0 > 0$ so that, if $0 < \e \le \e_0$,
the almost complex structure $J_\e$ satisfies the Morse condition for all
closed $\Reeb_\e$--orbits of period less than $T_0$.
Moreover, every closed $\Reeb_\e$--orbit of period less than $T_0$
corresponds to a critical submanifold of $G\vert_{N_{T}}$ for some $T < T_0$.
\end{lemma}
Note that in general $(X_{\varepsilon})\hook\  d\l \neq 0$, so that
the perturbed cylindrical
almost complex manifold is not symmetric. However, this will not be an issue,
because it is very close to be symmetric. Hence, the proof of Proposition
\ref{prop:automatic-bound} can be repeated almost verbatim. In particular,
$J_\e$--holomorphic curves asymptotic to the closed $\Reeb_\e$--orbits
$\og_1, \ldots, \og_k$ and $\ug_1, \ldots, \ug_l$
corresponding to critical submanifolds of $G\vert_{N_T}$ and representing
a given homology class in $H_2(V,\bigcup \og_i \cup \ug_i)$ have  a
uniformly bounded energy.
On the other hand, observe that a $J_\e$--holomorphic curve $F$, that is
asymptotic
to at least one closed $\Reeb_\e$--orbit which is not a critical
submanifold of some
  $G\vert_{N_T}$, satisfies
$$
E_{\lambda}(F)> T_0.$$
In other words, the $J_\e$--holomorphic curves which satisfy the energy bound
$E(F)\leq T_0$ are  asymptotic only to the non-degenerate
trajectories obtained by the above perturbation.
Let $\g$ be a closed $\Reeb$--orbit on $N_T$ and $\Delta s \in \RR$.
Observe that the gradient flow $\psi^s$, $s\in {\mathbb R}$, of the
gradient vector field $\nabla G$ leaves invariant the  submanifolds
  $N_T$ for $T\leq T_0$. Moreover,
this flow leaves  invariant also the foliation of $N_T$ into
trajectories of $\Reeb$. Let $\g=\gamma(t)$
  be a periodic trajectory
of $\Reeb$ of period $T$ which is not critical for $G$. For $\Delta
s>0 $ we define the
{\it cylindrical gradient trajectory}
$\wt{G}_{\gamma, \Delta s}\co  [0,\Delta s] \times S^1 \to V$
by the formula
\begin{equation*}
\wt{G}_{\g, \Delta s}(s,t) =\psi^s(\gamma(t)),\; s\in[0,\Delta s],\;
t\in S^1\,.
\end{equation*}
Thus the cylinder $\wt{G}_{\gamma, \Delta s}$  is swept by
the gradient trajectories
of the function $G$ starting at  points of  $\g$ and having  length $\Delta s$.
Next, let us define the objects that will be obtained as limits of
$J_\e$--holomorphic maps, when $\e \to 0$. We start with the notion of
a generalized holomorphic building of height 1, with $k'$ sublevels.
The notion is defined by the following four conditions (i)--(iv).
Suppose we are given:
\begin{enumerate}
\item[(i)] $k'$ nodal holomorphic curves
$$
(a_i,f_i;S_i,M_i,D_i,\uZ_i \cup \oZ_i) \, ,\text{where}\:  i = 1, \ldots, k';
$$
\item[(ii)] $k'+1$ collections of cylindrical gradient trajectories
$$
\{ \wt{G}_{\g_{j,i},\Delta s_i} \, ,  j = 1, \ldots, p_i \} \ ,
\text{where}\:  i = 0, \ldots, k',
$$
with $\Delta s_0 = -\infty$, $\Delta s_{k'} = +\infty$ and
$0 \le \Delta s_i < \infty$ for $i = 1, \ldots k'$. We denote by
$\wt{S}_i$ their domains consisting of a collection of cylinders, and by
$\wt{\Gamma}^\pm_i$ the sets of boundary circles corresponding
to their positive and negative punctures;
\item[(iii)]  the number $p^+_i$ of positive punctures of $F_i$ and the number
$p^-_{i+1}$ of negative punctures of $F_{i+1}$ are equal to $p_i$ and for
$i = 1, \ldots, k'$, there are  orientation reversing
diffeomorphisms $\Phi_i \co  \Gamma^+_i \to \wt{\Gamma}^-_i$ and
$\Psi_{i-1} \co  \wt{\Gamma
}^+_{i-1} \to \Gamma^-_i$
which are orthogonal on each boundary component.
\item[(iv)] Gluing the $S^{Z_i}_i$ and $\wt{S}_i$ by means of  the
mappings $\Phi_i$ and the $\Psi_i$,
we obtain  a piecewise-smooth surface
$$
\overline{S} = \wt{S}_0 \mathop{\cup}\limits_{\Psi_0} S^{Z_1}_1 \mathop{\cup}
\limits_{\Phi_1}  \wt{S}_1 \mathop{\cup}\limits_{\Psi_1} \ldots
\mathop{\cup}\limits_{\Phi_{k'}} \wt{S}_{k'} \, .
$$
The last condition in the definition of a generalized holomorphic
building of height $1$
in $\RR \times V$, with $k'$ sublevels i requires that
the maps $\overline{f}_i$, for $i= 1, \ldots, k'$ and
$\wt{G}_{\gamma_{j},i, \Delta s_i}$, for $j = 1, \ldots, p_i$  and $i
= 0, \ldots, k'$,
fit together into a continuous map $\overline{f} \co  \overline{S} \to V$.
\end{enumerate}
Next we  extend this definition to generalized holomorphic buildings of height
$k$, with $k'_i$ sublevels in level $i = 1, \ldots, k$, by
concatenating $k$ generalized holomorphic buildings $F_i$ of height $1$, with
$k'_i$ sublevels, as in Section \ref{sec:height-k}. The stability condition for
generalized holomorphic buildings  means stability of all its
sublevels in the sense of Section
\ref{ssect2.3}.
\begin{figure}[ht!]\small
\centering
\psfrag{a1}{$(a_1,f_1)$}
\psfrag{a2}{$(a_2,f_2)$}
\psfrag{a3}{$(a_3,f_3)$}
\includegraphics[width=.7\textwidth]{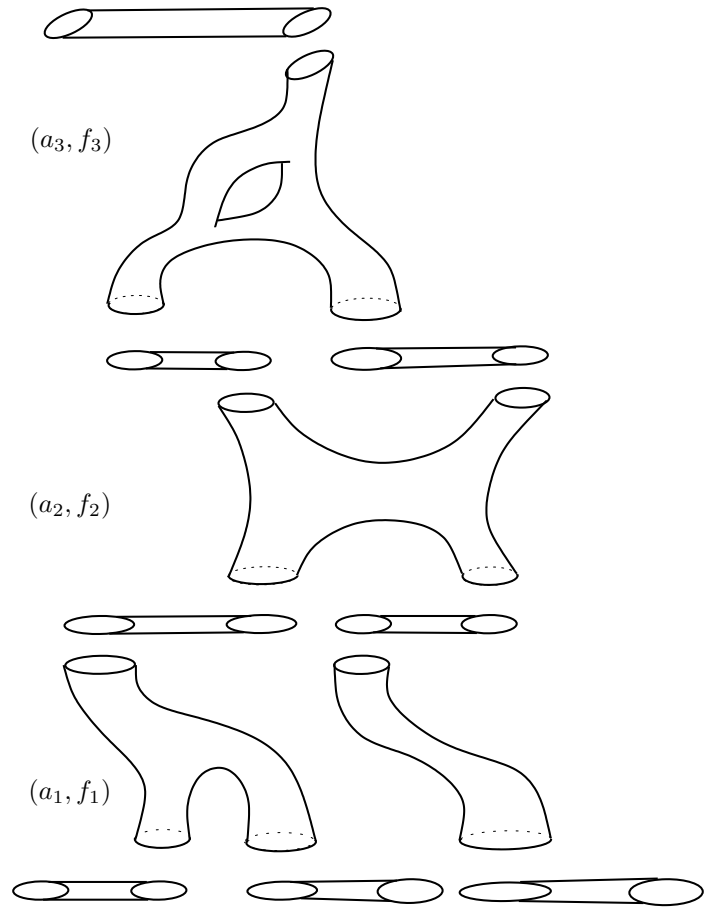}
\caption{Generalized holomorphic building of height $1$ with $3$
sublevels}\label{fig16}
\end{figure}
We now extend  the notion of convergence to a generalized holomorphic building.
Suppose we are given a sequence $(F_m, \Phi_m)$,  with $m\geq 1$,  of
$J_{\e_m}$--holomorphic buildings of height $\le k$, where  $\e_m \to 0$.
The sequence $(F_m,\Phi_m)$ converges to a generalized $J$--holomorphic
building $F$ of height $k$, with $k'_i$ sublevels in level $i = 1, \ldots, k$,
if there exist sequences $M'_m$ of extra sets of
marked points for the curves $(F_m,\Phi_m)$ and a set $M'$ of extra marked
points for $(F,\Phi)$ which have the same cardinality $N$ and which stabilize
the corresponding underlying Riemann surfaces such that the following
conditions
are satisfied.  Assume
$$
(\bS_{F_m},r_{F_m},\Phi_m) = (S_m,j_m,M_m \cup M'_m,D_m,r_m)
$$
and
$$
(\bS_{F},r_{F},\Phi) = (S,j,M \cup M',D,r)
$$
are  the decorated stable nodal Riemann surfaces which underly  the curves
diffeomorphisms  $\varphi_m \co  S^{D,r} \to S^{D_m,r_m}$ with
$\varphi_m(M) = M_m$ and
$\varphi_m(M') = M'_m$ which satisfy the conditions CRS1--CRS3 in the definition
of the Deligne--Mumford convergence of Riemann surfaces and such that,
in addition, the following conditions hold.
\begin{description}
\item{\bf CGHC1}\qua  For every  component $C$ of $S^{D,r} \setminus \bigcup \Gamma_i$
which  is not a cylinder corresponding to a gradient trajectory,
the sequence of projections $f_m \circ \varphi_m |_C \co C \to V$
converges to $f|_C \co  C \to V$ in $C^{\infty}_{\text{loc}}$.
\item{\bf CGHC2}\qua If  $C_{i,j}$  is the union of components of
$S^{D,r} \setminus \bigcup \Gamma_i$ which correspond to the same sublevel
$j = 1, \ldots, k'_i$ of the same level $i = 1, \ldots, k$ of the building $F$,
then there exist sequences $c^{i,j}_m$ for $i = 1, \ldots, k$ and  $j
= 1, \ldots,
k'_i$, and  $m\geq 1$, such that $(a_m \circ \varphi_m - a
-c^{i,j}_m)|_{C_{i,j}}
\to 0$ in $C^{\infty}_{\text{loc}}$.
\end{description}
With these definitions, the compactness theorem can be stated as follows.
\begin{theorem} \label{thm:compact-degen}
Assume $0 < \e_n \leq \e_0$ and  $\e_n \to 0$.
Let $$\bF_n=\left(F_n=(a_n,f_n);S_n,j_n,M_n,
\uZ_n\cup \oZ_n\right) ,\;$$
$$ \oZ_n=
\left\{\left(\oz_1\right)_n, \dots, \left(\oz_{p^+}\right)_n\right\}\;\;,
  \uZ_n=
\left\{\left(\uz_1\right)_n, \dots, \left(\uz_{p^-}\right)_n\right\},$$
be a sequence of  smooth holomorphic curves
in $(W=\RR \times V,J_{\e_n})$ of the
same signature $(g,\mu, p^-,p^+)$ and asymptotic
at the  corresponding punctures to the same orbits in $\Pc_{\e_n}$.
Then there exists a subsequence
that converges to a stable generalized holomorphic building $F$ of height $k$.
\end{theorem}
This theorem is proven in \cite{B} and can be used to compute some symplectic
and
contact invariants in the Morse--Bott case.
  \subsection{Compactness  results in the Relative Symplectic Field Theory}
In the relative Gromov--Witten theory one studies moduli spaces of holomorphic
maps $\left((S,\partial S), j\right)\to\left((W,L),J\right)$ where
$L$ is a totally real submanifold of the middle
dimension. In case $J$ is tamed by a symplectic form $\omega$ one
usually requires, in order to ensure
compactness results, that $L$ is Lagrangian with respect to $\omega$,
ie, $\omega|_L=0$.
In the context of the Symplectic Field Theory we consider   a
symmetric cylindrical
manifold $(\R\times V,J)$ with $\dim\, V=2n-1$ and with $J$ adjusted to a
maximal rank closed $2$--form $\omega$ on $V$. Suppose that $L$ is a
$(n-1)$--dimensional isotropic
submanifold for $\omega$, ie, $\omega|_L=0$, and also integral for
$\xi$, ie, $\l|_L=0$.
In the contact case, when $\omega=d\lambda$, the second condition
implies the first one
and just requires that $L$ is a Legendrian
submanifold for the contact structure $\xi$.  Set $\widehat
L=\R\times L\subset\R\times V$.
Given a Riemann surface $(S,j)$ with boundary $\partial S$ and with
  two sets of punctures $Z=\{z_1,\dots,z_s\}\subset \Int S$ and
$B=\{b_1,\dots,b_t\}\subset\partial S$,
  we consider the moduli
  space of $(j,J)$--holomorphic  maps $(S\setminus Z,\partial
S\setminus B)\to (\R\times V,\widehat L)$ such that
  near interior punctures from $Z$ the maps are asymptotic to periodic
orbits of the Reeb vector field
   $\Reeb$, and near the boundary punctures
  from $Z$,  the maps are asymptotic to chords of $\Reeb$ connecting
two points on $L$,
    either at $+\infty$ or at $-\infty$.
  One can define the moduli spaces of holomorphic buildings of height
$k$ and prove the corresponding
  compactness theorem, similar to Theorem \ref{thm:compact-cyl},
   following the scheme of  the absolute case. In fact, in the real
analytic situation the relative
   compactness  theorem can be formally deduced from the absolute one
by using Gromov's doubling trick, see
   \cite{Gr}. See also  the paper \cite{EES} which is devoted to the
Legendrian contact homology theory, and where
   the necessary compactness results are  proven  in a somewhat different setup.
Suppose now  $(W=\overline W\cup E, J) $ is an almost complex manifold with
  symmetric cylindrical end(s) where $J$ is adjusted to a symplectic
  form $\omega$ on $\overline  W$. Let us consider a Lagrangian
  submanifold $\Lambda\subset \overline W$ which intersects
  the boundary $V=\partial\overline W$ transversally
   along a submanifold $L=\overline W\cup\Lambda$, integral for the distribution
   $\xi=\{\l=0\}$ on $V$. Then $\Lambda$ can be extended to the end
$E$ as the cylindrical manifold
    $\widehat L=\R\times L$.
   We abbreviate  $\widehat \Lambda=\Lambda\cup\widehat L$ and
consider the moduli space of $(j,J)$--holomorphic
    maps $(S\setminus Z,\partial S\setminus B)\to (W,\widehat
\Lambda)$ such that
  near interior punctures from $Z$ the maps are asymptotic  to
periodic orbits of $\Reeb$,
  and near the boundary punctures
  from $Z$ the maps are asymptotic to chords of $\Reeb$ connecting two
points on $L$,   at one of the ends of $W$.
  The definition of the moduli spaces of holomorphic buildings of
height $k_-|1|k_+$, the formulation  and
  the  proof of  the relative  analogue  of Theorem \ref{thm:compact-ends} can
  be done following the same scheme, as in the
  absolute case.
\appendix
\section{Appendix : Asymptotic convergence estimates} \label{ap:estimates}
We shall  first describe the structure of the manifold $V$ near a periodic
orbit of Morse--Bott type.
\begin{lemma}\label{lm:Bott-neighborhood}
Assuming the Morse--Bott situation,  let $N$ denote a
  component of the set $N_T\subset V$ covered by periodic orbits
of period $T$
of the vector field $\Reeb$. Let  $\g$ be one of the orbits from $N$. Then
\begin{itemize}
\item[\rm a)] if $T$ is the minimal period of $\g$ then there exists a neighborhood
$U\supset\g$ in $V$ such that
$U\cap N$ is invariant under the flow of $\Reeb$ and  one  finds coordinates
$$(\vartheta,x_1,\dots, x_{n-1}, y_1,\dots,y_{n-1})$$
such that $$N=\{x_1,\dots, x_p=0,y_1,\dots, y_q=0\},\;\;\hbox{for}\;\;
0\leq p,q\leq n-1,$$ $$\Reeb|_N=\frac{\partial}{\partial\vartheta}\, ,$$  and
$$\omega|_N=\omega_0|_N\;\;\hbox{where}\;\;
\omega_0=\sum\limits_1^{n-1}dx_i\wedge dy_i\,.$$
\item[\rm b)] if  $\gamma$ is  a $m$--multiple of a trajectory $P$ of
a minimal period  $\frac Tm$
there exists a tubular neighborhood $ \wt U$ of $P$ such that
its $m$--multiple cover $ U$ together
with all the structures induced by the covering map $\wt U\to U$
from the corresponding objects on $U$
satisfy the properties of the part a).
\end{itemize}
\end{lemma}
\begin{proof} In the case a) the orbit $\g$ has a neighborhood $U$
such that $U\cap N$ is fibrated by trajectories of $\Reeb$.
  since $\Reeb\hook \omega = 0$ and $d\omega = 0$, the $2$--form $\omega|_N$
descends to the quotient by $S^1$. By the Morse--Bott assumption the rank of
of the form induced on this quotient  is constant, and hence  ${N}/S^1$ is
foliated by isotropic submanifolds. The statement then follows using an
appropriate version of Darboux theorem. In the situation of b) by taking
the $m$--multiple cover of a neighborhood of $P$ we return to the
situation considered
in a).
\end{proof}
Let $N$ be,   as in the above lemma, a component consisting of
periodic orbits of $\Reeb$
of period $T$. Let $\g$ is one of the periodic orbits from $N$.
Let $F = (a,f) \co  [-R,+R] \times S^1 \rightarrow \RR \times V$ be a
$J$--holomorphic cylinder in a neighborhood of  $\gamma$. If the
minimal period of $\g$
is equal to $T/m$ then
by taking the $m$--folded cover of a neighborhood of $\g$ and lifting
there our holomorphic
map we are in
the situation in which  $T$ is the minimal period of $\g$, and hence
a neighborhood of
$\g$ in $N$ is fibrated
by the closed orbits of $\g$.
We  introduce in a neighborhood $U$ of  $\g$ the local coordinates
  $(\vartheta, z_{in}, z_{out})
\in \RR/\ZZ \times \RR^k \times \RR^{2n-2-k}$,
where $k =2n-2-(p+q)$ and
\begin{equation*}
z_{in}=(x_{p+1},\dots,x_{n-1},y_{q+1},\dots,y_{n-1}),\qquad  z_{out}
=(x_1,\dots,x_p,y_1,\dots,y_q).
\end{equation*}
We will also set $z = (z_{in}, z_{out})$ and $I_0 = [-R,+R]$,
$\theta_0 = I_0 \times
S^1$.
\begin{lemma}  \label{lem:CR}
Near a closed $\Reeb$--orbit, the Cauchy--Riemann equations can be
written as follows,
\begin{eqnarray}
z_s + Mz_t + S z_{out} &=& 0 \label{eqn:z} \\
a_s - T \vartheta_t + B z_{out} + B' z_t  &=& 0 \label{eqn:as} \\
a_t + T \vartheta_s + C z_{out} + C' z_t &=& 0 . \label{eqn:at}
\end{eqnarray}
\end{lemma}
\begin{proof}
The Cauchy--Riemann equations are given by
$$
\frac{\partial F}{\partial s} + J \frac{\partial F}{\partial t} = 0.
$$
Let us first extract the $z$--components of this equation. The first term
gives $z_s$, while the second term gives 3 contributions:
\begin{enumerate}
\item $M z_t$, the $z$--components of $J z_t$, where $M = \pi_z \circ J|_z$ and
$\pi_z$ is the projection on the $z$ coordinates,
\item $\pi_z (J \frac{\partial}{\partial \vartheta}) \vartheta_t$. Since
$\frac{\partial}{\partial \vartheta} = T \Reeb$ along $N_T$, this term has
the form $S' z_{out}$.
\item $\pi_z (J \frac{\partial}{\partial t}) a_t$. Since $\pi_z \Reeb
= 0$ along
$N_T$, this term has the form $S'' z_{out}$.
\end{enumerate}
Combining all terms, we obtain the equation
$$
z_s + M z_t + S z_{out} = 0
$$
where $S = S' + S''$. \\
Let us now extract the $a$ component of the Cauchy--Riemann equation. The first
term gives $a_s$, while the second term gives 3 contributions:
\begin{enumerate}
\item $\pi_a \circ J z_t = B' z_t$, where $\pi_a$ is the projection on the
$a$ coordinate.
\item $\pi_a(J \frac{\partial}{\partial \vartheta}) \vartheta_t =
-T \vartheta_t + B z_{out}$, since $J \frac{\partial}{\partial \vartheta}
= -T \frac{\partial}{\partial a}$ along $N_T$.
\item $\pi_a(J \frac{\partial}{\partial t}) a_t = 0$.
\end{enumerate}
Combining all terms, we obtain the equation
$$
a_s - T \vartheta_t + B z_{out} + B' z_t = 0 .
$$
Finally, we apply $J$ to the Cauchy--Riemann equation and extract the $a$
component from the resulting equation. The second term gives $- a_t$,
while the first term gives 3 contributions:
\begin{enumerate}
\item $\pi_a(J z_s) = C_1 z_{out} + C' z_t$ using the $z$--components of the
Cauchy--Riemann equations.
\item $\pi_a(J \frac{\partial}{\partial \vartheta}) \vartheta_s =
-T \vartheta_t + C_2 z_{out}$, since $J \frac{\partial}{\partial \vartheta}
= -T \frac{\partial}{\partial a}$ along $N_T$.
\item $\pi_a(J \frac{\partial}{\partial t}) a_s = 0$.
\end{enumerate}
Combining all terms and changing the sign, we obtain the equation
$$
a_t + T \vartheta_t + C z_{out} + C' z_t = 0
$$
where $C = C_1 + C_2$.
\end{proof}
Define the linear operators
$$
A(s) \co  H^1(S^1,\RR^{2n-2}) \subset L^2(S^1,\RR^{2n-2})
\rightarrow L^2(S^1,\RR^{2n-2})
$$
by the formula
$$
\big( A(s) z \big) (t) = -M \big( u(s,t) \big) z_t(t) -
S \big( u(s,t) \big) z_{out}(t) ,
$$
where we have abbreviated $u(s, t)=\bigl(a(s,t), f(s,t)\bigr)$.
Then  equation \eqref{eqn:z} becomes $z_s(s,\cdot) = A(s)z(s,\cdot)$.
Note that the explicit expression for this operator depends on the
$J$--holomorphic
map $F = (a,f)$ through the matrices $M$ and $S$.
If we substitute instead $\vartheta(s,t) = \vartheta(s_0,0) + t$,
$a(s,t) = Ts$,
$z_{out}(s,t) = 0$ and $z_{in}(s,t) = z_{in}(s_0,0)$, for some $s_0 \in I_0$,
we obtain another operator $\wt A(s)$.
We will denote the limit $\lim_{s \rightarrow \infty} \wt A(s)$ by $A_0$.
We write $\big( A_0 u \big) (t) = -M_0(t) u_t(t)  -S_0(t) u_{out}(t)$.
This operator corresponds to the linearized Cauchy--Riemann equation along the
closed $\Reeb$--orbit $\g_{s_0}(t) = (\vartheta(s_0,0)+t,z_{in}(s_0,0),0)$.
Revisiting the proof of Lemma \ref{lem:CR}, we can see that $M_0 S_0 z_{out}$
is the linearization of the $z$--component of the vector field $\Reeb$
along $\g_{s_0}$.
Hence, the matrices $M_0(t), S_0(t)$ are in the symplectic algebra and
the operator $A_0$ is self-adjoint with respect to the inner product
$$
\langle u,v \rangle_0 = \int_0^1 \langle u,-J_0 M_0 v \rangle dt
$$
where $J_0$ is the standard complex structure on $\RR^{2n-2}$ and $\langle
\cdot,\cdot \rangle = \omega_0( \cdot, J_0 \cdot)$.
The kernel of $A_0$ is independent of $s_0$ and is generated by
constant loops with values in the tangent space to $N$. Let $P_0$ be the
orthonormal projection to $\ker A_0$ with respect to
$\langle \cdot, \cdot \rangle_0$,
and $Q_0 = I - P_0$. The operator $Q_0$ clearly has the
following properties: $(Q_0 z)_t = z_t$, $(Q_0 z)_s = Q_0 z_s$,
$(Q_0 z)_{out} = z_{out}$ and $Q_0 A_0 = A_0 Q_0$.
We will first obtain some estimates for the decaying rate of $z_{out}$.
Abbreviate $g_0(s): = \frac12 \| Q_0 z(s) \|_0^2$ and
$\beta_0(s) = (\vartheta(s_0,0)-\vartheta(s,0),z_{in}(s_0,0)-z_{in}(s,0))$.
\begin{lemma}  \label{zout}
There exist $\delta > 0$ and ${\bar \beta} > 0$ such that, if
$$
\sup_{(s,t) \in \theta_0} | \partial^\alpha z_{out}(s,t) |
\le \delta
$$
for multi-indices $\alpha$ with $|\alpha| \le 2$,  and
\begin{eqnarray*}
\sup_{(s,t) \in \theta_0} | \partial^\alpha z_{in}(s,t) |
&\le& \delta  \\
\sup_{(s,t) \in \theta_0} | \partial^\alpha (\vartheta(s,t) - t) |
&\le& \delta
\end{eqnarray*}
for those multi-indices $\alpha$ satisfying  $0 < |\alpha| \le 2$,
then for $s \in I_0$ satisfying $|\beta_0(s)| \le {\bar \beta}$, we have
$$
g''_0(s) \ge c_1^2 g_0(s),
$$
where $c_1 > 0$ is a constant independent of $s_0$.
\end{lemma}
\proof
Clearly,
$$
g_0''(s) \ge \langle Q_0 z_{ss},Q_0 z \rangle_0 .
$$
Let us compute the right hand side. First,
\begin{eqnarray*}
Q_0 z_s &=& Q_0 A(s) z(s) \\
&=& Q_0 A_0 z + Q_0 [ A(s) - A_0 ] z \\
&=& Q_0 A_0 z + Q_0 [ \Delta_0 z_t + \hat{\Delta}_0 z_{out}
+ (\tilde{\Delta}_0 \beta_0) z_t + ({\bar \Delta}_0 \beta_0) z_{out}] \\
&=& A_0 Q_0 z + Q_0 \Delta_0 (Q_0 z)_t + Q_0 \hat{\Delta}_0 (Q_0 z)_{out}
+ Q_0 (\tilde{\Delta}_0 \beta_0)(Q_0 z)_t \\
&& + Q_0 ({\bar \Delta}_0 \beta_0) (Q_0 z)_{out} ,
\end{eqnarray*}
where
\begin{eqnarray*}
\Delta_0 &=& M(\vartheta(s,0)+t,z_{in}(s,0),0) - M(\vartheta(s,t),z(s,t)) \\
\hat{\Delta}_0 &=& S(\vartheta(s,0)+t,z_{in}(s,0),0) -
S(\vartheta(s,t),z(s,t)) ,
\end{eqnarray*}
and  where the matrices $\tilde{\Delta}_0(s,t), {\bar \Delta}_0 (s,t)
$ are defined via the mean value theorem
applied to $M$ between points
$(\vartheta(s_0,0)+t,z_{in}(s_0,0))$ and
$(\vartheta(s,0)+t,z_{in}(s,0))$, and hence  we have
\begin{eqnarray*}
\tilde{\Delta}_0 \beta_0&=& M_0 - M(\vartheta(s,0)+t,z_{in}(s,0),0) \\
{\bar \Delta}_0 \beta_0 &=& S_0 - S(\vartheta(s,0)+t,z_{in}(s,0),0) .
\end{eqnarray*}
The expressions $\Delta_0$ and $\hat{\Delta}_0$ contain the dependence
in $z_{out}$ so that
\begin{eqnarray*}
\sup_{(s,t) \in \theta_0} | \partial^\alpha \Delta_0(s,t)| &\le& C \delta  \\
\sup_{(s,t) \in \theta_0} | \partial^\alpha \hat{\Delta}_0(s,t)| &\le&
  C  \delta ,
\end{eqnarray*}
for multi-indices $\alpha$ with $|\alpha| \le 1$.
On the other hand, the expressions $\tilde{\Delta}_0$ and ${\bar \Delta}_0$
contain the dependence in $z_{in}$ and $\vartheta$. Therefore,
\begin{eqnarray*}
\sup_{(s,t) \in \theta_0} | \partial^\alpha \tilde{\Delta}_0(s,t)| &\le& C \\
\sup_{(s,t) \in \theta_0} | \partial^\alpha {\bar \Delta}_0(s,t)| &\le& C ,
\end{eqnarray*}
for multi-indices $\alpha$ with $|\alpha| \le 1$.
Taking the derivative once more, we obtain
\begin{eqnarray*}
Q_0 z_{ss} &=& A_0(Q_0z)_s + Q_0(\frac{\partial}{\partial s} \Delta_0)(Q_0 z)_t
+ Q_0 \Delta_0(Q_0 z_s)_t  \\
&& + Q_0(\frac{\partial}{\partial s} \hat{\Delta}_0)(Q_0 z)_{out} + Q_0
\hat{\Delta}_0 ((Q_0 z)_{out})_s \\
&& + Q_0(\frac{\partial}{\partial s}\tilde{\Delta}_0)\beta_0 (Q_0 z)_t
+ Q_0(\tilde{\Delta}_0 \frac{d}{ds}\beta_0)(Q_0 z)_t + Q_0(\tilde{\Delta}_0
\beta_0)(Q_0 z_s)_t \\
&& + Q_0 (\frac{\partial}{\partial s}{\bar \Delta}_0) \beta_0 (Q_0 z)_{out}
+ Q_0({\bar \Delta}_0 \frac{d}{ds} \beta_0)(Q_0 z)_{out} \\
&& + Q_0({\bar \Delta}_0 \beta_0) (Q_0 z_s)_{out} .
\end{eqnarray*}
Taking the inner product with $Q_0 z$, we obtain
\begin{eqnarray*}
\langle Q_0 z_{ss},Q_0 z \rangle_0 &=& \langle Q_0 z_s, A_0(Q_0z) \rangle_0
+ \langle (\frac{\partial}{\partial s} \Delta_0)(Q_0 z)_t, Q_0 z \rangle_0 \\
&& + \langle \Delta_0(Q_0 z_s)_t, Q_0 z \rangle_0  \\
&& + \langle (\frac{\partial}{\partial s} \hat{\Delta}_0)(Q_0 z)_{out},
Q_0 z \rangle_0
+ \langle \hat{\Delta}_0 ((Q_0 z)_{out})_s, Q_0 z \rangle_0 \\
&& + \langle (\frac{\partial}{\partial s}\tilde{\Delta}_0)\beta_0 (Q_0 z)_t,
Q_0 z \rangle_0
+ \langle (\tilde{\Delta}_0 \frac{d}{ds}\beta_0)(Q_0 z)_t, Q_0 z \rangle_0 \\
&& + \langle (\tilde{\Delta}_0 \beta_0)(Q_0 z_s)_t, Q_0 z \rangle_0 \\
&& + \langle (\frac{\partial}{\partial s}{\bar \Delta}_0) \beta_0
(Q_0 z)_{out},
Q_0 z \rangle_0 + \langle ({\bar \Delta}_0 \frac{d}{ds} \beta_0)
(Q_0 z)_{out}, Q_0 z \rangle_0  \\
&& + \langle ({\bar \Delta}_0 \beta_0) (Q_0 z_s)_{out}, Q_0 z \rangle_0 .
\end{eqnarray*}
Let us denote the 11 terms of the right hand side by
$T_1, \ldots, T_{11}$.
Substituting  $Q_0 z_s$ by its value in $T_1$ we find
\begin{eqnarray*}
T_1 &=& \| A _0 Q_0 z \|_0^2 + \langle Q_0 \Delta_0 Q_0 z_t, A_0 Q_0
z \rangle_0
+ \langle Q_0 \hat{\Delta}_0  (Q_0 z)_{out}, A_0 Q_0 z \rangle_0 \\
&& + \langle Q_0 (\tilde{\Delta}_0 \beta_0)Q_0 z_t, A_0 Q_0 z \rangle_0
+ \langle Q_0 ({\bar \Delta}_0 \beta_0) (Q_0 z)_{out}, A_0 Q_0 z \rangle_0 .
\end{eqnarray*}
By integration by parts in $T_3$ and $T_8$,
\begin{eqnarray*}
T_3 &=& \langle \Delta_0 (Q_0 z_s)_t, Q_0 z \rangle_0 \\
&=& \int_0^1 \langle (Q_0 z_s)_t, -\Delta^*_0 J_0 M Q_0 z \rangle dt \\
&=& - \int_0^1 \langle Q_0 z_s, (-\frac{\partial}{\partial t}
\Delta^*_0  J_0M) Q_0 z \rangle dt
- \int_0^1 \langle Q_0 z_s, -\Delta^*_0 J_0 M Q_0 z_t \rangle dt .
\end{eqnarray*}
Similarly,
$$
T_8 = -\int_0^1 \langle Q_0 z_s, (-\frac{\partial}{\partial t}
(\tilde{\Delta}_0 \beta_0) J_0 M) Q_0 z \rangle dt
-\int_0^1 \langle Q_0 z_s, -(\tilde{\Delta}_0 \beta_0)
J_0 M Q_0 z_t \rangle dt .
$$
Applying  the Cauchy--Schwarz inequality to all terms $T_i$ and taking
into account the bounds for $\Delta_0$, $\hat{\Delta}_0$ and for
$\tilde{\Delta}_0$, ${\bar \Delta}_0$, we obtain
\begin{eqnarray*}
T_1 &\ge&  \| A_0 Q_0 z \|_0^2 - c \delta \| Q_0 z_t \|_0 \|A_0 Q_0 z \|_0
- c \delta \| Q_0 z \|_0 \| A_0 Q_0 z \|_0 \\
&& - c |\beta_0| \, \| Q_0 z_t \|_0
\| A_0 Q_0 z \|_0 - c |\beta_0| \, \| Q_0 z \|_0 \| A_0 Q_0 z \|_0.  \\
T_2 &\ge& -c \delta \| Q_0 z_t \|_0 \| Q_0 z\|_0 . \\
T_3 &\ge& -c \delta \| Q_0 z_s \|_0 \| Q_0 z \|_0 -c \delta
\| Q_0 z_s \|_0 \| Q_0 z_t \|_0.   \\
T_4 &\ge& -c \delta \| Q_0 z \|_0^2.  \\
T_5 &\ge& -c \delta \| Q_0 z_s \|_0 \| Q_0 z \|_0 .\\
T_6 &\ge& -c |\beta_0| \, \| Q_0 z_t \|_0 \| Q_0 z \|_0.  \\
T_7 &\ge& -c \delta \| Q_0 z_t \|_0 \| Q_0 z\|_0 .\\
T_8 &\ge& -c |\beta_0| \, \| Q_0 z_s \|_0 \| Q_0 z \|_0 -c |\beta_0| \,
\| Q_0 z_s \|_0  \| Q_0 z_t \|_0 . \\
T_9 &\ge& -c \, |\beta_0| \| Q_0 z \|_0^2.  \\
T_{10} &\ge& -c \delta \| Q_0 z \|_0^2 . \\
T_{11} &\ge& -c |\beta_0| \, \| Q_0 z_s \|_0 \| Q_0 z \|_0 .
\end{eqnarray*}
Using the expression for $Q_0 z_s$ we find
$$
\| Q_0 z_s \|_0 \le \| A_0 Q_0 z \|_0 + c \delta \| Q_0 z_t \|_0
+ c \delta \| Q_0 z \|_0 + c |\beta_0| \, \| Q_0 z_t \|_0 + c |\beta_0| \,
\| Q_0 z \|_0 .
$$
On the other hand, it is clear from the definition of $Q_0$ that
$$
\| A_0 Q_0 z \|_0 \ge c_1 (\| (Q_0 z)_t \|_0^2 + \| Q_0 z \|_0^2)^\frac12 .
$$
Using the last  2 inequalities to eliminate $Q_0 z$, $Q_0 z_s$ and $Q_0 z_t$
from the estimates for the $T_i$, we end up with
$$
\langle Q_0 z_{ss}, Q_0 z \rangle_0 \ge (1 - c \delta - c |\beta_0|)
\| A_0 Q_0 z \|_0^2 .
$$
Therefore, if $\delta$ and ${\bar \beta}$ are
sufficiently small, then for $|\beta_0(s)| < {\bar \beta}$
we will have
$$
\langle Q_0 z''(s), Q_0 z(s) \rangle_0 \ge \frac12 \| A_0 Q_0 z(s) \|_0^2 .
$$
 From this we deduce the desired estimate
\begin{align*}
g''_0(s) &\ge \langle Q_0 z''(s), Q_0 z(s) \rangle_0 \\
&\ge \frac12 \| A_0 Q_0 z(s) \|_0^2 \\
&\ge \frac{c_1^2}2 \| Q_0 z \|_0^2 = c_1^2 g_0(s) .\tag*{\qed}
\end{align*}

Define $s_1 = \sup\{ s \in I_0 \, | \, |\beta_0(s')| \le
{\bar \beta} \textrm{ for all } s' \in [s_0,s] \}$.
Then Lemma \ref{zout} implies the following estimate.
\begin{corollary}
\begin{equation} \label{eqn:g0}
g_0(s) \le  \max(g_0(s_0),g_0(s_1))
\frac{\cosh(c_1(s-\frac{s_0+s_1}2))}{\cosh(c_1\frac{s_0-s_1}2)}
\end{equation}
for $s \in [s_0, s_1]$.
\end{corollary}
\begin{proof} Let us assume for determinacy that the maximum
$\max(g_0(s_0),g_0(s_1))$ is achieved
at $s_0$, so that $g_0(s_0)\geq g_0(s_1)$. Let us
denote the right hand side in (\ref{eqn:g0}) by $g_1(s)$.
This function satisfies $g''_1(s) = c_1^2 g_1(s)$, $g_1(s_0) = g_0(s_0)$ and
$g_1(s_1)\geq g_0(s_1)$. Therefore, the difference $g(s) = g_0(s) - g_1(s)$
satisfies $g''(s) \ge c_1^2 g(s)$ for $s \in [s_0,s_1]$,  vanishes
at $s_0$ and is non-positive at $s_1$.
  The differential inequality implies that $g(s)$ cannot have a positive
local maximum for $s_0 < s < s_1$. Indeed, if $g$ is positive on the
open  sub-interval
$\Delta\subset (s_0,s_1)$ and vanishes on its boundary, then
there is point $a\in\Delta$ at which  $g''(a)<0$, while $g(a)>0$,
which contradicts the differential inequality.
  Hence, $g$ is non-positive on $[s_0,s_1]$.
  \end{proof}
Next we  derive some estimates for $z_{in}$.
Let $e$ be a unit vector in $\RR^{2n-2}$ with $e_{out} = 0$.
\begin{lemma} \label{zin}
Under the assumptions of Lemma \ref{zout} and for $s \in [s_0,s_1]$, we have
$$
|\langle z(s),e \rangle_0 - \langle z(s_0),e \rangle_0|
\le \frac{4d}{c_1} \max(\| Q_0z(s_0)\|_0, \| Q_0z(s_1) \|_0) .
$$
\end{lemma}
\proof
The inner product of the Cauchy--Riemann equation (\ref{eqn:z}) with $e$ gives
$$
\frac{d}{ds}\langle z,e \rangle_0 + \langle M z_t,e \rangle_0
+ \langle Sz_{out},e \rangle_0 = 0 .
$$
But we have
\begin{eqnarray*}
\langle Mz_t,e \rangle_0 &=& \int_0^1
\langle M(Q_0z)_t,-J_0 M_0 e \rangle dt \\
&=& \int_0^1 \langle Q_0 z, \frac{d}{dt} (M^* J_0 M_0) e \rangle dt ,
\end{eqnarray*}
so that
$$
| \langle Mz_t,e \rangle_0| \le d_1 \| Q_0 z \|_0 .
$$
Similarly,
$$
\langle S(Q_0 z)_{out},e \rangle_0 =
\int_0^1 \langle (Q_0 z)_{out}, S^*(-J_0) M_0 e \rangle dt ,
$$
so that
$$
|\langle Sz_{out},e \rangle_0| \le d_2 \|Q_0 z\|_0 .
$$
Therefore,
$$
\langle z(s),e \rangle_0 - \langle z(s_0),e \rangle_0 \le
d\int_{s_0}^s \| Q_0 z(\sigma)\|_0 d\sigma
$$
for $s \in I_0$. By equation (\ref{eqn:g0}),
$$
\| Q_0 z(\sigma) \|_0 \le \max(\| Q_0 z(s_0) \|_0,\| Q_0 z(s_1)\|_0)
\sqrt{\frac{\cosh(c_1(\sigma-\frac{s_0+s_1}2))}{\cosh(c_1\frac{s_0-s_1}2)}}
$$
for $\sigma \in [s_0,s_1]$. Hence, substituting in the integral and using
the fact that $\sqrt{\cosh u} < \sqrt{2} \cosh \frac{u}2$, we obtain
$$
|\langle z(s),e \rangle_0 - \langle z(s_0),e \rangle_0| \le
\frac{4d}{c_1} \max(\| Q_0z(s_0)\|_0,\| Q_0z(s_1)\|_0) ,
$$
since
$$
\sqrt{2} \frac{\sinh(c_1\frac{s_0-s_1}4)}
{\sqrt{\cosh(c_1\frac{s_0-s_1}2)}} \le 1 .\eqno{\qed}
$$

Next, we shall derive  estimates for the derivatives of $z$.
\begin{lemma}  \label{z}
There exist $\delta > 0$ and ${\bar \beta} > 0$ such that, if
$$
\sup_{(s,t) \in \theta_0} | \partial^\alpha z_{out}(s,t) |
\le \delta
$$
for multi-indices $\alpha$ with $|\alpha| \le k + 2$ and
\begin{eqnarray*}
\sup_{(s,t) \in \theta_0} | \partial^\alpha z_{in}(s,t) |
&\le& \delta  \\
\sup_{(s,t) \in \theta_0} | \partial^\alpha (\vartheta(s,t) - t) |
&\le& \delta
\end{eqnarray*}
for multi-indices $\alpha$ with $0 < |\alpha| \le k + 2$, then
\begin{equation}  \label{eqn:zest}
\begin{split}
\| \partial^\alpha &z(s)\|_0\\
&\le  C_\alpha
\max\limits_{|\alpha'| \le |\alpha|}(\|Q_0 \partial^{\alpha'} z(s_0)\|_0,
\|Q_0 \partial^{\alpha'} z(s_1)\|_0)
\sqrt{\frac{\cosh(c_1(s-\frac{s_0+s_1}2))}{\cosh(c_1\frac{s_0-s_1}2)}} ,
\end{split}
\end{equation}
for all $s \in [s_0,s_1]$ and for every multi-index $\alpha$ with $1 \le
|\alpha| \le k$.
\end{lemma}
\begin{proof}
The Cauchy--Riemann equation (\ref{eqn:z}) can be written as
\begin{eqnarray*}
z_s &=& A(s) z \\
&=& A_\infty z + {\bar \Delta} z_t + \Delta z_{out},
\end{eqnarray*}
where $\Delta = S_0 - S$ and ${\bar \Delta} = M_0 - M$.
Applying the projection $Q_0$ to this equation, we obtain
for $w=\Q_0 z$,
$$
w_s = A_0 w + Q_0 {\bar \Delta} w_t + Q_0 \Delta w_{out}.
$$
where $w = Q_0 z$. Let $W$ be the vector obtained by catenating
$(\frac{\partial}{\partial s})^i (A_0)^j w$ for $0 \le i,j \le k$.
Then $W$ satisfies an equation of the same type,
$$
W_s = {\cal A}_0 W + {\cal Q}_0 \tilde{\Delta} W_t +
{\cal Q}_0 \hat{\Delta} W_{out},
$$
where ${\cal A}_0 = {\rm diag}(A_0, \ldots, A_0)$ and
${\cal Q}_0 = {\rm diag}(Q_0, \ldots, Q_0)$. Therefore, using the same
estimates as in Lemma \ref{zout}, we obtain
$$
\| W(s) \|_0 \le C \max(\|W(s_0)\|_0,\|W(s_1)\|_0)
\sqrt{\frac{\cosh(c_1(s-\frac{s_0+s_1}2))}{\cosh(c_1\frac{s_0-s_1}2)}} .
$$
Next we estimate $P_0 z$ and its derivatives. Applying $P_0$ to the
Cauchy--Riemann equation (\ref{eqn:z}), we get
$$
(P_0 z)_s = P_0 {\bar \Delta}(Q_0 z)_t + P_0 \Delta (Q_0 z)_{out} .
$$
We can apply $(\frac{\partial}{\partial s})^i$, for $i = 0, \ldots, k-1$ to
this equation, and express the derivatives of $P_0 z$ in terms of
components of $W$,  to obtain  the desired estimate.
\end{proof}
We now derive  estimates for $\vartheta$.
\begin{lemma}  \label{theta}
Under the assumptions of Lemma \ref{z} with $k = 1$,
$$
\int_0^1 |\vartheta(s,\tau) - \vartheta(s_0,\tau)| d\tau
\le C  \max\limits_{|\alpha| \le 1}(\| Q_0 \partial^\alpha z(s_0)\|_0,
\| Q_0 \partial^\alpha z(s_1) \|_0)
$$
for all $s \in [s_0,s_1]$.
\end{lemma}
\proof
Consider the Cauchy--Riemann equations (\ref{eqn:as}) and (\ref{eqn:at})
for $a$ and $\vartheta$:
$$
\left\{ \begin{array}{rcl}
a_s - T \vartheta_t &=& -B (Q_0 z)_{out} - B' (Q_0 z)_t   \\
a_t + T \vartheta_s &=& -C (Q_0 z)_{out} - C' (Q_0 z)_t  .
\end{array} \right.
$$
If $s \in [s_0,s_1]$, the right hand side is bounded in norm as in
equation (\ref{eqn:zest}). Therefore, integrating the second equation over $t$,
we obtain
$$
\int_0^1 \vartheta_s dt \le C \max\limits_{|\alpha| \le 1}
(\| Q_0 \partial^\alpha z(s_0) \|_0, \| Q_0 \partial^\alpha z(s_1) \|_0)
\sqrt{\frac{\cosh(c_1(s-\frac{s_0+s_1}2))}{\cosh(c_1\frac{s_0-s_1}2)}} .
$$
Integrating over $s$, as in Lemma \ref{zin}, we find
$$
\int_0^1 |\vartheta(s,\tau) - \vartheta(s_0,\tau)| d\tau \le
C'  \max\limits_{|\alpha| \le 1}(\| Q_0 \partial^\alpha z(s_0) \|_0,
\| Q_0 \partial^\alpha z(s_1) \|_0) .\eqno{\qed}
$$

We are now in position to prove Proposition \ref{prop:twist}.
\begin{proof}[Proof of Proposition \ref{prop:twist}]
Let $\hbar$ is chosen as in Lemma \ref{lm:bubble}. Suppose also that
the neighborhood $U$ is chosen
so small that it satisfies the condition of Lemma
\ref{lm:Bott-neighborhood} and besides that it contains no
periodic orbits of $\Reeb$ of period $<2T$, other than those which
form the component $N=N_T\ni\gamma$.
By contradiction, assume that there exists a sequence $F_n = (a_n,f_n) :
[-n,+n] \times S^1 \rightarrow \RR \times V$ of holomorphic cylinders,
satisfying $E(F_n) \le E_0$ and $E_\omega(F_n) \le \hbar$, and a sequence
$c_n\to\infty,\,c_n<n$
such that
$f_n(s_n,t) \notin B_\varepsilon(f_n(0,t))$ for some $s_n \in [-k_n,k_n]$,
$k_n=n-c_n$.
By Lemma \ref{lm:bubble}, the gradient of $F_n$ is uniformly bounded on each
compact subset, otherwise we would obtain a bubble with the
$\omega$--energy exceeding $\hbar$.
Hence, by Ascoli--Arzela, we can extract a subsequence converging to a
holomorphic cylinder $F =\co \RR \times S^1 \rightarrow \RR \times V$,
which is necessarily a
trivial vertical cylinder over an orbit $\hat\g\in N$. Indeed, both
asymptotic limits
$\hat\g^\pm$ of this cylinder have to belong to $N_T$. This  forces
the equality $E_\omega(F)=0$ which then implies
that $F$ is a trivial vertical cylinder.
Let us show that
\begin{equation}  \label{eqn:sup1}
\sup_{(s,t) \in [-k_n,k_n] \times S^1} |\partial^\alpha z_{out,n}(s,t)| \to 0
\end{equation}
for multi-indices $\alpha$ with $|\alpha| \ge 0$ and
\begin{eqnarray}
\sup_{(s,t) \in [-k_n,k_n] \times S^1} | \partial^\alpha z_{in,n}(s,t) |
&\to& 0  \label{eqn:sup2} \\
\sup_{(s,t) \in [-k_n,k_n] \times S^1} | \partial^\alpha
(\vartheta_n(s,t) - t) |
&\to& 0 \label{eqn:sup3}
\end{eqnarray}
for multi-indices $\alpha$ with $|\alpha| \ge 1$, when $n \to \infty$.
If this were not the case, we could translate the coordinates to center them on
a sequence of points violating one of these properties.
The sequence of cylinders obtained in this way
would converge, as shown above, but the limit could not be a vertical cylinder,
giving a contradiction.
Hence, for $n$ sufficiently large, the suprema in equations (\ref{eqn:sup1}),
(\ref{eqn:sup2}) and (\ref{eqn:sup3}) for $|\alpha| \le 3$
(and $|\alpha|\geq 1$ in cases (\ref{eqn:sup2}) and (\ref{eqn:sup3}))
will be smaller than
a given $\delta > 0$. Taking $R = k_n$, we can then apply the Lemmas
\ref{zout},
\ref{zin} and \ref{z} with $k=1$ to the cylinders $F_n \co  [-k_n,+k_n] \times
S^1 \to \RR \times V$.
Since equations (\ref{eqn:sup1}) and (\ref{eqn:sup2}) imply that
$\|Q_0 z_n(s)\| \to 0$ as $n \rightarrow \infty$,
Lemma \ref{z} gives uniform convergence of $\|\partial^\alpha z_n(s)\|_0$ to
$0$, for $s$ satisfying $|\beta_0(s)| \le {\bar \beta}$. By the Sobolev
embedding theorem, these norm bounds imply pointwise bounds for
$|z(s,t)-z(0,t)|$. Moreover, Lemma \ref{theta} shows that
$|\vartheta_n(s,t)-\vartheta_n(0,t)| \to 0$ uniformly where
$\beta_0(s) \le {\bar \beta}$, when $n \to \infty$.
Hence, for $n$ sufficiently large, the property $\beta_0(s) \le {\bar \beta}$
will be satisfied for all $s \in [-k_n,k_n]$, so the pointwise estimates hold
for the whole $[-k_n,k_n] \times S^1$. But this contradicts the
initial assumption that $f_n(s_n,t) \notin B_\varepsilon(f_n(0,t))$ for some
$s_n \in [-k_n,k_n]$.
\end{proof}
\section{Appendix: Metric structures on the compactified moduli spaces}
\label{ap:metrics}
\subsection{Metrics  on $\overline{\Mc}^{\$}_{g,\mu}$
  and $\overline{\Mc}_{g,\mu}$ }
\label{ap:metric-surfaces}

\noindent  The space $\overline\Mc_{g,\mu}$ is the standard Deligne--Mumford
compactification
  of the moduli space $\Mc_{g,\mu}$. As it was shown by Wolpert
  (see \cite{Wolpert}), the completion of $\Mc_{g,\mu}$ in the
Weil--Petersson metric on $\Mc_{g,\mu}$  coincides
  with
  $\overline\Mc_{g,\mu}$, and thus the completed metric can be used to
metrize  the space $\overline\Mc_{g,\mu}$.
     As a topological space,  $\overline\Mc^{\$}_{g,\mu}$ can be
defined  as an oriented blow-up
  of    $\overline\Mc_{g,\mu}$ along the divisor corresponding to
nodal surfaces. Hence it also can be metrized, though
  not in a canonical way.
In this  Appendix   we   define different  metrics on
$\overline\Mc_{g,\mu}$ and  $\overline\Mc^{\$}_{g,\mu}$ compatible however
with the  same topologies on these  moduli spaces. These metrics
are  more suitable for our further considerations of moduli spaces
of holomorphic curves.
Choose an $\e\leq\frac{\e_0}2$. Given (the equivalency classes of)
two decorated nodal surfaces
$$(\bS,r)=(S,j,M,D,r)\quad\hbox{ and }
\quad(\bS',r')=(S',j',M',D',r')$$  we  take their
deformations $S^{D,r}$ and $(S')^{D',r'}$, and consider the $\e$--thick parts
  $\Tk_\e(\bS)$ and
$\Tk_\e(\bS')$  of $\bS$ and $\bS'$, viewed
as subsets of $\dot S^{D,r}$ and $(\dot S')^{D',r'}$.  Thus each
{\it compact}, ie,  not adjacent to the
punctures, component,
  $C_i$ of   the thin part of $\dot S^{D,r}$, for
$i=1,\dots, N_\e$ (resp.  $C'_i$ of the thin part of  $(\dot
S')^{D',r'}$, for $i=1,\dots, N'_\e$),
  contains   the circle
   $\Gamma_i=\Gamma_{C_i}$ (resp.
$\Gamma_i'=\Gamma_{C'_i}) $
  which is  either a closed geodesic or  one of the special circles
corresponding to
  the double points from $D$ (resp. $D'$).
Let us denote by $\Hc_\e(\bS,\bS')$
  the  (possibly empty)
  set of homeomorphisms $\dot S^{D,r}\to (\dot S')^{D',r'}$  which map
  $\Tk_\e(\bS)$ quasi-conformally onto $\Tk_\e(\bS')$.
  The homeomorphism $\varphi$  must preserve the ordering of cusps of
$\bS$ and $\bS'$ which correspond to
   the sets of marked points $M$ and $M'$, while the ordering of
circles $\Gamma_i$ and $\Gamma_i'$ is irrelevant.
For $\varphi\in\Hc_\e(\bS,\bS')$   we denote by $K_\e(\varphi)$
the maximal conformal distortion of $\varphi$ restricted to
${{\Tk_\e(\bS)}}$. If $\varphi$ is smooth,
\begin{equation}\label{eq:distortion}
K_\e(\varphi)=\max\limits_{x\in {\Tk_\e(\bS)} }
|\log\lambda_1(x)-\log\lambda_2(x)|,
\end{equation}
where ${\lambda_1(x)},{\lambda_2(x)} $  are the eigenvalues of
$(d\varphi(x))^*\circ d\varphi(x)$.
  Notice that given a homeomorphism
   $\varphi\in\Hc_\e(\bS,\bS')$ the image  $C^\varphi=\varphi(C)$
   of  a compact  component  $C$ of
  $\overline{\Tn_\e(\bS)}\subset \dot S^{D,r}$
  is a compact component of
   $\overline{\Tn_\e(\bS')}\subset (\dot S')^{D',r'}$. Let us consider
the components $\wh C$
  of $\overline{\Tn_{\e_0/2}(\bS)}\subset \dot S^{D,r}$ and and $\wh C'$ of
  $\overline{\Tn_{\e_0/2}(\bS')}\subset (\dot S')^{D',r'}$ which
contain $C$ and $C^{\varphi}$, respectively.
   Let $S_{\pm}$  $S'_{\pm}$ be the boundary circles  of $\wh C$ and $\wh C'$,
   and $\Gamma,\Gamma'$ be their central geodesics or special curves.  Let
    $\pi\co C\to\Gamma$ and $\pi'\co C'\to \Gamma'$
   be  the projections along the geodesics  orthogonal to $\Gamma$ and
$\Gamma'$.
   For any point $x\in\Gamma$ take the points
    $x_\pm\in S_\pm$ such that $\pi(x_\pm)=x$. Similarly, we define point
      $x'_\pm\in S_\pm'$ for any point $x'\in \Gamma'$.
    Let $\delta(x)$ denote the distance between the points
$\pi'(\varphi(x_+))$ and $\pi'(\varphi(x_-))$ in $\Gamma'$
    measured with
    respect to the arc length metric on $\Gamma'$.
    Similarly, $\delta'(x')$  denote the distance between the points
$\pi(\varphi^{-1}(x'_+))$ and $
    \pi(\varphi^{-1}(x'_-))$ in $\Gamma$ measured with
    respect to the arc length metric on $\Gamma$. We assume here that
the total length of  $\Gamma$ and
    $\Gamma'$ is normalized by $1$.
    Next, set
    \begin{equation}
    \label{eq:shift2}
     \delta_C(\varphi)=\mathop{\sup}\limits_{x\in\Gamma}(\delta(x))+
    \mathop{\sup}\limits_{x'\in\Gamma}(\delta(x'))\,.
    \end{equation}
  Clearly,  $\delta_C(\varphi)$ is independent of the ambiguity in choosing
   the ordering of the boundary circles of $\wh C$ and $\wh C'$.
        We  set
    \begin{equation}\label{eq:DM-metric1}
d_\e^{\varphi}\left((\bS,r),(\bS',r')\right)= K_\e(\varphi)+
\sum\limits_C\delta_C(\varphi)\;,
\end{equation}
where the sum is taken over all compact components $C$  of
$\overline\Tn_\e(\bS) \subset \dot S^{D,r}$,
           and
\begin{equation}\label{eq:DM-premetric2}
d_\e\left((\bS,r),(\bS',r')\right)=\min\big\{1,\inf\limits_{\varphi\in\Hc_\e(\bS,\bS')}
d_\e^{\varphi}\left((\bS,r),(\bS',r')\right)\big\} \;.
\end{equation}
Next, we define
the distance function $d\left((\bS,r),(\bS',r')\right)$ by the formula
\begin{equation}\label{eq:DM-metric2}
d\left((\bS,r),(\bS',r')\right)=\sum\limits_{i=1}^\infty
\frac{d_{1/2^i}\left((\bS,r),(\bS',r')\right)}{2^i}\,.
\end{equation}
\begin{proposition}\label{prop:DM-metric}
   Formula (\ref{eq:DM-metric2}) defines
a metric on $\overline\Mc^{\$}_{g,\mu}$.
\end{proposition}
\proof{} The only non-obvious properties  which we need to
  check
is the triangle inequality and the non-degeneracy of the distance
function $d$.
   
{\bf Triangle inequality}\qua Let us verify the triangle inequality
for $d_\e$. The inequality for $d$ then follows.
    Suppose that we are given three stable nodal Riemann surfaces
$(\bS_1,r_1),\;(\bS_2,r_2),\;(\bS_3,r_3)$. Denote
$\Hc_{ij}=\Hc_\e(\bS_i,\bS_j)$.
   If $\Hc_{12}=\Hc_{13}=\Hc_{23}=\varnothing$,
or if $\Hc_{12}=\Hc_{13}=\varnothing$ but $\Hc_{23}\neq\varnothing$
then the triangle inequality is obviously satisfied. If
    $\varphi_{12}\in\Hc_{12}$ and $\varphi_{23}\in\Hc_{23}$
    then $\varphi_{13}=\varphi_1\circ\varphi_2\in\Hc_{13}$.
    Notice that
    \begin{equation}
    K_\e(\varphi_{13})\leq K_\e(\varphi_{12}) +
    K_\e(\varphi_{23})
    \end{equation}
    and,
    for every component $C$ of $\Tn_\e(\bS_1)\subset S^{D_1,r_1}$ ,
     \begin{equation}
   \delta_C(\varphi_{13})
   \leq\delta_C(\varphi_{12})+
   \delta_{\varphi_{12}(C)}(\varphi_{23})\,.
   \end{equation}
    Hence,
    \begin{equation}
    \begin{split}
  d_\e^{\varphi_{13}}\left((\bS_1,r_1),(\bS_3,r_3)\right) &\leq
  K_\e(\varphi_{13})+
\sum\limits_C\delta_C(\varphi_{13})
\cr &\leq
K_\e(\varphi_{12})+
\sum\limits_C\delta_C(\varphi_{12})
  \cr
&+ K_\e(\varphi_{23}) + \sum\limits_{C'}\delta_{C'}
(\varphi_{23})
\cr
&= d_\e^{\varphi_{12}}\left((\bS_1,r_1),(\bS_2,r_2)\right)+
  d_\e^{\varphi_{23}}\left((\bS_2,r_2),(\bS_3,r_3)\right)
\end{split}
\end{equation}
where the first two sums are taken over all compact components
$$C\subset\Tn_\e(\bS_1) \subset\dot S^{D_1,r_1},$$ and the third sum
is
  taken over all compact components $$C'\subset\Tn_\e(\bS_2)
\subset\dot S^{D_2,r_2}.$$
Hence
  \begin{equation}
  d_\e\left((\bS_1,r_1),(\bS_3,r_3)\right)\leq
  d_\e\left((\bS_1,r_1),(\bS_2,r_2)\right)+
  d_\e\left((\bS_2,r_2),(\bS_3,r_3)\right).
  \end{equation}
{\bf Non-degeneracy}\qua Suppose  that
$$d\left((\bS,r),(\bS',r')\right)=0\,.$$
Then $d_{1/2^i}\left((\bS,r),(\bS',r')\right)=0\,$ for all $i\geq 0$.
Taking into account that the space
  of quasi-conformal homeomorphisms $\Tk_\e(\bS)\to\Tk_\e(\bS')$  of
bounded distortion
  is compact  for every  $\e>0$ and by passing to a diagonal
subsequence  we conclude that there
exists a  sequence of homeomorphisms  $\varphi_k\in
\Hc_{1/2^k}(\bS,\bS')$, for  $k\geq 1$, which converges uniformly
on the thick parts  to a quasi-conformal
  homeomorphism of zero conformal
  distortion, and hence to a biholomorphism
  $$\varphi\co  \left(\dot S^{D,r}\setminus \bigcup\Gamma_i,j\right)\to
\left((\dot S')^{D,r}\setminus \bigcup\Gamma'_i,j'
  \right)\,. $$
The removal of singularities now allows us to extend the
biholomorphism  $\varphi$
  to an equivalence between the nodal Riemann surfaces $\bS=(S,j,M,D)$
and $\bS'=(S',j',M',D')$.
  On the other hand, we also have for any special circle $\Gamma_i$
$$\delta_{C_i}(\varphi_k)\xrightarrow[k\to \infty] {}0\,,$$
  where $C_i$ is the component of $\overline{\Tn_{1/2^k}(\bS )}\subset
\dot S^{D,r}$ which contains $\Gamma_i$. But that implies that
$\varphi$ is  an equivalence of
the {\it decorated} nodal Riemann surfaces $(\bS,r)$ and $(\bS',r')$.
\endproof

The metric on $\overline\Mc_{g,\mu}$ can be defined by a  formula
similar to (\ref{eq:DM-metric1})--(\ref{eq:DM-metric2})  but
without the term
  $\sum\limits_C\delta_C(\varphi)$. In other words,
  the metric on $\overline\Mc_{g,\mu}$
  is by definition  the push-forward of the metric
  on $\overline\Mc^{\$}_{g,\mu}$ under the canonical projection
  $\overline\Mc^{\$}_{g,\mu}\to\overline\Mc_{g,\mu}$.

\subsection{Metric on $\overline{\Mc}_{g,\mu,p_-,p_+}(V)$}
\label{ap:metric-cylindrical}

\noindent  We will define a metric  on
$\overline \Mc_{g,\mu,p_-,p_+}(V)$ similar to the way it was done
  in Appendix \ref{ap:metric-surfaces} above
   for the Deligne--Mumford compactification of the space of Riemann surfaces.
   Take  (the equivalency classes of) two stable holomorphic building
      $(F,\Phi)$ and $(F',\Phi')$  of the same signature
$(g,\mu,p^+,p^-)$, and of
    height $k$ and $k'\leq k$, respectively.
   Let $$(\bS=\bS_{F},r=r_{F,\Phi})
   \quad\hbox{ and}\quad(\bS'=\bS_{F'},r'=r_{F',\Phi'})$$
    be  the underlying decorated Riemann nodal surfaces.
        Suppose first
    that the decorated Riemann surfaces
$$\bS=\left(S,M,D,r\right)\quad\hbox{ and}
    \quad \bS'=\left(S',M',D',r'\right)$$ are stable and  consider their
    deformations  $S^{D,r}$ and  $(S')^{D',r'}$, as in Section
\ref{sec:Riemann-stable} and Appendix \ref{ap:metric-surfaces}
    above.
    The surfaces $\dot S^{D,r}=S^{D,r}\setminus M$ and
     $(\dot S')^{D',r'}=(S')^{D',r'}\setminus M'$
         are endowed with the uniformizing
      hyperbolic metrics,
       which are degenerate along
       the special circles corresponding to the double
    points.
    Fix an $\e<\e_0$ and consider the $\e$--thick parts $\Tk_\e(\bS)
    \subset \dot S^{D,r} $ and
    $\Tk_\e(\bS')\subset (\dot S')^{D',r'}$.
   Let the notation $\Hc_\e(\bS,\bS')$ and  $d_\e^\varphi\left((
\bS,r),(\bS',r')\right)$
   for $\varphi\in\Hc_\e(\bS,\bS')$ have the same meaning as in
     Appendix \ref{ap:metric-surfaces}.
     Let us recall that the $V$--comp\-on\-ents  $f$ and $f'$ of the maps
$F$ and $F'$
     continuously extend to the maps $\of\co \overline S\to V$ and
$\of'\co \overline S'\to V$,
     where $\overline S$ and $\overline S'$ are oriented blow-ups of
the surfaces
     $S^{D,r}$ and $(S')^{D',r'}$ along the sets $M$ and $M'$ of
marked points, as it was described in Section
     \ref{sec:blow-up} above.
     Let us  set
   $$\tilde d_\e^\varphi\left((F,\Phi),(F',\Phi')\right)=
\rho(\of,\of'\circ \varphi)\,,$$
   where $\rho$ is the $C^0$--distance between mappings of a compact surface
   $\overline S$ into the  manifold $V$ endowed with a Riemannian metric.
  \begin{figure}[ht!]\small
\centering
\psfrag{t}{$t$}
\psfraga <-5pt, -5pt> {xi}{$x_i\in M$}
\psfrag{mic1}{$\min a\vert_{C_1}$}
\psfrag{mic2}{$\min a\vert_{C_2}$}
\psfrag{mic3}{$\min a\vert_{C_3}$}
\psfrag{mac2}{$\max a\vert_{C_2}$}
\psfrag{mac3}{$\max a\vert_{C_3}$}
\psfrag{mac4}{$\max a\vert_{C_4}$}
\psfrag{C1}{$C_1$}
\psfrag{C2}{$C_2$}
\psfrag{C3}{$C_3$}
\psfrag{C4}{$C_4$}
\includegraphics[width=.9\textwidth]{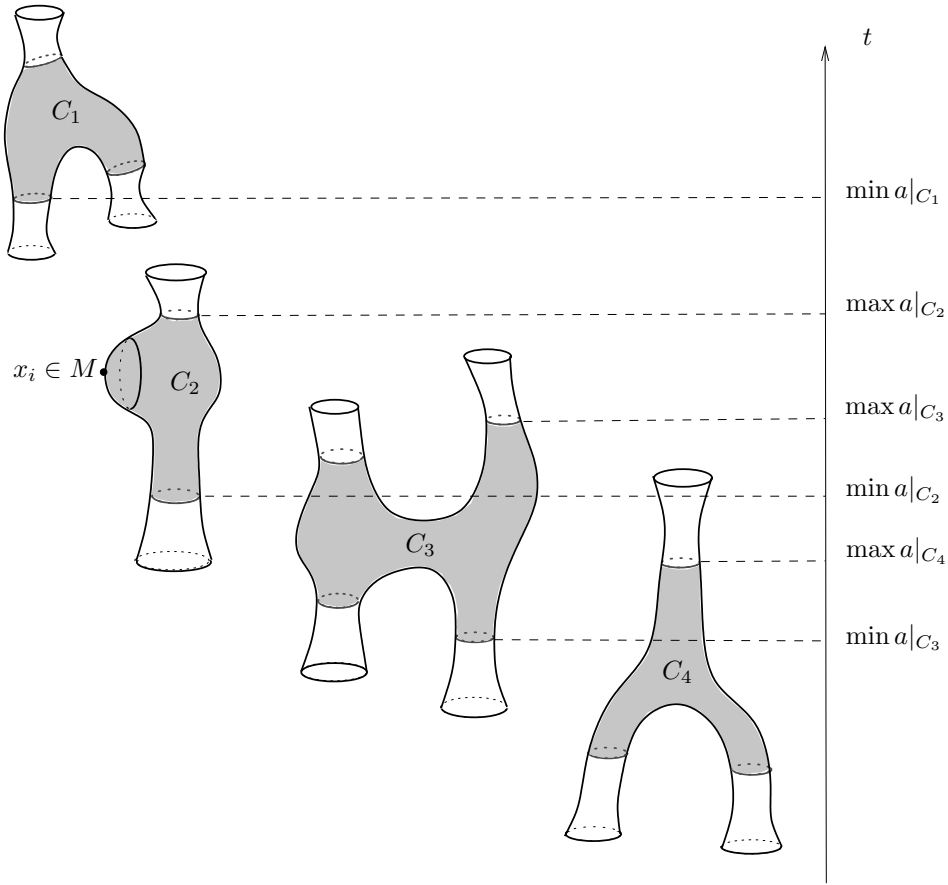}
\caption{$C_2, C_3, C_4\prec_{\varepsilon}C_1$; $C_2, C_3, C_4$ are
on the same $\varepsilon$--level}\label{fig17}
\end{figure}
   Next, we define the {\it $\e$--level} of a  connected component of
$\Tk_\e(\bS)$.
   Given two components $C,C'\subset  \Tk_\e(\bS)$, we say that  $C'$
has a bigger $\e$--level than  $C$, and write
     $C\prec_\e C'$
   if either
   \begin{itemize}
   \item{ } the level of   $C$ is less than the level of   $C'$ in the
building  $(F,\Phi)$, or
   \item{ } $C$ and $C'$ belong to the same level of the building
$(F,\Phi)$, and
     $$\min a|_{C'}> \max a|_{C},$$ and one cannot find a sequence of components
     $$C_1=C,C_2,\dots,C_k=C'$$ which belong to the same level of the
building $(F,\Phi)$ and satisfy
     $$\max a|_{C_i}\geq\max a|_{C_{i-1}}\;\;\hbox{but}\;\;
     \min a_{C_i}\leq\max a|_{C_{i-1}}\;\;\hbox{ for}\;\; i=2,\dots,k\,.$$
   \end{itemize}
   We define the {\it $\e$--level} of the component $C\subset  \Tk_\e(\bS)$
   as
   $$\max \{k\;|\hbox{there exist} \; C_1,\dots, C_{k-1}\subset \Tk_\e(\bS),\;
   C_1\prec_\e\dots\prec_\e C_{k-1}\prec_\e C\}\,.$$
    Similarly we define the $\e$--level of components
   of   $\Tk_\e(\bS')$.
   Let us denote by $C^\e_i$ ,  for $i=1,\dots,k,$  the union of components
   of $\Tk_\e(\bS)$ of  $\e$--level  $i$.
   We set   $$\hat d_\e^\varphi\left(F,F'\right)=1$$ if there exists
at least one component
   $C\subset  \Tk_\e(\bS)$ which has a different $\e$--level than the component
   $\varphi(C)\subset \Tk_\e(\bS')$. Otherwise   we set
    $$\hat d_\e^\varphi\left(F,F'\right)=\min\left\{1,
   \sum\limits_1^k\inf\limits_{c\in\R}||c+a_i-a'_i\circ
\varphi||_{C^0(C^\e_i)}\right\}\,,$$
     and define
   \begin{equation}
   \begin{split}
  &D_\e \left((F,\Phi),(F',\Phi')\right)\\
& =\min\left(1,\inf\limits_{\varphi\in\Hc_{\e}(\bS,\bS')}\big(
   d_\e^\varphi\!\left(( \bS,r),(\bS',r')\right)+\tilde
d_\e^\varphi\!\left((F,\Phi),(F',\Phi')\right)
   +\hat d_\e^\varphi\!\left(F,F'\right)\big)\right).
   \end{split}
   \end{equation}
   Next, we  introduce a distance function
$D\left((F,\Phi),(F',\Phi')\right)$ by the formula
   \begin{equation}\label{eq:dist-unstable}
   D\left((F,\Phi),(F',\Phi')\right)=\sum\limits_1^\infty
   \frac{1}{2^j}D_{1/2^j}\left((F,\Phi),(F',\Phi')\right)\,.
   \end{equation}
    Let us denote by $U_0$ the subset
    of curves from $\overline\Mc_{g,\mu,p_-,p_+}(V)$  for which the underlying
    Riemann surfaces are stable.  It is now straightforward to verify
the following proposition.
        \begin{proposition}\label{prop:dist1}
    The distance function $D$ is   a metric on $U_0$.
       \end{proposition}
    If the surfaces $\bS$ and/or  $\bS'$
       are unstable we first add extra sets $L$ and $L'$,
    $\#L=\#L'=l$, of
    marked points  to each of the holomorphic buildings
   to stabilize their underlying surfaces.
   Let us denote by $U_l$ the subset
    of buildings from $\overline\Mc_{g,\mu,p_-,p_+}(V)$ which can be
stabilized by adding
    $\leq l$ marked points,  and
       define for two curves
        $(F,\Phi),(F',\Phi')\in U_l\subset \overline\Mc_{g,\mu,p_-,p_+}(V)$
   the distance    $D^l \left((F,\Phi),(F',\Phi')\right)$
     by the formula
   \begin{equation}\label{eq:l-dist}
     D^l \left((F,\Phi),(F',\Phi')\right)=\min\big\{1,\mathop{\inf
     \limits_{L,L'}}\limits_{ \#L=\#L'=l} d(F^L, (F')^{L'})\big\}\,,
     \end{equation}
     where the infimum is taken over  all sets $L$ and $L'$
     of cardinality  $l$ which stabilize the surfaces $\bS$ and $\bS'$.
     \begin{proposition}\label{prop:dist2}
     $D^l$ is a metric on $U_l$.
     \end{proposition}
     \proof{}
     We only need to verify the non-degeneracy of $D^l$.
     Suppose that   $$D^l \left((F,\Phi),(F',\Phi')\right)\,\to\,0.$$
      Then
     there exist sequences of extra sets of marked points $L_k$ on $S$
and $L'_k$
     on $S'$, $k=1,\dots,$ such that
        $D \left((F^{L_k},\Phi),((F')^{L_k'},\Phi')\right)=0$.
        In view of the compactness  of the moduli space of stable
        nodal Riemann  surfaces (see Theorem \ref{thm:DM}) we can assume, after
        possibly passing to a subsequence, that
      $$(F^{L_j},\Phi)\to (F^{\overline L},\Phi)$$
        and    $$((F')^{L'_j},\Phi')\to ((F')^{\overline L'},\Phi').$$
        The stable Riemann surfaces  $\overline\bS$ and
$\overline\bS'$ which underly
        the holomorphic buildings $(F^{\overline L},\Phi)$  and
$((F')^{\overline L'},\Phi')$
         are limits of  sequences
        of stable Riemann surfaces $\bS_k$ and $\bS_k'$,  for $k\geq 1$,
        underlying $F^{L_k}$ and
        $(F')^{L_k'}$.
        We have
         $$D \left((F^{\overline L},\Phi),
        ((F')^{\overline L'},\Phi)\right)=0\,,$$
        and therefore  can apply Proposition \ref{prop:dist1} to
finish the proof.  \endproof

      Next, we extend the distance function $D^l$ to
      $\overline\Mc_{g,\mu,p_-,p_+}(V)$ by setting
    \begin{equation*}
     D^l \left((F,\Phi),(F',\Phi')\right)=
     \begin{cases}
     0,&\text{if both $(F,\Phi)$ and $F',\Phi')$ are not in $U_l$};\\
     1,&\text{if one of the curves is in $U_l$ and the other one}\\
     \phantom{1,}&\text{is not. }
       \end{cases}
     \end{equation*}
  Of course,   the distance function   $D^l$ extended this way to the whole
  moduli space   $\overline\Mc_{g,\mu,p_-,p_+}(V)$ is degenerate,
  and    hence is not a metric. However, it is still a {\it
  pseudo-metric}, and in particular satisfies the triangle inequality.
   Let us also note that  $D^k\leq D^l$ for  $k\geq l$.
   Finally we define   the required  metric
   $D^{\st}$ on    $\overline\Mc_{g,\mu,p_-,p_+}(V)$
   by the formula
     \begin{equation} \label{eq:dist-stab}
     D^{\st}=\sum_1^\infty\frac{D^l}{2^l}\,.
     \end{equation}
     \begin{proposition}\label{eq:dist3}
     The distance function $D^{\st}$ is a metric on
      $\overline\Mc_{g,\mu,p_-,p_+}(V)$.
      \end{proposition}
      \proof{} The only thing to check is that $D^{\st}$ is non-degenerate.
      Suppose  $D^{\st} \left((F,\Phi),(F',\Phi')\right)=0$.
      Then $D^l \left((F,\Phi),(F',\Phi')\right)=0$ for all $l\geq 0$.
      Suppose that both buildings $(F,\Phi)$ and $(F',\Phi')$ can be
stabilized by sets
     of cardinality $k$.  Then it
     follows from Proposition \ref{prop:dist2}
     that      $D^k \left((F,\Phi),(F',\Phi')\right)=0$  implies
$(F,\Phi)=(F',\Phi')$.
     \qed

\end{document}